\documentclass[11pt]{article}

\usepackage{hyperref}
\usepackage[T1]{fontenc}
\usepackage[latin1]{inputenc}
\usepackage{amsmath,amsthm,amssymb}
\usepackage{mathabx}
\usepackage{xcolor}
\usepackage{authblk}
\usepackage[margin=1in]{geometry}
\usepackage{eufrak}
\usepackage{mathrsfs}
\usepackage{enumitem} 

\providecommand{\U}[1]{\protect\rule{.1in}{.1in}}

\newtheorem{prop}{Proposition}[section]
\newtheorem{cor}[prop]{Corollary}
\newtheorem{defi}[prop]{Definition}

\newtheorem{rmk}[prop]{Remark}

\newtheorem{lem}[prop]{Lemma}

\newtheorem{theo}[prop]{Theorem}

    \newcommand{\vertiii}[1]{{\vert\kern-0.25ex\vert\kern-0.25ex\vert #1
    \vert\kern-0.25ex\vert\kern-0.25ex\vert}}

\def\tr{\mbox{\rm Tr}}

\newcommand{\EE}{\mathbb{E}}
\newcommand{\FF}{\mathbb{F}}
\newcommand{\HH}{\mathbb{H}}
\newcommand{\II}{\mathbb{I}}
\newcommand{\JJ}{\mathbb{J}}
\newcommand{\KK}{\mathbb{K}}
\newcommand{\LL}{\mathbb{L}}

\newcommand{\PP}{\mathbb{P}}

\newcommand{\RR}{\mathbb{R}}
\newcommand{\SSS}{\mathbb{S}}

\newcommand{\Ba}{ {\cal B }}
\newcommand{\Ca}{ {\cal C }}
\newcommand{\Da}{ {\cal D }}
\newcommand{\La}{ {\cal L }}

\newcommand{\Ka}{ {\cal K }}

\newcommand{\Ea}{ {\cal E }}
\newcommand{\Sa}{ {\cal S }}
\newcommand{\Ra}{ {\cal R }}

\newcommand{\Fa}{ {\cal F }}

\newcommand{\Qa}{ {\cal Q }}

\newcommand{\Ia}{ {\cal I }}

\newcommand{\Ta}{ {\cal T}}

\newcommand{\Ja}{ {\cal J }}
\newcommand{\Pa}{ {\cal P }}

\newcommand{\Ya}{ {\cal Y }}

\newcommand{\point}{\mbox{\LARGE .}}

\newcommand{\cqfd}{\hfill\blbx \\}
\def\blbx{\hbox{\vrule height 5pt width 5pt depth 0pt}\medskip}
\def \KK{\mathbb{K}}
\def \PP{\mathbb{P}}
\def \RR{\mathbb{R}}
\def \EE{\mathbb{E}}
\def \EE{\mathbb{E}}

\def \LL{\mathbb{L}}

\def \WW{\mathbb{W}}
\def  \BB{\mathbb{B}}

\numberwithin{equation}{section}

\makeatletter
\DeclareRobustCommand\frownotimes{\mathbin{\mathpalette\frown@otimes\relax}}
\newcommand{\frown@otimes}[2]{
  \vbox{
    \ialign{##\cr
      \hidewidth$\m@th#1{}_\frown$\kern-\scriptspace\hidewidth\cr
      \noalign{\nointerlineskip\kern-.1pt}
      $\m@th#1\otimes$\cr
    }
  }
}
\makeatother

\newcommand\quotient[2]{
        \mathchoice
            {
                \text{\raise0ex\hbox{$#1$}/\lower0ex\hbox{$#2$}}%
            }
            {
                #1\,/\,#2
            }
            {
                #1\,/\,#2
            }
            {
                #1\,/\,#2
            }
    }

\begin{document}

\title{A second order analysis of McKean-Vlasov semigroups}
\author[$1$]{M. Arnaudon}
\author[$2$]{P. Del Moral
}
\affil[$1$]{{\small Univ. Bordeaux, CNRS, Bordeaux INP, IMB, UMR 5251,  F-33400, Talence, France }}
\affil[$2$]{{\small INRIA, Bordeaux Research Center, Talence, France \& CMAP, Polytechnique Palaiseau, France}}
\date{}

\maketitle

\begin{abstract}
We propose a second order differential calculus to analyze the regularity and the stability properties of the distribution semigroup
associated with McKean-Vlasov diffusions. This methodology provides second order Taylor type expansions with remainder for both the evolution semigroup as well as the stochastic flow associated with this class of nonlinear diffusions. Bismut-Elworthy-Li formulae for the gradient and the Hessian of the integro-differential operators associated with these expansions are also presented. 

The article also provides explicit Dyson-Phillips expansions and a refined analysis
 of the norm of these integro-differential operators. Under some natural and easily verifiable regularity conditions we derive a series of exponential decays inequalities with respect to the time horizon. We illustrate the impact of these results with a second order extension of  the 
 Alekseev-Gr\"obner lemma to nonlinear measure valued semigroups and interacting diffusion flows. This second order 
perturbation analysis provides direct proofs of several uniform propagation of chaos properties w.r.t. the time parameter, including bias,
 fluctuation error estimate as well as exponential concentration inequalities.\\



\end{abstract}

\section{Introduction}

\subsection{Description of the models}\label{ref-desciption}

For any $n\geq 1$ we let $P_n(\RR^d)$ be the convex set of probability measures $\eta,\mu$ on $\RR^d$ with absolute $n$-th moment
and equipped  with the Wasserstein distance  of order $n$ denoted by
$\WW_n(\eta,\mu)$.
Also let  $b_t(x_1,x_2)$  be some Lipschitz function from  $\RR^{2d}$ into $\RR^d$ and 
 let  $W_t$  be an $d$-dimensional Brownian motion defined on some filtered probability space $(\Omega,(\FF_t)_{t\geq 0},\PP)$. 
 We also consider the Hilbert space
  $\HH_t(\RR^d):=\LL_2((\Omega,\FF_t,\PP),\RR^d)$ equipped with the $\LL_2$ inner product $\langle\point,\point\rangle_{\HH_t(\RR^d)}$. Up to a probability space enlargement there is no loss of generality to assume that $\HH_t(\RR^d)$ contains square integrable $\RR^d$-valued 
  variables independent of the Brownian motion. 
  
 For any $\mu\in P_2(\RR^d)$ and any time horizon $s\geq 0$ we denote by
  $X_{s,t}^{\mu}(x)$ the stochastic flow defined for any $t\in [s,\infty[$ and any starting point $x\in \RR^d$ by the McKean-Vlasov diffusion 
    \begin{equation}\label{diff-st-ref-general}
    dX^{\mu}_{s,t}(x)=b_t\left(X^{\mu}_{s,t}(x),\phi_{s,t}(\mu)\right)~dt+
    dW_t\quad \mbox{\rm with}\quad     b_t\left(x,\mu\right):=\int~\mu(dy)~b_t(x,y)
    \end{equation}
In the above display, $\phi_{s,t}$ stands for the evolution semigroup on $P_2(\RR^d)$ defined by the formulae
    $$
     \phi_{s,t}(\mu)(dy)=\mu P^{\mu}_{s,t}(dy):=\int~\mu(dx)~P^{\mu}_{s,t}(x,dy)\quad\mbox{\rm with}\quad P^{\mu}_{s,t}(x,dy):=\PP(X^{\mu}_{s,t}(x)\in dy)
    $$
      {We denote by $L_{t,\phi_{s,t}(\mu)}$ the generator of the  stochastic flow $X^{\mu}_{s,t}(x)$. }
 The existence of the stochastic flow $X^{\mu}_{s,t}(x)$ is ensured by the Lipschitz property of the drift function see for instance~\cite{graham,huang}. 
 To analyze the smoothness of the semigroup $\phi_{s,t}$ we need to strengthen this condition.

 We shall assume that  the function 
 $b_t(x_1,x_2)$ is differentiable at any order with uniformly bounded derivatives.
 In addition, the partial differential matrices
w.r.t. 
 the first and the second coordinate are uniformly bounded; that is for any $i=1,2$ we have
 \begin{equation}\label{def-H}
  \Vert b^{[i]}\Vert_2:=\sup_{t\geq 0}~\sup_{(x_1,x_2)\in\RR^{2d}}\Vert b_t^{[i]}(x_1,x_2)\Vert_2<\infty\quad\mbox{\rm with}\quad b^{[i]}_t(x_1,x_2):=\nabla_{x_i}b_t
  (x_1,x_2)
   \end{equation}  
   
  In the above display, $\Vert A\Vert_{2}:=\lambda_{\tiny max}(AA^{\prime})^{1/2}$ stands for  the spectral norm  of some matrix $A$,
  where $A^{\prime}$ stands for the transpose of $A$, $\lambda_{\tiny max}(\point)$ and  $\lambda_{\tiny min}(\point)$  the maximal and minimal eigenvalue. 
In the further development of the article, we shall also denote by $A_{\tiny sym}=(A+A^{\prime})/2$ the symmetric part of a matrix $A$.
  { In the further development of the article we represent the gradient of a real valued function as a column vector, or equivalently as the transpose of the differential-Jacobian operator which is, as any cotangent vector, represented by a row vector.
The gradient and the Hessian of a column vector valued function as tensors of type $(1,1)$ and $(2,1)$, see for instance (\ref{grad-def}). }
  
  The mean field particle interpretation of the nonlinear diffusion (\ref{diff-st-ref-general}) is described by a system of $N$-interacting diffusions
 $\xi_t=(\xi^i_t)_{1\leq i\leq N}$ defined by the stochastic differential equations
     \begin{equation}\label{diff-st-ref-general-mf}
 d\xi_t^i=b_t(\xi^i_t,m(\xi_t))~dt+ dW^i_t\quad \mbox{\rm with}\quad 1\leq i\leq N \quad \mbox{\rm and}\quad m(\xi_t):=\frac{1}{N}\sum_{1\leq j\leq N}\delta_{\xi^i_t}
  \end{equation}  
 In the above display, $\xi_0^i$ stands for $N$ independent random variables $\xi_0^i$ with common distribution $\mu_0$, and 
 $W^i_t$ are  $N$ independent copie of the Brownian motion $W_t$.

  McKean-Vlasov diffusions and their mean field type particle interpretations arise in a variety of application domains, including in porous media and granular
 flows~\cite{bene-1,bene-2,cattiaux,toscani},
 fluid mechanics~\cite{mckean-1,mckean-2,otto,villani}, data 
 assimilation~\cite{Bishop/DelMoral:2016,d-2013,DelMoral/Tugaut:2016},
 and more recently in mean field game theory~\cite{benoussan,cardaliaguet,peng-17,carmona-delarue,carmona-delarue-2,carmona-delarue-3,huang-0,gueant}, and many others.

 The origins of this subject certainly go back to the beginning of the 1950s with the article by Harris and Kahn~\cite{harris} using mean field type splitting techniques for estimating particle transmission energies. We also refer to the pioneering article by 
 Kac~\cite{kac-1,kac-2} on particle interpretations of Boltzmann and Vlasov equations, and the seminal articles by McKean~\cite{mckean-1,mckean-2} on mean field particle interpretations of nonlinear parabolic equations arising in fluid mechanics.
 Since this period, the analysis of this class of mean field type nonlinear
 diffusions and their discrete time versions have been developed in various directions. For a survey on these developments we refer 
 to~\cite{carmona-delarue,d-2013,sznitman}, and the references therein.

  The McKean-Vlasov diffusions discussed in this article belong to the class of nonlinear Markov processes.
One of the most important and difficult research questions  concerns  the regularity analysis and more particularly the stability and the long time behavior of these 
  stochastic models.  
  
   In contrast with conventional Markov processes, one of the main difficulty  of these Markov processes comes from the fact that the evolution semigroup $\phi_{s,t}(\mu)$
 is nonlinear w.r.t. the initial condition $\mu$ of the system. The additional complexity in the analysis of these models
is that their state space is the convex set of probability measures, thus conventional functional analysis and differential calculus on Banach space  cannot be directly applied.
 
The main contribution of this article is the development of a second order differential calculus to analyze the regularity and the stability properties of the distribution semigroup
associated with McKean-Vlasov diffusions. This methodology provides second order Taylor type expansions with remainder for both the evolution semigroup as well as the stochastic flow associated with this class of nonlinear diffusions. 
We also provide a refined analysis
 of the norm of these integro-differential operators with a series of exponential decays inequalities with respect to the time horizon.
 
 The article is organized as follows:
 
  The main contributions of this article are briefly discussed in section~\ref{sec-statement-intro}.
 The main theorems are stated in some detailed in section~\ref{sec-statements}. Section~\ref{spr-sec} provides some pivotal results on tensor integral operators
 and on integro-differential operators associated with the second order Taylor expansions of the semigroup $\phi_{s,t}(\mu)$. Section~\ref{tangent-sec} is dedicated to the analysis of the tangent process associated with the nonlinear diffusion flow.
We  presents explicit Dyson-Phillips expansions as well as some spectral estimates.  The last section, section~\ref{tangent-sec} is mainly concerned with the proofs  of the first and second order  Taylor expansions. The proof of some technical results are collected in the appendix.   Detailed comparisons with existing literature on this subject
 are also provided in section~\ref{comparison-sec}.   
  
    \subsection{Statement of some main results}\label{sec-statement-intro}
  One of the main contribution of the present article is the derivation of a second order Taylor expansion with remainder of the semigroup $\phi_{s,t}$ on probability spaces. For any pair of  measures $\mu_0,\mu_1\in P_2(\RR^d)$, these expansions take basically the  following form: 
\begin{equation}\label{TT-intro}
\phi_{s,t}(\mu_1)\simeq \phi_{s,t}(\mu_0)+(\mu_1-\mu_0)D_{\mu_0}\phi_{s,t}+\frac{1}{2}~(\mu_1-\mu_0)^{\otimes 2}D^2_{\mu_0}\phi_{s,t}
\end{equation}
In the above display, $D^k_{\mu_0}\phi_{s,t}$ stands some first and second order operators, with $k=1,2$.  
A more precise description of these expansions and the remainder terms is provided in section~\ref{sec-taylor-intro}. 

Section~\ref{ae-expansion},
also provides an almost sure second order Taylor expansions with remainder  of the random state $X^{\mu}_{s,t}(x)$ of the McKean diffusion 
 w.r.t. the initial distribution $\mu$.  These almost sure expansions take basically the following form
\begin{equation}\label{TT-intro-a-e}
 \begin{array}{l}
 \displaystyle X^{\mu_1}_{s,t}(x)-X^{\mu_0}_{s,t}(x)\simeq \int (\mu_1-\mu_0)(dy)~D_{\mu_0}X^{\mu_0}_{s,t}(x,y)+\frac{1}{2}~~\int~(\mu_1-\mu_0)^{\otimes 2}(dz)~D^2_{\mu_0}X^{\mu_0}_{s,t}(x,z)
\end{array}
\end{equation}
 for some random functions $D^k_{\mu_0}X^{\mu_0}_{s,t}$ from $\RR^{(1+k)d}$ into $\RR^d$, with $k=1,2$. A more precise description of these almost sure expansions is provided in section~\ref{ae-expansion} (see for instance (\ref{a-e-T1}) and theorem~\ref{theo-ae-taylor}). 
 
 Given  some random variable $Y\in \HH_s(\RR^d)$ with distribution $\mu\in P_2(\RR^d)$, observe that the stochastic flow $ \psi_{s,t}(Y):=X^{\mu}_{s,t}(Y)$
   satisfies the $\HH_t(\RR^d)$-valued stochastic differential equation
 \begin{equation}\label{eq-diff}
 d \psi_{s,t}(Y):=B_t( \psi_{s,t}(Y))~dt+dW_t
 \end{equation}
In the above display, $B_t$ stands for the drift function  from $\HH_t(\RR^d)$ into itself
defined by the formula
  $$
 B_t(X):={\EE}\left(b_t(X,\overline{X})~|~X\right)
 $$
 In the above display, $\overline{X}$ stands for an independent copy of $X$. The above Hilbert space valued representation of the 
 McKean-Vlasov diffusion (\ref{diff-st-ref-general}) readily implies that for any $Y_1,Y_0\in \HH_s(\RR^d)$ we have the exponential contraction inequality
 $$
 \Vert \psi_{s,t}(Y_1)-\psi_{s,t}(Y_0)\Vert_{\,\HH_t(\RR^d)}\leq e^{-\lambda (t-s)}~\Vert Y_1-Y_0\Vert_{\,\HH_t(\RR^d)}
 $$
   for some $\lambda>0$, as soon as the following condition is satisfied
 \begin{equation}\label{Hilbert-condition}
 \left\langle X_1-X_0, B_t(X_1)- B_t(X_0)\right\rangle_{\,\HH_t(\RR^d)}\leq-2\lambda~ \Vert X_1-X_0\Vert_{\,\HH_t(\RR^d)}^2
 \end{equation}
for any $t\geq 0$ and any $X_1,X_0\in \HH_t(\RR^d)$.
In addition,  in this framework the first order differential  $ \partial \psi_{s,t}(Y)$
of the stochastic flow coincides with the conventional Fr\'echet derivative of functions from an Hilbert space into another. 
In addition, we shall see that the gradient of first order operator $D_{\mu}\phi_{s,t}$  coincides with the dual of the tangent process  associated with the Hilbert space-valued representation (\ref{eq-diff}) of the McKean-Vlasov diffusion (\ref{diff-st-ref-general}); that is, for any smooth function $f$ we have that the dual tangent formula
 \begin{equation}\label{ref-intro-dual-D}
 \partial \psi_{s,t}(Y)^{\star}\cdot\nabla f(\psi_{s,t}(Y))=\nabla D_{\mu}\phi_{s,t}(f)(Y)
\end{equation}
A more precise description of the Fr\'echet differential  $ \partial \psi_{s,t}(Y)$ and the dual operator is provided in section~\ref{sec-H-tangent-intro} and section~\ref{tangent-sec}. A proof of the above formula is provided in theorem~\ref{theo-0}.

The Taylor expansions discussed above are valid under fairly general and easily verifiable conditions on the drift function.
  For instance, the regularity condition (\ref{def-H}) is clearly satisfied for linear drift functions. 
As it is well known, dynamical systems and hence stochastic models involving drift functions with quadratic growth require additional regularity 
conditions to 
  ensure non explosion of the solution in finite time.

  Of course the expansions (\ref{TT-intro}) and (\ref{TT-intro-a-e}) will be of rather poor practical interest without a better understanding of the 
  differential operators and the remainder terms. To get some useful approximations, we need to quantify with some precision
the norm of these operators. A important part of the article is concerned with developing a series of quantitative estimates of the differential 
operators $D^k_{\mu_0}\phi_{s,t}$ and the remainder term; see for instance theorem~\ref{theo-intro-3} and theorem~\ref{theo-intro-4}.

To avoid estimates that grow exponentially fast with respect to the time horizon,
we need to estimate with some precision the operator norms of the differential operators in (\ref{TT-intro}). To this end, we shall consider an additional regularity condition:\\
  
{\em $(H)$ :   There exists some $\lambda_0>0$ and  $\lambda_1>\Vert b^{[2]}\Vert_2$ such that for any $(x_1,x_2)\in\RR^{2d}$ and any time horizon
$t\geq 0$ we have
 \begin{equation}\label{def-HS-0}
 A_t(x_1,x_2)_{\tiny sym}\leq -\lambda_0~I\quad \mbox{and}\quad b_t^{[1]}(x_1,x_2)_{\tiny sym}\leq -\lambda_1~I
 \end{equation}
In the above display,  $I$ stands for the identity matrix and $A_t$ the matrix-valued function defined by
 \begin{equation}\label{def-lambda-1-2}
A_t(x_1,x_2):=\left[\begin{array}{cc}
b_t^{[1]}(x_1,x_2)&b_t^{[2]}(x_2,x_1)\\
b_t^{[2]}(x_1,x_2)&b_t^{[1]}(x_2,x_1)
\end{array}
\right]\quad \mbox{and we set}\quad \lambda_{1,2}:= \lambda_1-\Vert b^{[2]}\Vert_2
 \end{equation}
 }  
 
 Whenever (\ref{def-HS-0}) and (\ref{def-lambda-1-2}) are met for some parameters $\lambda_0$ and $\lambda_1\in \RR$  all the exponential estimates
stated in the article remains valid but they grow exponentially fast with respect to the time horizon.
More detailed comments on the above regularity conditions, including illustrations  for linear drift and
gradient flow models, as well as comparisons with related conditions used in the literature on this subject are also provided in section~\ref{sec-H-comments-intro}.

Under the above condition, we shall develop
several exponential decays inequalities for the norm of the differential operators $D_{\mu_0}^k\phi_{s,t}$ as well as for the  remainder terms in the Taylor expansions. The first order estimates are given in (\ref{est-nabla-D-intro}), 
the ones on the
  Bismut-Elworthy-Li gradient and Hessian extension formulae are provided in  (\ref{def-lambda-hat-bismut-intro}) and (\ref{def-lambda-hat-bismut-intro-hessian}). Second and third order estimates can also be found in (\ref{estimate-D2-nabla}) and (\ref{D3-estimation}). 

The second order differential calculus discussed above provides a natural theoretical basis to analyze the stability properties 
of the semigroup $\phi_{s,t}$ and the one of the mean field particle system discussed in (\ref{diff-st-ref-general-mf}).
 
For instance, a first order Taylor expansion of the form (\ref{TT-intro}) already indicates that the sensitivity properties of the semigroup w.r.t. the initial condition $\mu$ are encapsulated in the first order differential operator $D_{\mu}\phi_{s,t}$. Roughly speaking, whenever  $(H)$ is satisfied, we show that there exists some parameter $\lambda>0$ such that
 \begin{equation}\label{R-intro-0}
 \vee_{k=1,2}\vertiii{D_{\mu_0}^k\phi_{s,t}}\simeq e^{-\lambda (t-s)}\quad \mbox{\rm and therefore}\quad
\vertiii{\phi_{s,t}(\mu_1)- \phi_{s,t}(\mu_0)}\simeq e^{-\lambda (t-s)}
 \end{equation} 
 for some operator norms $\vertiii{\point}$.
For a more precise statement we refer to theorem~\ref{theo-intro-2} and the discussion following the theorem. 

The second order expansion (\ref{TT-intro}) also provides a natural basis to quantify the propagation of chaos properties of the mean field particle model (\ref{diff-st-ref-general-mf}). Combining these Taylor expansions with a backward semigroup analysis  we derive a  
a variety of uniform mean error estimates w.r.t. the time horizon. 
 This backward second order analysis can be seen a second order extension of the Alekseev-Gr\"obner lemma~\cite{alekseev,grobner} to nonlinear measure valued and stochastic semigroups. For a more precise statement we refer to theorem~\ref{theo-ag-mf}.  
 As in (\ref{R-intro-0}), one of the  main feature of the expansion (\ref{TT-intro}) is that it allows to enter the stability properties of the limiting semigroup $\phi_{s,t}$ into the analysis of the flow of empirical measures $m(\xi_t)$. 
 
Roughly speaking, this backward perturbation analysis can be interpreted as a second order variation-of-constants technique applied to nonlinear equations in distribution spaces.  
As in the Ito's lemma, the second order term is essential to capture the quadratic variation of the processes, see for instance the recent articles~\cite{dm-singh,hudde}
in the context of conventional stochastic differential equation, as well as in~\cite{mp-18,dm-2003} in the context of interacting jump models.

The discrete time version of this backward perturbation semigroup methodology can also be found in  chapter 7 in~\cite{d-2004}, a well as in the articles~\cite{dm-g-99,guionnet,dm-2000} and~\cite{dmrio-09,dm-hu-wu} for general classes of mean field particle systems.  

The central idea is to consider the telescoping sum on some time mesh $t_n\leq t_{n+1}$ given by the  interpolating formula
$$
m_{t_n}-\phi_{t_0,t_n}(m_{t_0})=\sum_{1\leq k\leq n}\left[\phi_{t_k,t_n}(m_{t_k})-\phi_{t_{k},t_n}\left(\phi_{t_{k-1},t_k}(m_{t_{k-1}})\right)\right]
\quad\mbox{\rm with}\quad
m_{t_k}:=m(\xi_{t_k})
$$
Applying (\ref{TT-intro}) and whenever $(t_{k}-t_{k-1})\simeq 0$ we have the second order approximation
$$
m_{t_n}-\phi_{t_0,t_n}(m_{t_0})\simeq\frac{1}{\sqrt{N}}\sum_{1\leq k\leq n} \Delta M_{t_{k}} D_{m_{t_{k-1}}}\phi_{t_{k},t_n}+\frac{1}{2N}\sum_{1\leq k\leq n}
(\Delta M_{t_{k}})^{\otimes 2}D^2_{m_{t_{k-1}}}\phi_{t_k,t}
$$
with the local fluctuation random fields 
$$
\Delta M_{t_{k}}:=\sqrt{N}~\left(m_{t_k}-\overline{m}_{t_{k}}\right)\quad \mbox{\rm and}\quad
\overline{m}_{t_{k}}:=\phi_{t_{k-1},t_k}\left(m_{t_{k-1}}\right)\simeq m_{t_{k-1}}
$$
For discrete generation particle systems, $\xi_{t_k}^i$ are defined by $N$ conditionally  independent variables given the system $\xi_{t_{k-1}}$. For a more rigorous analysis we refer to section~\ref{ips-sec}.

The above decomposition shows that
  the first order operator $D_{\mu}\phi_{s,t}$ reflects the fluctuation errors of the particle measures, while the second order term 
 encapsulates their bias. In other words, estimating the norm of second order operator $D^2_{\mu}\phi_{s,t}$ allows to quantify the bias induced by the interaction function, while the estimation of first order term is used to derive central limit theorems as well as $\LL_p$-mean error estimates.  
 
As in (\ref{R-intro-0}), these  estimates  take basically the following form.  For $n\geq 1$ and any sufficiently regular function $f$ we have
\begin{equation}\label{ln-est}
\vertiii{D_{\mu_0}\phi_{s,t}}\simeq e^{-\lambda (t-s)}\Longrightarrow 
\vert\EE\left[\Vert m_t(f)-\phi_{0,t}(m_0)(f)\Vert^n\right]^{1/n}\vert\leq c_n/\sqrt{N}
\end{equation}
In addition, we have the uniform bias estimate w.r.t. the time horizon
\begin{equation}\label{bias-intro-est}
\vertiii{D^2_{\mu_0}\phi_{s,t}}\simeq e^{-\lambda (t-s)}\Longrightarrow
\vert\EE\left[m_t(f)-\phi_{0,t}(m_0)(f)\right]\vert\leq c/N
\end{equation}
In the above display, $\vertiii{\point}$ stands   for some operator norm, and $(c,c_n)$ stands for some finite constants whose values doesn't depend on the time horizon. We emphasize that the above results are direct consequence of a second order extension of the Alekseev-Gr\"obner type lemma for particle density profiles.
For more precise statements we refer to theorem~\ref{theo-ag-mf} and the discussion following the theorem.

 \subsection{Some basic notation}\label{sec-notation}

   Let $\mbox{\rm Lin}(\Ba_1,\Ba_2)$ be the set of bounded linear operators from a normed space $\Ba_1$ into a possibly different normed space $\Ba_2$ equipped with the operator norm $\vertiii{\point}_{\Ba_1\rightarrow \Ba_2}$. When $\Ba_1=\Ba_2$  we write $\mbox{\rm Lin}(\Ba_1)$
 instead of $\mbox{\rm Lin}(\Ba_1,\Ba_1)$.

 With a slight abuse of notation, we denote by $I$ the identity $(d\times d)$-matrix, for any $d\geq 1$, as well as the identity operator  in $\mbox{\rm Lin}(\Ba_1,\Ba_1)$. 
 We also denote by $\Vert\point\Vert$ any (equivalent) norm on some finite dimensional vector space over $\RR$. 
 
 We also use the conventional notation $\partial_{\epsilon}$, $\partial_{x_i}$, $\partial_s$, $\partial_t$ and so on  for the partial derivatives 
  w.r.t. some real valued parameters $\epsilon$, $x_i$, $s$ and $t$.
 
We let $\nabla f(x)=\left[\partial_{x_i}f(x)\right]_{1\leq i\leq d}$ be the gradient column vector associated with some smooth function $f(x)$ from $\RR^d$ into $\RR$.
 Given some smooth function $h(x)$ from $\RR^d$ into $\RR^d$
we denote by $\nabla h=\left[\nabla h^1,\ldots,\nabla h^d\right]$ the gradient matrix associated with the column vector
 function $h=(h^i)_{1\leq i\leq d}$. We also let $(\nabla\otimes\nabla)$ be the second order differential operator
 defined for any twice differentiable function $g(x_1,x_2)$ on $\RR^{2d}$ by the Hessian-type formula
\begin{eqnarray}
\left((\nabla\otimes\nabla)g\right)_{i,j}&=&(\nabla_{x_1}\otimes\nabla_{x_2})(g)_{i,j}=(\nabla_{x_2}\otimes\nabla_{x_1})(g)_{j,i}=\partial_{x^i_1}\partial_{x^j_2}g\label{def-nabla-nabla}
\end{eqnarray}
 
 We consider the  space $\Ca^n(\RR^d)$ of $n$-differentiable functions and we denote by $\Ca^n_m(\RR^d)$ the subspace of functions $f$ such that 
$$
\sup_{0\leq k\leq n}{\Vert\nabla^k f(x)\Vert}\leq c~w_m(x)\quad\mbox{\rm with the weight function}\quad w_m(x)=(1+\Vert x\Vert)^{m}\quad \mbox{\rm for some}\quad m\geq 0.
$$
 We equip $\Ca^n_m(\RR^d)$ with the norm
$$
\Vert f\Vert_{\Ca^n_m(\RR^d)}:=\sum_{0\leq k\leq n}\Vert \nabla^kf/w_{m}\Vert_{\infty}
\quad\mbox{\rm with}\quad \Vert \nabla^kf/w_{m}\Vert_{\infty}=\sup_{x\in \RR^d}\Vert \nabla^kf(x)/w_{m}(x)\Vert
$$
When there are no confusions, we drop to lower symbol $\Vert\point\Vert_{\infty}$ and we write $\Vert f\Vert$ instead of $\Vert f\Vert_{\infty}$
the supremum norm of some real valued function.
We let $e(x):=x$ be the identify function on $\RR^d$ and for any $\mu\in P_n(\RR^d)$ and $n\geq 1$ we set
$$
\Vert e\Vert_{\mu,n}:=\left[\int~\Vert x\Vert^n~\mu(dx)\right]^{1/n}
$$

For any $\mu_1,\mu_2\in P_n(\RR^d)$, we also denote by $\rho_n(\mu_1,\mu_2)$ some polynomial function of $\Vert e\Vert_{\mu_i,n}$ with $i=1,2$. When $\mu_1=\mu_2$ we write $\rho_n(\mu_1)$ instead of $\rho_n(\mu_1,\mu_1)$.

Under our regularity conditions on the drift function, using elementary stochastic calculus for any $n\geq 2$ and $\mu\in P_n(\RR^d)$  we check the following estimates
   \begin{equation}\label{ref-moments}
\EE\left(\Vert X_{s,t}^{\mu}(x)\Vert^n\right)^{1/n}\leq  c_{n}(t)~\left(\Vert x\Vert+\Vert e\Vert_{\mu,2}\right)
~~ \mbox{\rm which implies that}~~
\phi_{s,t}(\mu)(\Vert e\Vert^n)^{1/n}~\leq c_{n}(t)~\Vert e\Vert_{\mu,n}
 \end{equation} 
 In the above display and throughout the rest of the article, we write $c(t), c_{\epsilon}(t), c_n(t), c_{n,\epsilon}(t), c_{\epsilon,n}(t)$ and $c_{m,n}(t)$ with $m,n\geq 0$  and $\epsilon\in [0,1]$
some collection of non decreasing  and non negative functions of  the time parameter $t$ whose values may vary from line to line, 
but which only depend on the parameters $m,n, \epsilon$,
as well as on the drift function $b_t$.   Importantly these contants do not depend on the probability measures $\mu$.
We also write $c,c_{\epsilon},c_n,c_{n,\epsilon},$ and $c_{m,n}$ when the constant do not depend on the time horizon.

 \section{Statement of the main theorems}\label{sec-statements}

\subsection{First variational equation on Hilbert spaces}\label{sec-H-tangent-intro}
 As expected, the Fr\'echet differential $\partial \psi_{s,t}(Y)$ of the stochastic flow $\psi_{s,t}(Y)$ associated with the stochastic differential equation (\ref{eq-diff})
 satisfies an Hilbert space-valued linear equation (cf. (\ref{eq-diff-lin-eq})). The drift-matrix of this evolution equation
 is given  by the Fr\'echet differential $\partial B_{t}(\psi_{s,t}(Y))$
 of the drift function $B_t$ evaluated along the solution of the flow.
Mimicking the exponential notation of the solution of conventional homogeneous linear systems, the evolution semigroup 
(a.k.a. propagator) associated with the first variational equation is written as follows 
$$
\partial \psi_{s,t}(Y)=e^{\oint_s^t\partial B_{u}(\psi_{s,u}(Y))\,du}\in  \mbox{\rm Lin}(\HH_s(\RR^d),\HH_t(\RR^d))
$$
 The above exponential is understood as an operator valued Peano-Baker series~\cite{Peano}. A more detailed presentation of these models is provided in section~\ref{tangent-sec}. 
 
 The $\HH_t(\RR^d)$-log-norm of an operator $T_t\in  \mbox{\rm Lin}(\HH_t(\RR^d),\HH_t(\RR^d))$ is defined by 
 $$
\gamma(T_t):= \sup_{\Vert Z\Vert_{\HH_t(\RR^d)}=1}{\langle Z, (T_t+T^{\star}_t)/2\cdot Z\rangle_{\,\HH_t(\RR^d)}}
 $$
 
 Our first main result is an extension of an inequality of Coppel~\cite{coppel} to tangent processes associated with Hilbert-space valued stochastic flows.
  \begin{theo}\label{theo-intro-1}
For any time horizon $t\geq s$ and any $Y\in\HH_s(\RR^d)$ we have the log-norm estimate
 \begin{equation}\label{theo-intro-1-eq-log-norm}
- \int_s^t\gamma\left(-\partial B_{u}(\psi_{s,u}(Y))\right)~du\leq  \frac{1}{t}\log{
 \vertiii{ e^{\oint_s^t\partial B_{u}(\psi_{s,u}(Y))\,du}}_{\HH_t(\RR^d)\rightarrow \HH_t(\RR^d)}}\leq \int_s^t\gamma\left(\partial B_{u}(\psi_{s,u}(Y))\right)~du
  \end{equation}
In addition, we have
 \begin{equation}\label{theo-intro-1-eq}
 (H)\Longrightarrow
  \partial B_t(X)_{\tiny sym}\leq -\lambda_0~I\Longrightarrow\frac{1}{t}\log{
 \vertiii{ e^{\oint_s^t\partial B_{u}(\psi_{s,u}(Y))\,du}}_{\HH_t(\RR^d)\rightarrow \HH_t(\RR^d)}}\leq -\lambda_0
 \end{equation}
 \end{theo}
 The proof of the above theorem in provided in section~\ref{spectral-sec}.

Let $Y_0,Y_1\in \HH_s(\RR^d)$ be a pair of random variables with distributions $(\mu_0,\mu_1)\in P_2(\RR^d)^2$.
Also let $\mu_{\epsilon}$ be the probability distribution of the random variable
\begin{equation}\label{def-Y-epsilon}
Y_{\epsilon}:=(1-\epsilon)~Y_0+\epsilon ~Y_1\Longrightarrow\partial_{\epsilon}\psi_{s,t}(Y_{\epsilon})=e^{\oint_s^t\partial B_{u}(\psi_{s,u}(Y))\,du}\cdot(Y_1-Y_0)
\end{equation}
This observation combined with the above theorem yields an alternative and more direct proof of an exponential Wasserstein contraction estimate obtained in~\cite{mp-var-18}. Namely, using (\ref{theo-intro-1-eq})
we readily check the $\WW_2$-exponential contraction inequality
\begin{equation}\label{ref-stab-W-1}
{  \partial B_t(X)_{\tiny sym}\leq -\lambda_0~I\quad \Longrightarrow}\quad\WW_2(\phi_{s,t}(\mu_1),\phi_{s,t}(\mu_0))\leq e^{-\lambda_0 (t-s)}~\WW_2(\mu_0,\mu_1)
\end{equation}
For any function $f\in \Ca^1(\RR^d)$  with bounded derivative we also quote the first order  expansion
$$
\left[\phi_{s,t}(\mu_1)-\phi_{s,t}(\mu_0)\right](f)=\int_0^1~\langle \partial\psi_{s,t}(Y_{\epsilon})^{\star}\cdot\nabla f(\psi_{s,t}(Y_{\epsilon})), (Y_1-Y_0)\rangle_{\,\HH_t(\RR^d)}~d\epsilon
$$
In the above display, $\langle \point,\point\rangle_{\,\HH_t(\RR^d)}$ stands for the conventional inner product on $\LL_2((\Omega,\FF_t,\PP),\RR^d)$.
The above assertion is a direct consequence of theorem~\ref{theo-0}.
 \subsection{Taylor expansions with remainder}\label{sec-taylor-intro}
 
 The first expansion presented in this section is a first order linearization of the measure valued mapping $\phi_{s,t}$  in terms of a semigroup of linear integro-differential
 operators.
  \begin{theo}\label{theo-intro-2}
  For any  $m,n\geq 1$ and 
 $\mu_0,\mu_1\in P_{m\vee 2}(\RR^d)$,
  there exists a semigroup of linear operators $D_{\mu_1,\mu_0}\phi_{s,t}$ from $\Ca^n_m(\RR^d)$ into itself such that
  \begin{equation}\label{a-1}
\phi_{s,t}(\mu_1)=  \phi_{s,t}(\mu_0)+(\mu_1-\mu_0)D_{\mu_1,\mu_0}\phi_{s,t}
\end{equation}
In addition, when $(H)$ is satisfied we have the gradient estimate
  \begin{equation}\label{est-nabla-D-intro}
\Vert \nabla D_{\mu_1,\mu_0}\phi_{s,t}(f)\Vert\leq c~e^{-\lambda(t-s)}~\Vert\nabla f\Vert\quad \mbox{for some $\lambda>0$}
\end{equation}
 \end{theo}
 
  The proof of the above theorem with a more explicit description of the first order operators $D_{\mu_1,\mu_0}\phi_{s,t}$ are provided in section~\ref{g-sg-sec}. {In (\ref{est-nabla-D-intro}) we can choose $\lambda=\lambda_{1,2}$, with the parameter $\lambda_{1,2}$ introduced in (\ref{def-lambda-1-2})}.  The semigroup property is a consequence of theorem~\ref{theo-1} and the gradient estimates is a reformulation of the operator norm estimate discussed in (\ref{commutation-D-etimate}).
 
We also provide Bismut-Elworthy-Li-type formulae that allow to extend the gradient and Hessian operators $\nabla^k D_{\mu_1,\mu_0}\phi_{s,t}$ with $k=1,2$  to 
measurable and bounded functions.  When the condition $(H)$ is satisfied we show the following exponential estimates
 \begin{equation}\label{def-lambda-hat-bismut-intro}
\Vert  \nabla D_{\mu_1,\mu_0}\phi_{s,t}(f)\Vert\leq c~\left(1\vee1/\sqrt{t-s}\right)~e^{- \lambda(t-s)}~\Vert f\Vert\quad \mbox{for some $\lambda>0$}
\end{equation}
In addition, we have the Hessian estimate
  \begin{equation}\label{def-lambda-hat-bismut-intro-hessian}
\Vert  \nabla^2 D_{\mu_1,\mu_0}\phi_{s,t}(f)\Vert\leq c~\left(1\vee 1/(t-s)\right)~e^{- \lambda(t-s)}~\Vert f\Vert
\quad \mbox{for some $\lambda>0$}
\end{equation}
  The proof of the first assertion can be found in  remark~\ref{rmk-1} on page~\pageref{rmk-1}. The proof of the Hessian estimates 
  is a consequence of the  decomposition of $\nabla^2 D_{\mu_0,\mu_1}\phi_{s,t}$ discussed in (\ref{dD-W}) and the Hessian estimates (\ref{bismut-est-P2}) and (\ref{def-lambda-hat-bismut}).
 
It is worth mentioning that  the semigroup property is equivalent to the chain rule formula
\begin{equation}\label{sg-chain-rule}
 D_{\mu_1,\mu_0}\phi_{s,t}=D_{\mu_1,\mu_0}\phi_{s,u}\circ D_{\phi_{s,u}(\mu_1),\phi_{s,u}(\mu_0)}\phi_{u,t}
\end{equation}
which is valid for any $s\leq u\leq t$. Without further work,
theorem~\ref{theo-intro-2} also yields the exponential $\WW_1$-contraction inequality
\begin{equation}\label{ref-stab-W-1-1}
\WW_1(\phi_{s,t}(\mu_1),\phi_{s,t}(\mu_0))\leq c~e^{-\lambda (t-s)}~\WW_1(\mu_0,\mu_1)
\end{equation}
with the same parameter $\lambda$ a in (\ref{est-nabla-D-intro}).
In the same vein, the estimate (\ref{def-lambda-hat-bismut-intro}) yields the total variation estimate
$$
\Vert \phi_{s,t}(\mu_1)-\phi_{s,t}(\mu_0)\Vert_{\tiny tv}\leq c~\left(1\vee1/\sqrt{t-s}\right)~e^{- \lambda(t-s)}~\Vert\mu_0-\mu_1\Vert_{\tiny tv}
$$
with the same parameter $\lambda$ a in (\ref{def-lambda-hat-bismut-intro}). {In all the inequalities discussed above we can choose any 
parameter $\lambda>0$ such that $\lambda<\lambda_{1,2}$, with the parameter $\lambda_{1,2}$ introduced in (\ref{def-lambda-1-2}). In the $\WW_1$-contraction inequality (\ref{ref-stab-W-1-1}) we can choose $\lambda=\lambda_{1,2}$. A more refined estimate is provided in section~\ref{sec-H-comments-intro}.}

Next theorem provides a first order Taylor expansion with remainder.

  \begin{theo}\label{theo-intro-3}
  For any  $m,n\geq 0$ and 
 $\mu_0,\mu_1\in P_{m+2}(\RR^d)$,
  there exists a  linear operators $D^2_{\mu_1,\mu_0}\phi_{s,t}$ from $\Ca^{n+2}_m(\RR^d)$ into $\Ca^{n}_{m+2}(\RR^{2d})$ such that
  \begin{equation}\label{s-o-r}
\phi_{s,t}(\mu_1)=  \phi_{s,t}(\mu_0)+(\mu_1-\mu_0)D_{\mu_0}\phi_{s,t}+\frac{1}{2}~(\mu_1-\mu_0)^{\otimes 2}D^2_{\mu_1,\mu_0}\phi_{s,t}
\end{equation}
with the first order operator $ D_{\mu_0}\phi_{s,t}:=D_{\mu_0,\mu_0}\phi_{s,t}$ introduced in theorem~\ref{theo-intro-2}. 
In addition, when $(H)$ is satisfied we also have the estimate
\begin{equation}\label{estimate-D2-nabla}
\Vert (\nabla\otimes\nabla) D^2_{\mu_1,\mu_0}\phi_{s,t}(f)\Vert\leq c~e^{-\lambda(t-s)}~\sup_{i=1,2}\Vert \nabla^i f\Vert\quad \mbox{for some $\lambda>0$}
\end{equation}

 \end{theo}
 { The proof of the above theorem in provided in section~\ref{sec-first-order-remainder}. A more precise description of the second order operator $D^2_{\mu_1,\mu_0}\phi_{s,t}$ is provided in (\ref{def-D2}) and (\ref{nabla-2-D-2}). Using (\ref{s-o-r})  and arguing as in the proof of proposition 2.1 in~\cite{mp-18}, for any twice differentiable function $f$ with bounded derivatives we check the backward evolution equation
 \begin{equation}\label{ref-backward-eq}
\partial_s\phi_{s,t}(\mu)(f)=-\mu L_{s,\mu}\left(D_{\mu}\phi_{s,t}(f)\right)
\end{equation}
with the first order operator $ D_{\mu}\phi_{s,t}$ introduced in theorem~\ref{theo-intro-3}. 
The above equation is a central tool to derive an extended version of the Alekseev-Gr\"obner lemma~\cite{alekseev,grobner} to measure valued semigroups and interacting diffusions (cf. theorem~\ref{theo-ag-mf}).}

Next theorem provides a second order Taylor expansion with remainder.

  \begin{theo}\label{theo-intro-4}
  For any  $m,n\geq 1$ and 
 $\mu_0,\mu_1\in P_{m+4}(\RR^d)$,
  there exists a  linear operators $D^3_{\mu_1,\mu_0}\phi_{s,t}$ from $\Ca^{n+3}_m(\RR^d)$ into $\Ca^{n}_{m+4}(\RR^{3d})$ such that
\begin{equation}\label{taylor-2-intro}
\begin{array}{l}
\phi_{s,t}(\mu_1)- \phi_{s,t}(\mu_0)\\
\\
\displaystyle=(\mu_1-\mu_0)D_{\mu_0}\phi_{s,t}+\frac{1}{2}
(\mu_1-\mu_0)^{\otimes 2}D^2_{\mu_0}\phi_{s,t}+(\mu_1-\mu_0)^{\otimes 3}D^3_{\mu_0,\mu_1}\phi_{s,t}
\end{array}
\end{equation}
with the second order operator $ D^2_{\mu_0}\phi_{s,t}:=D^2_{\mu_0,\mu_0}\phi_{s,t}$ introduced in theorem~\ref{theo-intro-3}. In addition, 
when $(H)$ is satisfied we have the third order estimate
\begin{equation}\label{D3-estimation}
\begin{array}{l}
\vert (\mu_1-\mu_0)^{\otimes 3}D^3_{\mu_0,\mu_1}\phi_{s,t}(f)\vert\\
\\
\leq c~ e^{-\lambda(t-s)}~\left(\vee_{i=1,2,3}\Vert\nabla^i f\Vert\right)~\WW_2(\mu_0,\mu_1)^3\quad \mbox{for some $\lambda>0$}
\end{array}
\end{equation}

 \end{theo}
  The proof of the first part of the above theorem in provided in section~\ref{sec-o-analysis}. {We can choose in (\ref{D3-estimation}) any 
parameter $\lambda>0$ such that $\lambda<\lambda_{1,2}$, with the parameter $\lambda_{1,2}$ introduced in (\ref{def-lambda-1-2})}.
The proof of the third order estimate (\ref{D3-estimation}) is rather technical, thus it is provided in the appendix, on page~\pageref{D3-estimation-lem-proof}.

 \subsection{Illustrations}\label{illustration-sec}

The first part of this section states with more details the almost sure expansions discussed in (\ref{TT-intro-a-e}). Up to some differential calculus technicalities, this result is a more or less direct consequence of the Taylor expansions with remainder presented in theorem~\ref{theo-intro-3} and theorem~\ref{theo-intro-4} combining with a backward formula presented in~\cite{mp-var-18}.
 
 The second part of this section is concerned with a second order extension of  the 
 Alekseev-Gr\"obner lemma to nonlinear measure valued semigroups and interacting diffusion flows. This second order stochastic perturbation analysis 
 is also mainly based on the second order Taylor expansion with remainder presented in theorem~\ref{theo-intro-4} .
 
  In the further development of this section without further mention we shall assume that condition $(H)$ is satisfied.

  \subsubsection{Almost sure expansions}\label{ae-expansion}

We recall the backward formula
\begin{equation}\label{stoch-interpolation-bis}
 X^{\mu_1}_{s,t}(x)-X^{\mu_0}_{s,t}(x)=\int_s^t\left[\nabla X^{\phi_{s,u}(\mu_0)}_{u,t}\right]({X}^{\mu_1}_{s,u}(x))^{\prime}~\left[\phi_{s,u}(\mu_1)-\phi_{s,u}(\mu_0)\right](b_u(X^{\mu_1}_{s,u}(x),\point))
~du
\end{equation}
  The above formula combined with (\ref{ref-stab-W-1}) and the tangent process estimates presented in section~\ref{var-eq-sec} yields  the  uniform almost sure estimates
\begin{equation}\label{estimate-X-mu}
\Vert  X^{\mu_1}_{s,t}(x)-X^{\mu_0}_{s,t}(x)\Vert\leq e^{-(\lambda_0\wedge\lambda_1)(t-s)}~\WW_2(\mu_0,\mu_1)
\end{equation}
The above estimate is a consequence of (\ref{ref-stab-W-1}) and conventional exponential estimates of the tangent process  $\nabla X^{\mu}_{s,t}$ (cf. for instance (\ref{tau-1-estimates})). A detailed proof of this claim and the backward formula (\ref{stoch-interpolation-bis}) can be found in~\cite{mp-var-18}.

  We extend the operators $D^k_{\mu}\phi_{s,t}$ introduced in theorem~\ref{theo-intro-4} to tensor valued functions $f=(f_{i})_{i\in [n]}$ with $i=(i_1,\ldots,i_n)\in [n]:=\{1,\ldots,d\}^n$ by considering the same type tensor function
with entries
\begin{equation}\label{def-D-k-extension-intro}
D_{\mu}^k\phi_{s,t}(f)_i:=D_{\mu}^k\phi_{s,t}(f_{i})\quad\mbox{\rm and we set}\quad d_{s,t}^{\mu}(x,y):=D_{\mu}\phi_{s,t}(b_t(x,\point))(y)
\end{equation}
for any $(x,y)\in \RR^{2d}$. A brief review on tensor spaces is provided in section~\ref{not-tensor-sec}. We also consider the function
$$
D_{\mu}X^{\mu}_{s,t}(x,y):=\int_s^t\left[\nabla X^{\phi_{s,u}(\mu)}_{u,t}\right]({X}^{\mu}_{s,u}(x))^{\prime}~d_{s,u}^{\mu}(X^{\mu}_{s,u}(x),y)~du
$$
Combining the first order formulae stated in theorem~\ref{theo-intro-3} with conventional Taylor expansions we check the following theorem.
\begin{theo}
For any $x\in\RR^d$, $\mu_0,\mu_1\in P_2(\RR^d)$ and $s\leq t$ we have the almost sure expansion
\begin{equation}\label{a-e-T1}
 \displaystyle X^{\mu_1}_{s,t}(x)-X^{\mu_0}_{s,t}(x)=\int~(\mu_1-\mu_0)(dy)~D_{\mu_0}X^{\mu_0}_{s,t}(x,y)+\Delta^{[2],\mu_0,\mu_1}_{s,t}(x)
\end{equation}
with
 the second order remainder function $\Delta^{[2],\mu_0,\mu_1}_{s,t}$ such that
$$
\Vert \Delta^{[2],\mu_0,\mu_1}_{s,t}\Vert\leq c~e^{-\lambda(t-s)}~\WW_2(\mu_0,\mu_1)^2
\quad \mbox{for some $\lambda>0$}
$$
\end{theo}
The detailed proof of the above theorem is provided in the appendix, on page~\pageref{theo-ae-taylor-proof}.

Second order expansions are expressed in terms of the functions defined for any $(x,y)\in \RR^{2d}$ and for any $z\in \RR^{2d}$ by the formulae
$$
 d_{s,t}^{[1,1],\mu}(x,y):=D_{\mu}\phi_{s,t}(b_t^{[1]}(x,\point)^{\prime})(y)\quad\mbox{\rm and}\quad d_{s,t}^{[2],\mu}(x,z):=D_{\mu}^2\phi_{s,t}(b_t(x,\point))(z)
$$
 We associate with these objects the function  $D^2_{\mu_0}X^{\mu_0}_{s,t}$ defined
 by
\begin{eqnarray*}
 \displaystyle D^2_{\mu}X^{\mu}_{s,t}(x,z)&:=&\int_s^t\left[\nabla X^{\phi_{s,u}(\mu)}_{u,t}\right]({X}^{\mu}_{s,u}(x))^{\prime}~\left[
 d_{s,u}^{[2],\mu}(X^{\mu}_{s,u}(x),z)+D^{[1,1]}_{\mu}X^{\mu}_{s,u}(x,z)\right]~du\\
&&\hskip4cm \displaystyle+\int_s^t\left[\nabla^2 X^{\phi_{s,u}(\mu)}_{u,t}\right]({X}^{\mu}_{s,u}(x))^{\prime}~D^{[2,1]}_{\mu}X^{\mu}_{s,u}(x,z)~du
\end{eqnarray*}
In the above display, $D^{[i,1]}_{\mu}X^{\mu}_{s,u}$ stands for the functions given by
$$
\begin{array}{l}
 \displaystyle D^{[1,1]}_{\mu}X^{\mu}_{s,u}(x,z):=\left[d_{s,u}^{[1,1],\mu}(X^{\mu}_{s,u}(x),z_2)~
D_{\mu}X^{\mu}_{s,u}(x,z_1)+d_{s,u}^{[1,1],\mu}(X^{\mu_0}_{s,u}(x),z_1)~
D_{\mu}X^{\mu}_{s,u}(x,z_2)\right]~\\
\\
 \displaystyle  D^{[2,1]}_{\mu_0}X^{\mu}_{s,u}(x,z)
 :=
\left[D_{\mu}X^{\mu}_{s,u}(x,z_1)~d_{s,u}^{\mu}(X^{\mu}_{s,u}(x),z_2)+D_{\mu}X^{\mu}_{s,u}(x,z_2)~d_{s,u}^{\mu}(X^{\mu}_{s,u}(x),z_1)\right]
\end{array}
$$
We are now in position to state the main result of this section.
\begin{theo}\label{theo-ae-taylor}
For any $x\in\RR^d$, $\mu_0,\mu_1\in P_2(\RR^d)$ and $s\leq t$ we have the almost sure expansion
\begin{equation}\label{ae-taylor}
\begin{array}{l}
 \displaystyle X^{\mu_1}_{s,t}(x)-X^{\mu_0}_{s,t}(x)\\
 \\
  \displaystyle=\int (\mu_1-\mu_0)(dy)~D_{\mu_0}X^{\mu_0}_{s,t}(x,y)+\frac{1}{2}~~\int~(\mu_1-\mu_0)^{\otimes 2}(dz)~D^2_{\mu_0}X^{\mu_0}_{s,t}(x,z)
+ \Delta^{[3],\mu_0,\mu_1}_{s,t}(x)
\end{array}
\end{equation}
with a third order remainder function $\Delta^{[3],\mu_1,\mu_0}_{s,t}$ such that 
$$
\Vert \Delta^{[3],\mu_0,\mu_1}_{s,t}\Vert\leq c~e^{-\lambda(t-s)}~\WW_2(\mu_0,\mu_1)^3\quad \mbox{for some $\lambda>0$}
$$
\end{theo}
The proof of the above theorem is provided in the appendix, on page~\pageref{theo-ae-taylor-proof}.  {In the remainder term estimates presented in the above theorems, we can choose any 
parameter $\lambda>0$ such that $\lambda<\lambda_{1,2}$, with the parameter $\lambda_{1,2}$ introduced in (\ref{def-lambda-1-2})}.

  \subsubsection{Interacting diffusions}\label{ips-sec}
For any $N\geq 2$,  the $N$-mean field particle interpretation associated with a collection of generators $L_{t,\eta}$ is defined by the Markov process
 $\xi_t=\left(\xi_t^i\right)_{1\leq i\leq N}\in  (\RR^d)^N$ with generators $ \Lambda_t$ given for any  sufficiently smooth function 
 $F$ and any $x=(x^i)_{1\leq i\leq N}\in  (\RR^d)^N$
 by
 \begin{equation}\label{ref-mean-field-FK}
\Lambda_t(F)(x)=\sum_{1\leq i\leq N}~L_{t,m(x)}(F_{x^{-i}})(x^i)     \end{equation}
 with the function
 $$
 F_{x^{-i}}(y):=F\left(x^1,\ldots,x^{i-1},y,x^{i+1},\ldots,x^N\right)\quad \mbox{\rm and the measure}\quad m(x)=\frac{1}{N}~\sum_{1\leq i\leq N}~\delta_{x^i}
 $$
 We extend $L_{t,\mu}$ to symmetric  functions $F(x^1,x^2)$ on $\RR^{2d}$ by setting
$$
L_{t,\mu}^{(2)}(F)(x^1,x^2):=L_{t,\mu}(F(x^1,\point))(x^2)+L_{t,\mu}(F(\point,x^2))(x^1)
$$
In this notation, in our context we readily check that
\begin{equation}\label{backward-D-phi-mf}
\begin{array}{rclcrcl}
\Fa(x)&=&m(x)(f)&\Longrightarrow&\displaystyle
\Lambda_t(\Fa)(x)&=& m(x)L_{t,m(x)}(f)\\
&& &&&&\\
\Fa(x)&=&m(x)^{\otimes 2}(F)&\Longrightarrow&\displaystyle
\Lambda_t(\Fa)(x)&=& \displaystyle m(x)^{\otimes 2}L_{t,m(x)}^{(2)}(F)+\frac{1}{N}~m(x)\left[\Gamma(F)\right]
\end{array}\end{equation}
for any symmetric function $F(x^1,x^2)=F(x^2,x^1)$, with the function $\Gamma(F)$ on $\RR^d$ defined for any $y\in\RR^d$ by
the formula
$$
\begin{array}{l}
\displaystyle\Gamma(F)(y)
:=\tr\left(\left(\left[\nabla\otimes\nabla\right]F\right)(y,y)\right)=\sum_{1\leq i\leq d}~\left(\partial_{x_1^i}\,\partial_{x_2^i}\, F\right)(y,y)\\
\\
\displaystyle\Longrightarrow
\Gamma\left(f\otimes g\right)(y)=\sum_{1\leq k\leq d}\partial_{y_k}f(y)~\partial_{y_k}g(y)=\tr\left(\nabla f(y)\nabla g(y)^{\prime}\right)
\end{array}$$
A proof of the above formula is provided in the appendix, on page~\pageref{backward-D-phi-mf-proof}.
Applying Ito's formula, for any smooth function $g:t\in [0,\infty[\mapsto g_t\in \Ca^2_b(\RR^d)$  we prove that
$$
m_t:=m(\xi_t)\Longrightarrow
dm_t(g_t)=\left[m_t\left(\partial_tg_t\right)+m_tL_{t,m_t}(g_t)\right]~dt+\frac{1}{\sqrt{N}}~dM_t(g)
$$
In the above display, $g\mapsto M_t(g)$ stands for a martingale random field with angle bracket
$$
\partial_t\langle M(f),M(g)\rangle_t:=m_t\left(\Gamma(f\otimes g)\right)\Longrightarrow
\partial_t\langle M(g)\rangle_t=\int~m_t(dx)~\Vert \nabla g(x)\Vert^2
$$
The above evolution equation is rather standard in mean field type interacting particle system theory, a detailed proof 
 can be found in~\cite{dm-geom-2013} (see for instance section 4.3). In the same vein, with some obvious abusive notation,
 using (\ref{backward-D-phi-mf}) we have
\begin{eqnarray*}
dm_s^{\otimes 2}(F)&=& \left[m_s\otimes dm_s+ dm_s\otimes m_s +(dm_s\otimes dm_s)\right](F)\\
&=& \left[m_s^{\otimes 2}L_{s,m_s}^{(2)}(F)+\frac{1}{N}~m_s\left[\Gamma(F)\right]\right]~ds+\mbox{\rm martingale increment}\\
&&\Longrightarrow
\left[dm_s\otimes dm_s\right](F)=\frac{1}{N}~m_s\left[\Gamma(F)\right]~ds
\end{eqnarray*}
 
We fix a final time horizon $t\geq 0$ and we denote by $$s\in [0,t]\mapsto M_{s}\left(D_{m_{\point}}\phi_{\point,t}(f)\right)$$  the martingale
associated with the predictable function
$$
s\in [0,t]\mapsto g_s=D_{m_s}\phi_{s,t}(f)
$$
Combining the It\^o formula with the tensor product formula (\ref{backward-D-phi-mf}) and with the backward formula 
(\ref{ref-backward-eq}) 
we obtain
$$
d\,\phi_{s,t}(m_s)(f)=-m_s L_{s,m_s}\left(D_{m_s}\phi_{s,t}(f)\right)~ds+(dm_s)\left(D_{m_s}\phi_{s,t}(f)\right)+\frac{1}{2}~(dm_s\otimes dm_s)(D^2_{m_s}\phi_{s,t}(f))~ds
$$
This implies that
$$
d\,\phi_{s,t}(m_s)(f)=\frac{1}{2}~(dm_s\otimes dm_s)(D^2_{m_s}\phi_{s,t}(f))~ds+ \frac{1}{\sqrt{N}}~dM_{s}\left(D_{m_{\point}}\phi_{\point,t}(f)\right)
$$

This yields the following theorem.

\begin{theo}\label{theo-ag-mf}
 For any time horizon $t\geq 0$, the interpolating semigroup $s\in [0,t]\mapsto\phi_{s,t}(m_s)$ satisfies for any  $f\in \Ca^2(\RR^d)$ with $\sup_{k=1,2}\Vert \nabla^kf\Vert\leq 1$ the evolution equation
\begin{equation}
d\,\phi_{s,t}(m_s)(f)=\frac{1}{2N}~m_s\left[\Gamma\left(D^2_{m_s}\phi_{s,t}(f)\right)\right]~ds+ \frac{1}{\sqrt{N}}~dM_{s}\left(D_{m_{\point}}\phi_{\point,t}(f)\right)
\end{equation}
\end{theo}

The above theorem can be seen as a second order extension of the Alekseev-Gr\"obner lemma~\cite{alekseev,grobner} to nonlinear measure valued and stochastic semigroups.  This result also  extends the perturbation theorem obtained in~\cite{mp-18} (cf. theorem 3.6) in the context of interacting jumps processes to McKean-Vlasov diffusions. The discrete time version of the backward perturbation analysis described above can also be found in~\cite{dm-g-99,guionnet,dm-2000}
in the context of Feynman-Kac particle models (see also~\cite{d-2004,d-2013,dm-2003}).

We end this section with some direct consequences of the above theorem. Firstly, using (\ref{est-nabla-D-intro}) and (\ref{estimate-D2-nabla}) we have the almost sure estimates
\begin{eqnarray*}
\vert\partial_s\langle M_{\point,t}\left(D_{m_{\point}}\phi_{\point,t}(f)\right)
\rangle_s\vert&\leq &c~e^{-2\lambda(t-s)}~\Vert\nabla f\Vert^2
\\
\mbox{\rm and}\quad\Vert m_s\left[\Gamma\left(D^2_{m_s}\phi_{s,t}(f)\right)\right]\Vert&\leq& c~e^{-\lambda(t-s)}~\sup_{i=1,2}\Vert \nabla^i f\Vert\quad \mbox{for some $\lambda>0$}
\end{eqnarray*}
Without further work, the above inequality yields the uniform bias estimate stated in the r.h.s. of  (\ref{bias-intro-est}),
 for any twice differentiable  function $f$ with bounded derivatives. Using well known martingale concentration inequalities (cf. for instance lemma 3.2 in~\cite{nishi}), there exists some finite parameter $c$ such that for any $t\geq 0$ and any $\delta\geq 1$ the probability of the following event
 $$
\vert m_t(f)- \phi_{0,t}(m_0)(f)-\frac{1}{2N}\int_0^t~m_s\left[\Gamma\left(D^2_{m_s}\phi_{s,t}(f)\right)\right]~ds\vert\leq c~\sqrt{\frac{\delta}{N}}
 $$
 is greater than $1-e^{-\delta}$.
In addition, using the  Burkholder-Davis-Gundy inequality,  for any $n\geq 1$ we obtain the time uniform estimates stated in the r.h.s. of (\ref{ln-est}).
On the other hand, using (\ref{a-1}) and (\ref{est-nabla-D-intro}) we have the almost sure exponential contraction inequality
$$
\WW_1(\phi_{0,t}(m_0),\phi_{0,t}(\mu_0))\leq c~e^{-\lambda t}~\WW_1(m_0,\mu_0)\quad \mbox{for some $\lambda>0$}
$$
This yields the bias estimates 
$$
\vert\EE\left[m_t(f)-\phi_{0,t}(\mu_0)(f)\right]\vert\leq \frac{c_1}{N}+\frac{c_2}{N^{1/d}}~e^{-\lambda t}~
$$
 for any twice differentiable  function $f$ with bounded derivatives. The r.h.s. estimate comes from well known estimates of the 
 average of the Wassertein distance for occupation measures, see for instance~\cite{dudley-69}
 and the more recent studies~\cite{fournier,lei-18}. The above inequality yields the following uniform bias estimate
 $$
 \sup_{t\geq \frac{d-1}{d\lambda}\log{N}}\vert\EE\left[m_t(f)-\phi_{0,t}(\mu_0)(f)\right]\vert\leq \frac{c}{N}
 $$

  \subsection{Comments on the regularity conditions}\label{sec-H-comments-intro}

We discuss in this section  the regularity condition $(H)$ introduced in (\ref{def-HS-0}). We illustrate these spectral conditions for linear-drift and
gradient flow models. Comparisons with related conditions presented in other works are also provided.

Firstly, we mention that the condition stated in (\ref{def-HS-0}) has been introduced in the article~\cite{mp-var-18} to derive several Wasserstein exponential contraction inequalities as well as uniform propagation of chaos estimates w.r.t. the time horizon.

Using the log-norm triangle inequality and recalling that the log-norm is dominated by the spectral norm 
 we check that
 $$
 \lambda_{\tiny max}(A_t(x_1,x_2)_{\tiny sym})\leq  \lambda_{\tiny max}(b_t^{[1]}(x_1,x_2)_{\tiny sym})+2^{-1}~\Vert b_t^{[2]}(x_2,x_1)+b_t^{[2]}(x_1,x_2)^{\prime}\Vert_2
 $$
 
Choosing $\lambda_0$ and $\lambda_1$ as the supremum of the maximal eigenvalue functional of the matrices $ A_t(x_1,x_2)_{\tiny sym}$
 and $b_t^{[1]}(x_1,x_2)_{\tiny sym}$,  the Cauchy interlacing theorem (see for instance~\cite{lancaster} on page 294) yields 
 $  \lambda_1\geq \lambda_0\geq \lambda_{1,2}
 $. 
 
 For linear drift functions 
  \begin{equation}\label{lin-case-intro}
 b_t(x_1,x_2)=B_1\,x_1+B_2\,x_2
  \end{equation}
  the matrix $A_t(x_1,x_2)_{\tiny sym}$ reduces to the two-by-two block partitioned matrix
 \begin{equation}\label{lin-case}
 A_t(x_1,x_2)_{\tiny sym}=\left[\begin{array}{cc}
 (B_1)_{\tiny sym}& (B_2)_{\tiny sym}\\
  (B_2)_{\tiny sym}& (B_1)_{\tiny sym}
 \end{array}\right]\Longrightarrow \lambda_0\geq \lambda_1=- \lambda_{\tiny max}((B_1)_{\tiny sym})\quad\mbox{\rm and}\quad
 \Vert b^{[2]}\Vert_2= \Vert B_2\Vert_2
 \end{equation}
{In this situation  the diffusion flow  $X_{s,t}^{\mu}(x)\in \RR^d$ 
is  given  by the formula
$$
 X^{\mu}_{s,t}(x)=e^{(t-s)B_1}(x-\mu(e))+e^{(t-s)[B_1+B_2]}~\mu(e)+\int_s^te^{B_1(t-u)}~dW_u
  $$
  In the one dimensional case we have
 $$
 B_1<0<B_2\quad \Longrightarrow\quad B_1=-\lambda_1\leq  B_1+B_2=-\lambda_{1,2}=-\lambda_0
 $$}
  
Nonlinear Langevin diffusions are associated with the drift function
$$\begin{array}{l}
b(x_1,x_2):=-\nabla U(x_1)-\nabla V(x_1-x_2)\\
\\

\Longrightarrow b^{[1]}(x_1,x_2)=-\nabla^2 U(x_1)-\nabla^2 V(x_1-x_2)
\quad\mbox{\rm and}\quad b^{[2]}(x_1,x_2)=\nabla^2 V(x_1-x_2)
 \end{array}$$
some confinement type potential function $U$ (a.k.a. the exterior potential) and some interaction potential function $V$.  In this context we have
$$
\begin{array}{l}
 -A_t(x_1,x_2)_{\tiny sym}=\left[\begin{array}{cc}
\nabla^2 U(x_1)& 0\\
0& \nabla^2 U(x_2)
 \end{array}\right]\\
 \\
 \hskip2cm+\left[\begin{array}{cc}
\nabla^2 V(x_1-x_2)& -(\nabla^2 V(x_2-x_1)+\nabla^2 V(x_1-x_2))/2\\
-(\nabla^2 V(x_2-x_1)+\nabla^2 V(x_1-x_2))/2& \nabla^2 V(x_2-x_1)
 \end{array}\right] \end{array}
 $$
  When the potential function $V$ is even and convex we have
  $$
   A_t(x_1,x_2)_{\tiny sym}\leq -\left[\begin{array}{cc}
\nabla^2 U(x_1)& 0\\
0& \nabla^2 U(x_2)
 \end{array}\right]
  $$
In the reverse angle, when the function $V$ is odd we have the formula
$$
 A_t(x_1,x_2)_{\tiny sym}=-\left[\begin{array}{cc}
\nabla^2 U(x_1)+\nabla^2 V(x_1-x_2)& 0\\
0& \nabla^2 U(x_2)+ \nabla^2 V(x_2-x_1)
 \end{array}\right]
$$
  In both situations, condition  $(H)$  is satisfied when the strength of the confinement type potential 
  dominates the one of the interaction potential; that is when we have that
  $$
  \nabla^2 U(x_1)+\nabla^2 V(x_2)\geq \lambda_1>\Vert\nabla^2 V\Vert_{2}
  $$
  
{  The decay rate $\lambda_0$ in the $\WW_2$-contraction inequality (\ref{ref-stab-W-1}) is larger than
  the decay rate $\lambda_{1,2}$ in the $\WW_1$-contraction inequality (\ref{ref-stab-W-1-1}).  In addition, the $\WW_1$-exponential stability requires that
  $\lambda_0$ dominates the spectral norm of the matrix $b^{[2]}$. Next we provide a more refined analysis based on  
 the proof of the $\WW_2$-contraction inequality presented in~\cite{mp-var-18}. Using the interpolating paths $(Y_{\epsilon},\mu_{\epsilon})$ introduced in (\ref{def-Y-epsilon}) we set
 \begin{equation}\label{def-mu-epsilon}
X^{\epsilon}_{s,t}:=X^{\mu_{\epsilon}}_{s,t}(Y_{\epsilon})\quad \mbox{\rm and}\quad
\overline{X}^{\epsilon}_{s,t}:=\overline{X}^{\mu_{\epsilon}}_{s,t}(\overline{Y}_{\epsilon})
 \end{equation}
 In the above display $(\overline{X}^{\mu_{\epsilon}}_{s,t}(x),\overline{Y}_{\epsilon})$ stands for an independent copy of $(X^{\mu_{\epsilon}}_{s,t}(x),Y_{\epsilon})$. Arguing as in~\cite{mp-var-18} we have
 $$
 \partial_t\EE(\Vert \partial_{\epsilon}X^{\epsilon}_{s,t}\Vert)=\EE\left[\Vert \partial_{\epsilon}X^{\epsilon}_{s,t}\Vert^{-1}\left( \langle
 \partial_{\epsilon}X^{\epsilon}_{s,t},b^{[1]}(X^{\epsilon}_{s,t},\overline{X}^{\epsilon}_{s,t}) \partial_{\epsilon}X^{\epsilon}_{s,t}\rangle+
 \langle
 \partial_{\epsilon}\overline{X}^{\epsilon}_{s,t},b^{[2]}(X^{\epsilon}_{s,t},\overline{X}^{\epsilon}_{s,t}) \partial_{\epsilon}X^{\epsilon}_{s,t}\rangle\right)\right]
 $$
 We consider the symmetric and anti-symmetric matrices
 \begin{eqnarray*}
 b^{[2]}_t(x_1,x_2)_{\tiny sym}&:=&\frac{1}{2}~ \left(b^{[2]}_t(x_1,x_2)+ b^{[2]}_t(x_2,x_1)^{\prime}\right)
 \\
  b^{[2]}_t(x_1,x_2)_{\tiny asym}&:=&\frac{1}{2}~ \left(b^{[2]}_t(x_1,x_2)- b^{[2]}_t(x_2,x_1)^{\prime}\right)
 \end{eqnarray*}
 and we set
$$
(U^{\epsilon}_{s,t},\overline{U}^{\epsilon}_{s,t}):=\left(\frac{\partial_{\epsilon}X^{\epsilon}_{s,t}}{\sqrt{\Vert \partial_{\epsilon}X^{\epsilon}_{s,t}\Vert}},\frac{\partial_{\epsilon}\overline{X}^{\epsilon}_{s,t}}{\sqrt{\Vert \partial_{\epsilon}\overline{X}^{\epsilon}_{s,t}\Vert}}\right)
\quad \mbox{\rm and}\quad
(V^{\epsilon}_{s,t},\overline{V}^{\epsilon}_{s,t}):=\left(\frac{\partial_{\epsilon}X^{\epsilon}_{s,t}}{\Vert \partial_{\epsilon}X^{\epsilon}_{s,t}\Vert},\frac{\partial_{\epsilon}\overline{X}^{\epsilon}_{s,t}}{\Vert \partial_{\epsilon}\overline{X}^{\epsilon}_{s,t}\Vert}\right)
$$  
By symmetry arguments and using some elementary manipulations we check the formula
$$
\begin{array}{l}
2\, \partial_t\,\EE(\Vert \partial_{\epsilon}X^{\epsilon}_{s,t}\Vert)
=\EE\left(\left\langle\left(
\begin{array}{c}
U^{\epsilon}_{s,t}\\
\overline{U}^{\epsilon}_{s,t}
\end{array}
\right),A_t(X^{\epsilon}_{s,t},\overline{X}^{\epsilon}_{s,t})\left(
\begin{array}{c}
U^{\epsilon}_{s,t}\\
\overline{U}^{\epsilon}_{s,t}
\end{array}
\right)\right\rangle\right.\\
\\
\hskip3cm\left.+
\left(\sqrt{\Vert \partial_{\epsilon}\overline{X}^{\epsilon}_{s,t}\Vert}-\sqrt{\Vert \partial_{\epsilon}{X}^{\epsilon}_{s,t}\Vert}\,\right)^2~
\left\langle \overline{V}^{\epsilon}_{s,t}, b^{[2]}_t(X^{\epsilon}_{s,t},\overline{X}^{\epsilon}_{s,t})_{\tiny sym}\,V^{\epsilon}_{s,t}\right\rangle\right.\\
\\
\hskip4cm\left.+\left(\Vert \partial_{\epsilon}\overline{X}^{\epsilon}_{s,t}\Vert-\Vert \partial_{\epsilon}{X}^{\epsilon}_{s,t}\Vert\,\right)~
\left\langle \overline{V}^{\epsilon}_{s,t}, b^{[2]}_t(X^{\epsilon}_{s,t},\overline{X}^{\epsilon}_{s,t})_{\tiny asym}\,V^{\epsilon}_{s,t}\right\rangle
\,\right)
\end{array}
$$
This shows that
$$
\partial_t\,\EE(\Vert \partial_{\epsilon}X^{\epsilon}_{s,t}\Vert)\leq -\widehat{\lambda}_{1,2}~\EE(\Vert \partial_{\epsilon}X^{\epsilon}_{s,t}\Vert)
$$
with the parameter  $\widehat{\lambda}_{1,2}$ given by
$$
-\widehat{\lambda}_{1,2}:=\sup_{x_1,x_2}{\left[\lambda_{\tiny max}(A_t(x_1,x_2))+\Vert  b^{[2]}_t(x_1,x_2)_{\tiny sym}\Vert_2+\Vert  b^{[2]}_t(x_1,x_2)_{\tiny asym}\Vert_2\right]}\leq -\lambda_{1,2}
$$
We conclude that  the $\WW_1$-contraction inequality (\ref{ref-stab-W-1-1}) is met with $\lambda=\widehat{\lambda}_{1,2}$.}

  In a more recent article~\cite{wang} the author presents some Wasserstein contraction inequalities of the same form as in (\ref{ref-stab-W-1}) with $\lambda_0$ replaced by 
  some parameter $\lambda^{-}_0=(\kappa_1-\kappa_2)$, under the assumption
  $$
  \langle x_1-y_1,b_t(x_1,\mu_1)-b_t(y_1,\mu_2)\rangle\leq -\kappa_1~\Vert x_1-y_1\Vert^2+\kappa_2~\WW_2(\mu_1,\mu_2)^2\quad\mbox{\rm for some}\quad \kappa_1>\kappa_2
  $$
Taking  Dirac measures $\mu_1=\delta_{x_2}$ and $\mu_2=\delta_{y_2}$ we check that the above condition is equivalent to the fact that
$$
\langle x_1-y_1,b_t(x_1,x_2)-b_t(y_1,y_2)\rangle\leq -\kappa_1~\Vert x_1-y_1\Vert^2+\kappa_2~\Vert x_2-y_2\Vert^2
$$
By symmetry arguments this implies that
 \begin{equation}\label{Hilbert-condition-sym}
\langle x_1-y_1,b_t(x_1,x_2)-b_t(y_1,y_2)\rangle+\langle x_2-y_2,b_t(x_2,x_1)-b_t(y_2,y_1)\rangle\leq -\lambda_0^-~[\Vert x_1-y_1\Vert^2+\Vert x_2-y_2\Vert^2]
 \end{equation}
For the  linear drift model discussed in (\ref{lin-case}) the above condition reads
$$
\left[\begin{array}{cc}
 (B_1)_{\tiny sym}& (B_2)_{\tiny sym}\\
  (B_2)_{\tiny sym}& (B_1)_{\tiny sym}
 \end{array}\right]\leq -\lambda_0^-~I\quad \mbox{\rm which is implies that}\quad \lambda_0\geq \lambda_0^-
$$
We also have $(\ref{Hilbert-condition-sym})\Longrightarrow (\ref{Hilbert-condition})$ with $\lambda=\lambda^-_0$.
  \subsection{Comparisons with existing literature}\label{comparison-sec}

The perturbation analysis developed in the article differs from the Otto differential calculus on $(P_2(\RR^d),\WW_2)$ introduced in~\cite{otto} and further developed by Ambrosio and his co-authors ~\cite{ambrosio-1,ambrosio-2} and Otto and Villani in~\cite{otto-2}. These sophisticated gradient flow techniques in Wasserstein metric spaces are based on optimal transport theory. 

The central idea is to interpret $P_2(\RR^d)$ as an infinite dimensional Riemannian manifold. In this context, the Benamou-Brenier formulation of the Wasserstein distance provides a natural way to define geodesics, 
 gradients and Hessians w.r.t. the Wasserstein distance. The details of these gradient flow techniques are beyond the scope of the semigroup perturbation analysis considered herein. 
 
 This methodology is mainly used
 to quantify the entropy dissipation of Langevin-type nonlinear diffusions. Thus, it cannot be used to derive any Taylor expansion of the form (\ref{TT-intro})
 nor to analyze the stability properties of more general classes of McKean-Vlasov diffusions.

Besides some interesting contact points, the methodology developed in the present article doesn't rely on the more recent differential calculus on $(P_2(\RR^d),\WW_2)$ developed by P.L. Lions and his co-authors in the seminal works on mean field game theory~\cite{cardaliaguet,gueant}. 
In this context, the  first order Lions differential of
a smooth function from $P_2(\RR^d)$ into $\RR$ is defined as the conventional derivative of lifted real valued function acting on the Hilbert space of 
square integrable random variables.  In this interpretation, for a given test function, say $f$ the gradient $\nabla D_{\mu}\phi_{s,t}(f)(Y)$ of the first order differential in (\ref{TT-intro})  can be seen as the Lions derivative $(\delta u_{s,t}/\delta\mu)(Y)$ of the lifted scalar
function $Y\mapsto u_{s,t}(Y):=\EE(f(X_{s,t}^{\mu}(Y)))$, for some random variable $Y$ with distribution $\mu$. 

In the recent book~\cite{carmona-delarue}, to distinguish these two notions, the authors called 
the random variable $D_{\mu}\phi_{s,t}(f)(Y)$ the linear functional derivative. For a more thorough discussion on the origins and the recent developments in mean field game theory, we refer to the book~\cite{carmona-delarue} as well as the more recent articles~\cite{peng-17,chassagneux,crisan-18} and the references therein. 

To the best of our knowledge, most of the literature on Lions' derivatives is concerned with existence theorems without a refined analysis of the exponential decays of these differentials w.r.t. the time parameter. Last but not least, from the practical point of view all differential estimates we found in the literature are rather quite deceiving since after carefully checking, they grow exponentially fast with respect to the time horizon (cf. for instance~\cite{peng-17,chassagneux,chaudru,crisan-18}).

Taylor 
 expansions of the form (\ref{TT-intro}) have already been discussed in the book~\cite{d-2013} for discrete time nonlinear measure valued semigroups (cf. for instance chapters 3 and 10). We also refer to the more recent article~\cite{mp-18} in the context of continuous time Feynman-Kac semigroups. 
In this context, we emphasize that the semigroup $\phi_{s,t}(\mu)$ is explicitly given by a normalization of a linear semigroup of positive operators. Thus, a fairly simple  Taylor expansion yields the second order formula (\ref{TT-intro}). In contrast with Feynman-Kac models, McKean-Vlasov semigroups don't have any explicit form nor an analytical description. As a result, none of above methodologies cannot be used to analyze nonlinear diffusions.

  The second order perturbation analysis discussed in this article has been used with success in~\cite{dm-g-99,guionnet,dm-2000} to analyze the stability properties of Feynman-Kac type particle models, as well as the fluctuations and the exponential concentration of this class of interacting jump processes; see also~\cite{dmrio-09,dm-hu-wu} for general classes of discrete generation mean field particle systems, a well as  chapter 7 in~\cite{d-2004} and~\cite{mp-18,dm-2003} for continuous time models.

  These second order perturbation techniques have also been extended  in the seminal book by V.N. Kolokoltsov~\cite{kolo-10} 
  to general classes of nonlinear Markov processes and kinetic equations. Chapter 8 in~\cite{kolo-10} is dedicated to the analysis of
the first and the second order derivatives  of nonlinear semigroups
with respect to initial data. The use of the first and the second order derivatives
in the analysis of central limit theorems and propagation of chaos properties respectively is developed 
 in Chapters 9 and Chapter 10 in~\cite{kolo-10}. We underline that these results are obtained for  diffusion processes as well as for 
  jump-type processes and their combinations, see also~\cite{kol-troeva-1,kol-troeva-2}.

  Nevertheless none of these studies apply to derive non asymptotic 
  Taylor expansions (\ref{taylor-2-intro}) and (\ref{ae-taylor}) with exponential decay-type remainder estimates
   for McKean-Vlasov diffusions nor to estimate the stability properties of the associated semigroups. In addition, to the best of our knowledge
   the stochastic perturbation theorem~\ref{theo-ag-mf} is the first result of this type for mean field type interacting diffusions.

Last but not least, the idea of considering the flow of empirical measures $m(\xi_t)$ of a mean field particle model
as a stochastic perturbation of the limiting flow $\phi_{0,t}(\mu_0)$ certainly goes back to the work by Dawson~\cite{dawson}, itself based on the martingale approach
developed by Papanicolaou, Stroock and Varadhan in~\cite{pap-77}, published in the end of the 1970's. These two works are mainly centered on fluctuation type limit
theorems. They don't discuss any Taylor expansion on the  limiting semigroup $\phi_{s,t}$ nor  any question related to the stability properties of the underlying processes.

\section{Some preliminary results}\label{spr-sec}

The first part of this section provides a review of tensor product theory and Fr\'echet differential on Hilbert spaces. 
Section~\ref{not-tensor-sec} is concerned with conventional tensor products and Fr\'echet derivatives. Section~\ref{tensor-integral-op-sec}
provides a short introduction to tensor integral operators.

In the second part of this section we review some basic tools of the theory of stochastic variational equations, including some differential properties of Markov semigroups. Section~\ref{var-eq-sec} is dedicated to variational equations. Section~\ref{bel-formula-sec} discusses Bismut-Elworthy-Li extension 
formulae. We also provide some exponential inequalities for the gradient and the Hessian operators on bounded measurable functions.

The differential operator arising in the Taylor expansions (\ref{TT-intro}) are defined in terms of tensor integral operators that depend on the gradient of the drift function $b_t(x_1,x_2)$ of the nonlinear diffusion. These  integro-differential operators are described in
section~\ref{integro-diff-sec}. The last section, section~\ref{s-diff-form-sec} provides some differential formulae as well as some exponential decays estimates of the norm of these operators w.r.t. the time horizon.

\subsection{Fr\'echet differential}\label{not-tensor-sec}

We let $[n]$ stands for the set of $n$ multiple indexes $i=(i_1,\ldots,i_n)\in \Ia^n$ over some finite set $\Ia$. Notice that
 $[n_1]\times[n_2]=[n_1+n_2]$.
We denote by $\Ta_{p,q}(\Ia)$ the space of $(p,q)$-tensor $X$ with real entries $(X_{i,j})_{(i,j)\in [p]\times [q]}$. Given a 
 $(p_1,q_1)$-tensor $X$ and a  $(p_2,q_2)$-tensor $Y$ we denote by $( X\otimes Y)$
 the $((p_1+q_1),(p_2+q_2))$-tensor defined by
 
 $$
( X\otimes Y)_{(i,j),(k,l)}:=X_{i,k}~Y_{j,l}
 $$

For a given $(p_1,q)$-tensor $X$ and  a given $(q,p_2)$ tensor $Y$,  the product $XY$ and 
the transposition $Y^{\prime}$ are the $(p_1,p_2)$ and $(p_2,q)$ tensors with entries
$$
\forall (i,j)\in [p_1]\times [p_2]\qquad
(XY)_{i,j}:=\sum_{k\in [q]}X_{i,k}Y_{k,j}\quad \mbox{\rm and}\quad Y^{\prime}_{j,k}=Y_{k,j}
$$
We equip $\Ta_{p,q}(\Ia)$ with the Frobenius inner product
$$
\langle X,Y\rangle:=\tr(XY^{\prime}):=\sum_{i\in [p]}(XY^{\prime})_{i,i}\quad \mbox{\rm and the norm}\quad \Vert X\Vert_{\tiny Frob}:=\sqrt{\tr(XX^{\prime})}
$$
Identifying $(1,0)$-tensors $\Ta_{1,0}(\Ia)=\RR^{\Ia}$ with column vectors $(X_{i})_{i\in \Ia}\in \RR^{\Ia}$ the above quantities coincide with the conventional Euclidian inner product and norm
on the product space $\RR^{\Ia}$. When $\Ia=\{1,\ldots,d\}$ we simplify notation and we set $\RR^{d}$ instead of $\RR^{\{1,\ldots,d\}}$.
For any tensors $X$ and $Y$ with appropriate dimensions, using Cauchy-Schwartz inequality we check that
$$
\langle X,Y\rangle^2\leq \Vert X\Vert_{\tiny Frob}~\Vert Y\Vert_{\tiny Frob}\quad \mbox{\rm and}\quad \Vert XY\Vert_{\tiny Frob}\leq \Vert X\Vert_{\tiny Frob}~\Vert Y\Vert_{\tiny Frob}
$$

Let $\HH(\Ta_{p,q}(\Ia)):=\LL_2((\Omega,\FF,\PP),\Ta_{p,q}(\Ia))$ be the Hilbert space of $\Ta_{p,q}(\Ia)$-valued random variables defined on some probability space $(\Omega,\FF,\PP)$, equipped with the inner product $$\langle X,Y\rangle_{\,\HH(\Ta_{p,q}(\Ia))}=\EE(\langle X,Y\rangle)\quad\mbox{\rm and the norm}\quad
\Vert X\Vert_{\,\HH(\Ta_{p,q}(\Ia))}:=\langle X,X\rangle_{\,\HH(\Ta_{p,q}(\Ia))}^{1/2}$$ induced by the inner product $\langle X,Y\rangle$ on $\Ta_{p,q}(\Ia)$. 
We denote by $\EE(X)=\EE(X_{i,j})_{(i,j)\in [p]\times [q]}$  the entry-wise 
expected value of a $(p,q)$-tensor.

When $\Ia=\{1,\ldots,d\}$ and $(p,q)=(1,0)$ the space $\HH(\Ta_{p,q}(\Ia))$ coincides with  be the Hilbert space $\HH(\RR^d)=\LL_2((\Omega,\FF,\PP),\RR^d)$ of square integrable $\RR^d$-valued and $\FF$-measurable random variables.

We denote by $$\HH_n(\Ta_{p,q}(\Ia)):=\LL_2((\Omega,\FF_n,\PP),\Ta_{p,q}(\Ia))$$ 
the non decreasing sequence of Hilbert spaces associated with some increasing filtration $\FF_n\subset \FF_{n+1}$.

In Landau notation, we recall that a function 
$$
F:X\in\HH_1(\Ta_{p_1,q_1}(\Ia))~\mapsto ~F(X)\in\HH_2(\Ta_{p_2,q_2}(\Ja))
$$
is said to be Fr\'echet differentiable at $X$ if there exists a continuous map $$
X\in\HH_1(\Ta_{p,q}(\Ia))~\mapsto ~ \partial F(X)\in \mbox{\rm Lin}(\HH_1(\Ta_{p_1,q_2}(\Ia)),\HH_2(\Ta_{p_2,q_2}(\Ja)))
$$  such that 
$$
F(X+Y)=F(X)+\partial F(X)\cdot Y+\mbox{\rm o}\left(Y\right)
$$
\subsection{Tensor integral operators}\label{tensor-integral-op-sec}

Let $\Ba(E,\Ta_{p,q}(\Ia))$ be the set of bounded measurable functions from a measurable space $E$ into some tensor space $\Ta_{p,q}(\Ia)$.
Signed measures $\mu$ on $E$ act on bounded measurable functions $g$ from $E$ into $\RR$. We extend these integral operators to tensor valued functions 
$g=(g_{i,j})_{(i,j)\in [p]\times [q])}\in\Ba(E,\Ta_{p,q}(\Ia))$ by setting for any $(i,j)\in [p]\times [q]$
$$
 \mu(g)_{i,j}= \mu(g_{i,j}):=\int~\mu(dx)~g_{i,j}(x)\quad \mbox{\rm and we set}\quad\mu(g): =\int~\mu(dx)~g(x)
$$
Let $(E,\Ea)$ and $(F,\Fa)$ be some pair of measurable spaces.
A $(p,q)$-tensor integral operator  
$$
\Qa~:~g\in\Ba(F,\Ta_{q,r}(\Ia))\mapsto \Qa(g)\in\Ba(E,\Ta_{p,r}(\Ia))$$ is defined for 
$r\geq 0$ and $g\in \Ba(\Fa,\Ta_{q,r}(\Ia))$ by the tensor valued and measurable function
 $\Qa(g)$ with entries given $x\in E$ and $(i,j)\in ([p]\times [r])$ by the integral formula
$$
\Qa(g)_{i,j}(x)=\sum_{k\in [q]}~\int_{F}~\Qa_{i,k}(x,d\overline{x})~g_{k,j}(\overline{x})
$$
for some collection of integral operators $\Qa_{i,k}(x_1,dx_2)$ from $\Ba(E,\RR)$ into $\Ba(F,\RR)$. We also consider the operator
norm
$$
\vertiii{\Qa}:=\sup_{\Vert g\Vert\leq 1}{\Vert \Qa(g)\Vert}\quad \mbox{\rm for some tensor norm $\Vert\point\Vert$}
$$

The tensor product $(\Qa^1\otimes \Qa^2)$ of a couple of $(p_i,q_i)$-tensor integral operators
$$
\Qa^i~:~g\in\Ba(F_i,\Ta_{q_i,r_i}(\Ia))\mapsto \Qa(g)\in\Ba(E_i,\Ta_{p_i,r_i}(\Ia))\quad \mbox{\rm with}\quad i=1,2
$$
is a  $(p,q)$-tensor integral operator
$$
\Qa^1\otimes \Qa^2~:~h\in\Ba(F,\Ta_{q,r}(\Ia))\mapsto \Qa(g)\in\Ba(E,\Ta_{p,q}(\Ia))
$$
with the product spaces
$$
 (E,F):=(E_1\times E_2,F_1\times F_2)\quad  \mbox{\rm and}\quad
( p,q,r)=(p_1+p_2,q_1+q_2, r_1+r_2)
$$
The entries of $(\Qa^1\otimes \Qa^2)(h)$ are given for any $x=(x_1,x_2)$ and any pair of multi-indices
 $i=(i_1,i_2)\in ([p_1]\times  [p_2])$, $j=(j_1,j_2) \in ([r_1]\times  [r_2])$ by the integral formula
$$
(\Qa^1\otimes \Qa^2)(h)_{i,j}(x)=\sum_{k\in ( [q_1]\times[q_2])}~\int_{F_1\times F_2}~(\Qa^1\otimes \Qa^2)_{i,k}(x,dy)~h_{k,j}(y)
$$
with the tensor product measures defined for any $k=(k_1,k_2)\in ([q_1]\times [q_2])$ and any $y=(y_1,y_2)$ by
$$
(\Qa^1\otimes \Qa^2)_{(i_1,i_2),(k_1,k_2)}((x_1,x_2),d(y_1,y_2)):=
\Qa^1_{i_1,k_1}(x_1,dy_1)~\Qa^2_{i_2,k_2}(x_2,dy_2)
$$

\subsection{Variational equations}\label{var-eq-sec}
The gradient and the Hessian of a multivariate smooth function $h(x)=(h_{i}(x))_{i\in [p]}$ is defined by the $(1,p)$ and $(2,p)$ tensors $\nabla h(x)$  
and $\nabla^2 h(x)$
with entries given for any $1\leq  k,l\leq d$ and $i\in [p]$ by the formula
\begin{equation}\label{grad-def}
\nabla h(x)_{k,i}=\partial_{x_k}h_{i}(x)\quad \mbox{\rm and}\quad \nabla^2 h(x)_{(k,l),i}=\partial_{x_k}\partial_{x_l}h_{i}(x)
\end{equation}
 We consider the tensor valued functions $ b_t^{[k_1,k_2]}$ and $b_t^{[k_1,k_2,k_3]}$  defined for any $k_1,k_2,k_3=1,2$ by 
$$
 b_t^{[k_1,k_2]}:=(\nabla_{x_{k_1}}\otimes\nabla_{x_{k_2}})b_t\quad\mbox{and}\quad 
 b_t^{[k_1,k_2,k_3]}:=(\nabla_{x_{k_1}}\otimes\nabla_{x_{k_2}}\otimes\nabla_{x_{k_3}})b_t
$$
with the $(2,1)$ and $(3,1)$-tensor valued functions
$$
\left( b_t^{[k_1,k_2]}\right)_{(i_1,i_2),j}=\partial_{x_{k_1}^{i_1}}\partial_{x_{k_2}^{i_2}}b^j_t
\quad\mbox{\rm and}\quad
\left( b_t^{[k_1,k_2,k_3]}\right)_{(i_1,i_2,i_3),j}=\partial_{x_{k_1}^{i_1}}\partial_{x_{k_2}^{i_2}}\partial_{x_{k_3}^{i_3}}b^j_t
$$
 In the above display, $\partial_{x_{k}^{i}}b^j_t(x_1,x_2)$ stands for the partial derivative of the scalar function
 $b_t^j(x_1,x_2)$  w.r.t. the coordinate $x_k^i$, with the drift function  $b_t(x_1,x_2)$ from  $\RR^{2d}$ into $\RR^d$ introduced in section~\ref{ref-desciption},
 In the same vein, $\partial_{x_{k_1}^{i_1}}\partial_{x_{k_2}^{i_2}}b^j_t(x_1,x_2)$ and $\partial_{x_{k_1}^{i_1}}\partial_{x_{k_2}^{i_2}}\partial_{x_{k_3}^{i_3}}b^j_t(x_1,x_2)$
 stands for the second and third  partial derivatives of
 $b_t^j(x_1,x_2)$  w.r.t. the coordinates $x_{k_1}^{i_1}$, $x_{k_2}^{i_2}$ and $x_{k_3}^{i_3}$ with $k_1,k_2,k_3\in \{1,2\}$. 
 
 For any $\mu\in P_2(\RR^d)$ and $x_1\in\RR^d$ we also consider the tensor functions 
 $$
b_t^{[1]}(x_1,\mu)_{i,j}:=\int\mu(dx_2)~\partial_{x_1^i} b_t^j(x_1,x_2)\qquad
b_t^{[1,1]}(x_1,\mu)_{(i_1,i_2),j}:=\int\mu(dx_2)~\partial_{x_{1}^{i_1}}\partial_{x_{1}^{i_2}}b^j_t(x_1,x_2)
$$

Recalling that $b_t(x,\phi_{s,t}(\mu))$  has continuous and uniformly bounded derivatives up to the third order, the stochastic flow $
 x\mapsto X_{s,t}^{\mu}(x)$ is a twice differentiable function of the initial state $x$. 
In addition, when $(H)$ holds
 the gradient $\nabla X^{\mu}_{s,t}(x)$ of the diffusion flow $X^{\mu}_{s,t}(x)$ satifies the $(d\times d)$-matrix valued
         stochastic diffusion equation
\begin{equation}\label{tau-1-estimates}
\partial_t \,\nabla X^{\mu}_{s,t}(x)=\nabla X^{\mu}_{s,t}(x)~b_t^{[1]}\left(X^{\mu}_{s,t}(x),\phi_{s,t}(\mu)\right)~\Longrightarrow~
\Vert\nabla X^{\mu}_{s,t}(x)\Vert_2\leq e^{-\lambda_1(t-s)}
\end{equation}
The above estimate is a direct consequence of well known log-norm estimates for exponential semigroups, see for instance~\cite{coppel} as well as
 section 1.3 in the recent article~\cite{Bishop/DelMoral:2019}.

We have the stochastic tensor evolution equation
$$
\begin{array}{l}
\partial_t\, \nabla^2 X^{\mu}_{s,t}(x)\\
\\
\displaystyle=\nabla^2 X^{\mu}_{s,t}(x)~b_t^{[1]}(X^{\mu}_{s,t}(x),\phi_{s,t}(\mu))
+\left[\nabla X^{\mu}_{s,t}(x)\otimes \nabla X^{\mu}_{s,t}(x)\right]~b_t^{[1,1]}(X^{\mu}_{s,t}(x),\phi_{s,t}(\mu))\end{array}
$$
{This implies that
$$
\partial_t\Vert\nabla^2 X^{\mu}_{s,t}(x)\Vert^2_{\tiny Frob}\leq~-2\lambda_1~\Vert\nabla^2 X^{\mu}_{s,t}(x)\Vert_{\tiny Frob}^2+2\Vert b^{[1,1]}\Vert_{\tiny Frob}~\Vert\nabla X^{\mu}_{s,t}(x)\Vert_{\tiny Frob}^2~\Vert\nabla^2 X^{\mu}_{s,t}(x)\Vert_{\tiny Frob}
$$ from which we check that
$$
\partial_t\Vert\nabla^2 X^{\mu}_{s,t}(x)\Vert_{\tiny Frob}\leq~-\lambda_1~\Vert\nabla^2 X^{\mu}_{s,t}(x)\Vert_{\tiny Frob}+\Vert b^{[1,1]}\Vert_{\tiny Frob}~\Vert\nabla X^{\mu}_{s,t}(x)\Vert_{\tiny Frob}^2
$$
Using (\ref{tau-1-estimates}), this yields the estimate
\begin{equation}\label{tau-2-estimates}
\Vert\nabla^2 X^{\mu}_{s,t}(x)\Vert_{\tiny Frob}\leq c_1~e^{-\lambda_1(t-s)}~\int_s^t~e^{\lambda_1(u-s)}~\Vert\nabla X^{\mu}_{s,u}(x)\Vert_{\tiny Frob}^2~du\leq c_2~e^{-\lambda_1(t-s)}
\end{equation}
More generally, using the multivariate version of the de Fa\`a di Bruno derivation formula~\cite{ma-2009} (see also formula
(\ref{faa-di-bruno-ap}) in the appendix), for any $n\geq 1$  we also check the uniform estimate
\begin{equation}\label{tau-n-estimates}
\Vert\nabla^n X^{\mu}_{s,t}(x)\Vert_{\tiny Frob}\leq c_{n}~e^{-\lambda_1(t-s)}
\end{equation}}
A detailed proof is provided in the appendix, on page~\pageref{tau-n-estimates-proof}.

\subsection{Differential of Markov semigroups}

We have the commutation formula
\begin{equation}\label{commutation-tensor}
\nabla \circ P_{s,t}^{\mu}=\Pa_{s,t}^{\mu}\circ \nabla 
\end{equation}
with the $(1,1)$-tensor integral operator $\Pa_{s,t}^{\mu}$  defined for any $x\in \RR^d$ and any differentiable function $f$ on $\RR^d$ by the formula
\begin{equation}\label{def-Pa}
\Pa_{s,t}^{\mu}(\nabla f)(x):=
\EE\left[\nabla X_{s,t}^{\mu}(x)~\nabla f(X_{s,t}^{\mu}(x))\right]
\end{equation}
The tensor product of $\Pa_{s,t}^{\mu}$ is also given by the $(2,2)$-tensor integral operator
$$
\left(\Pa_{s,t}^{\mu}\right)^{\otimes 2}(h)(x_1,x_2):=
\EE\left[\left[\nabla X_{s,t}^{\mu}(x_1)\otimes\nabla \overline{X}_{s,t}^{\mu}(x_2)\right]~h\left(X_{s,t}^{\mu}(x_1),\overline{X}_{s,t}^{\mu}(x_2)\right)\right]
$$
In the above display, $\overline{X}_{s,t}^{\mu}(x)$ stands for an independent copy of $X_{s,t}^{\mu}(x)$ and $h=(\nabla\otimes \nabla) g$ stands for the matrix valued function defined in (\ref{def-nabla-nabla}). 
We also have the commutation formula
$$
\left(\Pa_{s,t}^{\mu}\right)^{\otimes 2}\circ (\nabla\otimes \nabla)=(\nabla\otimes \nabla)\circ \left(P^{\mu_0}_{s,t}\right)^{\otimes 2}
$$

In the same vein, we have the second order differential formula
  \begin{equation}\label{ref-nabla-Pa}
 \nabla^2P^{\mu}_{s,t}(f)=\Pa^{[2,1],\mu}_{s,t}(\nabla f)+\Pa^{[2,2],\mu}_{s,t}(\nabla^2 f)
  \end{equation}
with  the  $(2,1)$ and $(2,2)$-tensor integral operators
\begin{eqnarray}
\Pa^{[2,1],\mu}_{s,t}(\nabla f)(x)&:=&  \EE\left[\nabla^2 X_{s,t}^{\mu}(x)~\nabla f(X_{s,t}^{\mu}(x))\right]\nonumber\\
\Pa^{[2,2],\mu}_{s,t}(\nabla^2 f)(x)&:=&  \EE\left[\left(\nabla X_{s,t}^{\mu}(x)\otimes \nabla X_{s,t}^{\mu}(x)\right)~\nabla^2 f(X_{s,t}^{\mu}(x))\right]
\label{def-U-V}
\end{eqnarray}
Iterating the above procedure, the $n$-th differential of $P^{\mu}_{s,t}(f)$ at any order $n\geq 1$ takes the form
$$
 \nabla^nP^{\mu}_{s,t}(f)=\sum_{1\leq k\leq n}\Pa^{[n,k],\mu}_{s,t}(\nabla^k f)
$$ 
for some integral operators $\Pa^{[n,k],\mu}_{s,t}$. 
For instance,  we have the third order differential formula
  \begin{equation}\label{ref-nabla-3-Pa}
 \nabla^3P^{\mu}_{s,t}(\nabla f)=\Pa^{[3,1],\mu}_{s,t}(\nabla f)+\Pa^{[3,2],\mu}_{s,t}(\nabla^2 f)+\Pa^{[3,3],\mu}_{s,t}(\nabla^3 f)
  \end{equation}
with  the  $(2,1)$ and $(2,2)$-tensor integral operators
\begin{eqnarray}
\Pa^{[3,1],\mu}_{s,t}(\nabla f)(x)&:=&  \EE\left[\nabla^3 X_{s,t}^{\mu}(x)~\nabla f(X_{s,t}^{\mu}(x))\right]\nonumber\\
\Pa^{[3,2],\mu}_{s,t}(\nabla^2 f)(x)&:=&  \EE\left[\left(\nabla^2 X_{s,t}^{\mu}(x)\frownotimes \nabla X_{s,t}^{\mu}(x)\right)~\nabla^2 f(X_{s,t}^{\mu}(x))\right]\nonumber\\
\Pa^{[3,3],\mu}_{s,t}(\nabla^3 f)(x)&:=&\EE\left[\left(\nabla X_{s,t}^{\mu}(x)\otimes \nabla X_{s,t}^{\mu}(x)\otimes \nabla X_{s,t}^{\mu}(x)\right)~\nabla^3 f(X_{s,t}^{\mu}(x))\right]
\label{def-nabla-3}
\end{eqnarray}
with the $\frownotimes$-tensor product of type $(3,2)$ given for any $i=(i_1,i_2,i_3)$ and $l=(l_1,l_2)$ by
$$
\begin{array}{l}
\left(\nabla^2 X_{s,t}^{\mu}(x)\frownotimes \nabla X_{s,t}^{\mu}(x)\right)_{i,l}
:=
\left(\nabla^2 X_{s,t}^{\mu}(x)\otimes \nabla X_{s,t}^{\mu}(x)\right)_{((i_1,i_2),i_3),l}\\
\\
\hskip3cm+
\left(\nabla^2 X_{s,t}^{\mu}(x)\otimes \nabla X_{s,t}^{\mu}(x)\right)_{((i_2,i_3),i_1),l}
+\left(\nabla^2 X_{s,t}^{\mu}(x)\otimes \nabla X_{s,t}^{\mu}(x)\right)_{((i_3,i_1),i_2),l}
\end{array}
$$
The above formulae remains valid for any column vector multivariate function $f=(f_{i})_{1\leq i\leq d}$. 
An explicit description of the integral operators $\Pa^{[n,k],\mu}_{s,t}$ for any $1\leq k\leq n$ can be obtained using 
multivariate derivations and combinatorial manipulations, see for instance  the multivariate version of the de Fa\`a di Bruno  derivation formulae
(\ref{faa-di-bruno-ap}) and
(\ref{faa-di-bruno-ap-2}) in the appendix.
Following the proof of (\ref{tau-n-estimates}) 
we also check the uniform estimates
\begin{equation}\label{estimates-U-V}
\sup_{1\leq k\leq n}{\vertiii{ \Pa^{[n,k],\mu}_{s,t}}}\leq c_n~e^{-\lambda_1(t-s)}
\end{equation}
Using the moment estimates (\ref{ref-moments}) for any  $\mu\in P_2(\RR^d)$, $m,n\geq 0$, 
 and any $s\leq t$, we also check the rather crude estimate
\begin{equation}\label{crude-est-rev}
\vertiii{P^{\mu}_{s,t}}_{\Ca^n_m(\RR^d)\rightarrow \Ca^n_m(\RR^d)}\vee \vertiii{ (P^{\mu}_{s,t})^{\otimes 2}}_{\Ca^n_m(\RR^{2d})\rightarrow \Ca^n_m(\RR^{2d})}\leq c_{m,n}(t)~\left[1+\Vert e\Vert_{\mu,2}\right]^{m}
\end{equation}
For instance, using the de Fa\`a di Bruno derivation formula (\ref{faa-di-bruno-ap-2}) for any function $f\in\Ca^n_m(\RR^d)$
such that $\Vert f\Vert_{\Ca^n_m(\RR^d)}\leq 1$ and for any  $0\leq k\leq n$ 
we check that
$$
\Vert\nabla^kP^{\mu}_{s,t}(f)(x)\Vert=\Vert\EE\left(\nabla^k(f\circ X^{\mu}_{s,t})(x)\right)\Vert\leq c_{n,m}(t)~\EE\left((1+\Vert X_{s,t}^{\mu}(x)\Vert)^m\right)
$$
The estimates (\ref{ref-moments}) implies that 
$$
\Vert\nabla^kP^{\mu}_{s,t}(f)(x)\Vert\leq c_{n,m}(t)~\left(\Vert x\Vert+\Vert e\Vert_{\mu,2}\right)^m\leq  c_{n,m}(t)~(1+\Vert x\Vert)^m~
\left(1\vee \Vert e\Vert_{\mu,2}\right)^m
$$
from which we conclude that
$$
\vertiii{P^{\mu}_{s,t}}_{\Ca^n_m(\RR^d)\rightarrow \Ca^n_m(\RR^d)}\leq c_{m,n}(t)~\left[1+\Vert e\Vert_{\mu,2}\right]^{m}
$$
\subsection{Bismut-Elworthy-Li extension formulae}\label{bel-formula-sec}
We have the Bismut-Elworthy-Li formula
    \begin{equation}\label{bismut-omega}
\nabla P^{\mu}_{s,t}(f)(x)=
    \EE\left(f(X_{s,t}^{\mu}(x))~  \tau^{\mu,\omega}_{s,t}(x)\right)\quad \mbox{\rm with}\quad   \tau^{\mu,\omega}_{s,t}(x):=\int_s^t~   \partial_u \omega_{s,t}(u)~\nabla X^{\mu}_{s,u}(x)~dW_u
    \end{equation}
The above formula is valid for any function $\omega_{s,t}:u\in [s,t]\mapsto \omega_{s,t}(u)\in \RR $ of the following form 
    \begin{equation}\label{bismut-omega-varphi}
\omega_{s,t}(u)=\varphi\left((u-s)/(t-s)\right)~\Longrightarrow   \partial_u \omega_{s,t}(u)=\frac{1}{t-s}~\partial\varphi\left((u-s)/(t-s)\right)~
      \end{equation}
   for some non decreasing differentiable function $\varphi$ on $[0,1]$ with bounded continuous derivatives and such that
  $$
  (\varphi(0),\varphi(1))=(0,1)\Longrightarrow \omega_{s,t}(t)-\omega_{s,t}(s)=1
  $$
  In the same vein, for any $s\leq u\leq t$ we have
      \begin{equation}\label{bismut-omega-2}
\nabla^2 P^{\mu}_{s,t}(f)(x)=
    \EE\left(f(X_{s,t}^{\mu}(x))~  \left[\tau^{[2],\mu,\omega}_{s,u}(x)+\nabla X_{s,u}^{\mu}(x)~\tau^{\phi_{s,u}(\mu),\omega}_{u,t}(X_{s,u}^{\mu}(x))\,\tau^{\mu,\omega}_{s,u}(x)^{\prime}\right]\right)    \end{equation}
    with the stochastic process
    $$
     \tau^{[2],\mu,\omega}_{s,t}(x):=\int_s^t~   \partial_u \omega_{s,t}(u)~\nabla^2 X^{\mu}_{s,u}(x)~dW_u
$$
 Besides the fact that $X^{\mu}_{s,t}(x)$ is a nonlinear diffusion, the proof of the above formula follows the same proof as the one provided in~\cite{aht-03,bismut,Elworthy,xm-li,thompson} in the context of diffusions on differentiable
  manifolds. For the convenience of the reader, a detailed proof is provided in the appendix on page~\pageref{bismut-omega-proof}.  Using (\ref{bismut-omega}), for any $f$ s.t. $\Vert f\Vert\leq 1$   we check that
  \begin{eqnarray*}
 \Vert \nabla P^{\mu}_{s,t}(f)\Vert^2&\leq& 
 \EE\left( \Vert\tau^{\mu,\omega}_{s,t}(x)\Vert^2 \right)\\
 &\leq&  \int_s^t~e^{-2\lambda_1 (u-s)}~
 \Vert \partial_u \omega^{s,t}(u)\Vert^2~du=~\frac{1}{t-s}~ \int_0^1~e^{-2\lambda_1 (t-s)v}~
 \left( \partial\varphi(v)\right)^2~dv
\end{eqnarray*}
Let $\varphi_{\epsilon}$ with $\epsilon\in ]0,1[$ be some differentiable function on $[0,1]$ null on $[0,1-\epsilon]$ and such that $\vert \partial\varphi_{\epsilon}(u)\vert\leq c/\epsilon$
and
$(\varphi_{\epsilon}(1-\epsilon),\varphi(1))=(0,1)$, for instance we can choose
$$
\varphi(u)=\left\{
\begin{array}{ccl}
0&\mbox{\rm if}&u\in [0,1-\epsilon]\\
\displaystyle1+\cos{\left(\left(1+\frac{1-u}{\epsilon}\right)\frac{\pi}{2}\right)}&\mbox{\rm if}&u\in [1-\epsilon,1]
\end{array}
\right.
$$
In this situation, we find the rather crude uniform estimate
\begin{equation}\label{bismut-est}
 \Vert \nabla P^{\mu}_{s,t}(f)\Vert^2\leq \left(\frac{c}{\epsilon}\right)^2\frac{1}{t-s}~ \int_{1-\epsilon}^1~e^{-2\lambda_1 (t-s)v}~dv
\Longrightarrow
 \Vert \nabla P^{\mu}_{s,t}(f)\Vert\leq \frac{c}{\epsilon}~\frac{1}{\sqrt{t-s}}~e^{-\lambda_1(1-\epsilon) (t-s)}
\end{equation}
In the same vein, combining (\ref{bismut-omega-2}) with the estimate (\ref{tau-2-estimates}) for any $\epsilon\in ]0,1[$ and $u\in ]s,t[$ we also check the rather crude uniform estimate
$$
 \Vert \nabla^2 P^{\mu}_{s,t}(f)\Vert\leq \frac{c_1}{\epsilon}~\frac{1}{\sqrt{u-s}}~e^{-\lambda_1 (u-s) (1-\epsilon)}+\frac{c_2}{\epsilon^2}\frac{1}{\sqrt{(t-u)(u-s)}}~e^{-\lambda_1 (u-s) }~e^{-\lambda_1 (t-s) (1-\epsilon)}
$$
Choosing $u=s+(1-\epsilon)(t-s)$ in the above display we readily check that  
\begin{equation}\label{bismut-est-P2}
 \Vert \nabla^2 P^{\mu}_{s,t}(f)\Vert\leq \frac{c_1}{\epsilon\sqrt{1-\epsilon}}~\frac{1}{\sqrt{t-s}}~e^{-\lambda_1 (1-\epsilon)^2(t-s)}+\frac{c_2}{\epsilon^2}\frac{1}{\sqrt{\epsilon(1-\epsilon)}}~\frac{1}{t-s}~ e^{-2\lambda_1 (t-s) (1-\epsilon)}
\end{equation}

\subsection{Integro-differential operators}\label{integro-diff-sec}

Let $\BB^{\mu}_{s,t}(x_0,x_1)$  be  the matrix-valued function defined for any $(x_0,x_1)\in \RR^{2d}$, $\mu\in P_2(\RR^d)$ and any $s\leq t$  by
the formulae
\begin{equation}\label{def-BB-s-t}
\BB^{\mu}_{s,t}(x_0,x_1):=\nabla_{x_0}b_{s,t}^{\mu}(x_0,x_1)\quad\mbox{\rm with}\quad
b_{s,t}^{\mu}(x_0,x_1):=\EE\left[b_{t}\left(x_1,X_{s,t}^{\mu}(x_0)\right)\right]
\end{equation}
For instance, for the linear model discussed in (\ref{lin-case-intro}) we have
$$
\BB^{\mu}_{s,t}(x_0,x_1)^{\prime}=B_2~e^{(t-s)B_1}\quad\mbox{\rm and}\quad
b_{s,t}^{\mu}(x_0,x_1)=B_1x_1+B_2\left[e^{(t-s)B_1}(x_0-\mu(e))+e^{(t-s)[B_1+B_2]}~\mu(e)\right]
$$

We also consider the collection Weyl chambers $ [s,t]_n$ defined for any $n\geq 1$ by
\begin{eqnarray*}
[s,t]_n&:=&\left\{u=(u_1,\ldots,u_n)\in [s,t]^n~:~s\leq u_1\leq\ldots\leq u_n\leq t\right\}\quad \mbox{\rm and set}\quad
du:=du_1\ldots du_n
\end{eqnarray*}
We consider the space-time Weyl chambers
\begin{eqnarray}
\Delta_{s,t}&:=&\cup_{n\geq 1}\Delta_{s,t}^n\quad \mbox{\rm with}\quad \Delta_{s,t}^n:=[s,t]_n\times \RR^{nd} \label{def-Delta-s-t}
\end{eqnarray}
The coordinates of a generic point $(u,y)\in \Delta_{s,t}^n$ for some $n\geq 1$ are denoted by $$u=(u_1,\ldots,u_n)\in [s,t]_n\quad \mbox{\rm and} \quad y=(y_1,\ldots,y_n)\in\RR^{nd}$$ We also use the convention $u_0=s$ and $u_{n+1}=t$.  We  consider the  measures  $\Phi_{s,u}(\mu)$ on $\Delta_{s,t}$ given on every set $\Delta_{s,t}^n$ and any $n\geq 1$ by 
$$
\Phi_{s,u}(\mu)(d(u,y))=\phi_{s,u}(\mu)(dy)~du
$$
with the tensor product measures
$$
\phi_{s,u}(\mu)(dy):=\phi_{s,u_1}(\mu)(dy_1)\ldots\phi_{s,u_n}(\mu)(dy_n)
$$
 
 \begin{defi}
Let $b^{\mu}_{s,u}(x,y)$ be the function defined for any  $\mu\in P_2(\RR^d)$, $x\in\RR^d$,
and any $(u,y)\in \Delta_{s,t}^n$ and $n\geq 1$ by the formula
 \begin{equation}\label{defi-tau}
b^{\mu}_{s,u}(x,y)^{\prime}:= b_{s,u_1}^{\mu}\left(x,y_1\right)^{\prime}~ \prod_{1\leq k<n}\BB^{\phi_{s,u_k}(\mu)}_{u_k,u_{k+1}}(y_k,y_{k+1})
\end{equation}
\end{defi}
In the above display the product of matrices is understood as a directed product from $k=1$ to $k=(n-1)$.
For instance, for the linear model discussed in (\ref{lin-case-intro}) we have
$$
b^{\mu}_{s,u}(x,y)=B_2~e^{(u_n-u_{n-1})B_1}\ldots B_2~e^{(u_2-u_1)B_1}~b^{\mu}_{s,u_1}(x,y_1)
$$
For any $x\in\RR^d$,
and any $(u,y)\in \Delta_{s,t}^n$ and $n\geq 1$ we also set
  \begin{equation}\label{def-T}
  \BB_{u,t}^{\phi_{s,u}(\mu)}(y,x):=\BB_{u_n,t}^{\phi_{s,u_n}(\mu)}(y_n,x)\quad \mbox{\rm and}\quad
 \Pa^{\,\phi_{s,u}(\mu)}_{u,t}(\nabla f)(y):= \Pa^{\,\phi_{s,u_n}(\mu)}_{u_n,t}(\nabla f)(y_n)
  \end{equation}

\begin{defi} 
For any $\mu_0,\mu_1\in P_2(\RR^d)$ and $s\leq t$ we let $Q^{\mu_1,\mu_0}_{s,t}$  be  
the operator  defined  on differentiable functions $f$ on $\RR^d$ by 
  \begin{equation}\label{def-Q-Qa}
Q^{\mu_1,\mu_0}_{s,t}(f):=
\Qa^{\mu_1,\mu_0}_{s,t}(\nabla f)
  \end{equation}
with the $(0,1)$-tensor integral operator $\Qa^{\mu_1,\mu_0}_{s,t}$  defined by the integral formula
$$
\Qa^{\mu_1,\mu_0}_{s,t}(\nabla f)(x):=\int_{\Delta_{s,t}}~\Phi_{s,u}(\mu_1)(d(u,y))~b^{\mu_0}_{s,u}(x,y)^{\prime}~\Pa^{\,\phi_{s,u}(\mu_0)}_{u,t}(\nabla f)(y)
$$
\end{defi}
Recall that $b_{t}(x,y)$ is differentiable at any order with uniformly bounded derivatives.
Thus, using the estimates (\ref{ref-moments}) and (\ref{tau-n-estimates}), for any  $m,n\geq 0$, $\mu_0,\mu_1\in P_{m\vee 2}(\RR^d)$ we have
\begin{equation}\label{est-Q-1}
\Vert Q^{\mu_1,\mu_0}_{s,t}\Vert_{\Ca^1_m(\RR^d)\rightarrow \Ca^n_1(\RR^d)}\leq c_{m,n}(t)~\rho_{m\vee 2}(\mu_0,\mu_1)
\end{equation}

\begin{defi}
Let $p^{\mu_1,\mu_0}_{s,t}$ be  the function defined for any $s\leq t$ and $x,z\in \RR^d$ by the formula
    \begin{eqnarray}
p_{s,t}^{\mu_1,\mu_0}(x,z)^{\prime}
&=&b^{\mu_0}_{s,t}(x,z)^{\prime}+\int_{\Delta_{s,t}}~\Phi_{s,u}(\mu_1)(d(u,y))~b^{\mu_0}_{s,u}(x,y)^{\prime}
~ \BB^{\,\phi_{s,u}(\mu_0)}_{u,t}(y,z)\label{p-def-al}
  \end{eqnarray}
  \end{defi}
In this notation, we readily check the following proposition.
\begin{prop}\label{prop-ref-rev}
 The $(0,1)$-tensor   integral operator $\Qa^{\mu_1,\mu_0}_{s,t}$ can be rewritten as follows:
$$
\Qa^{\mu_1,\mu_0}_{s,t}(\nabla f)(x)=\int_{\Delta_{s,t}^1}~\Phi_{s,u}(\mu_1)(d(u,y))~p^{\mu_1,\mu_0}_{s,u}(x,y)^{\prime}~\Pa^{\phi_{s,u}(\mu_0)}_{u,t}(\nabla f)(y)
$$ 
\end{prop}
 
 For instance, for the linear model discussed in (\ref{lin-case-intro}) 
  the function $p_{s,t}^{\mu_1,\mu_0}(x,z)$ defined in 
 (\ref{p-def-al}) reduces to
  \begin{equation}\label{linear-example-d}
\begin{array}{l}
\displaystyle 
p_{s,t}^{\mu_1,\mu_0}(x,z)=B_1\,z+B_2\,e^{(t-s)(B_1+B_2)}\,x\\
\\
\displaystyle +B_2~\left[\int_s^te^{(t-u)(B_1+B_2)}\,B_1\,e^{(u-s)(B_1+B_2)}~du~\mu_1(e)+\int_s^te^{(t-u)(B_1+B_2)}\,B_2\,e^{(u-s)(B_1+B_2)}~du~\mu_0(e)\right]
\end{array}
\end{equation}
We check this claim expanding in (\ref{p-def-al}) the exponential series coming from the integration over the set $\Delta_{s,t}$.
A detailed proof of the above formula is provided in the appendix on page~\pageref{linear-example-d-proof}.
 
\subsection{Some differential formulae}\label{s-diff-form-sec}

The matrix $\nabla_{y_0}b_{s,t}^{\mu}(y_0,y_1)$
defined in (\ref{def-BB-s-t}) can alternatively be written as follows
$$
\nabla_{y_0}b_{s,t}^{\mu}(y_0,y_1)=\Pa^{\mu}_{s,t}\left(b^{[2]}_{t}(y_1,\point)\right)(y_0)=\EE\left[\nabla X_{s,t}^{\mu}(y_0)~b^{[2]}_{t}(y_1,X_{s,t}^{\mu}(y_0))\right]
$$
We also have the $(2,1)$ and $(3,1)$-tensor formulae
\begin{eqnarray*}
\nabla_{y_0}^2b_{s,t}^{\mu}(y_0,y_1)&=&\Pa^{[2,1],\mu}_{s,t}(b^{[2]}_{t}(y_1,\point))(y_0)+\Pa^{[2,2],\mu}_{s,t}(b^{[2,2]}_{t}(y_1,\point))(y_0)\\
\nabla^3_{y_0}b_{s,t}^{\mu}(y_0,y_1)&=&\Pa^{[3,1],\mu}_{s,t}(b^{[2]}_{t}(y_1,\point))(y_0)+\Pa^{[3,2],\mu}_{s,t}(b^{[2,2]}_{t}(y_1,\point))(y_0)+\Pa^{[3,3],\mu}_{s,t}(b^{[2,2,2]}_{t}(y_1,\point))(y_0)
\end{eqnarray*}
For any $(u,y)\in \Delta_{s,t}^n$ with $n\geq 1$ and for  any $k\geq 1$ we have the $(k,1)$-tensor formulae
 \begin{equation}\label{def-BB-k}
 \nabla^{k}_{y_0}b^{\mu}_{s,u}(y_0,y)=\BB_{s,u}^{[k],\mu}(y_0,y):=\nabla^k_{y_0}b_{s,u_1}^{\mu}(y_0,y_1)\prod_{1\leq k<n}\BB^{\phi_{s,u_k}(\mu)}_{u_k,u_{k+1}}(y_k,y_{k+1})
 \end{equation} We consider the 
 $(n,1)$-tensor valued function
 \begin{eqnarray*}
q_{s,t}^{[n],\mu_1,\mu_0}(x,z)
&:=&
\BB^{[n],\mu_0}_{s,t}(x,z)
+\int_{\Delta_{s,t}}~\Phi_{s,u}(\mu_1)(d(u,y))~\BB_{s,u}^{[n],\mu_0}(x,y)
~ \BB^{\,\phi_{s,u}(\mu_0)}_{u,t}(y,z)
\end{eqnarray*}
and we use the convention 
$$\BB^{[0],\mu_0}_{s,t}(x,z)=b^{\mu_0}_{s,t}(x,z)^{\prime}
\quad
\mbox{\rm so that}\quad
q_{s,t}^{[0],\mu_1,\mu_0}(x,z)=p_{s,t}^{\mu_1,\mu_0}(x,z)^{\prime}
$$
  {For instance, for the linear model discussed in (\ref{lin-case-intro}) and (\ref{linear-example-d}) the above objects reduce to
$$
q_{s,t}^{[1],\mu_1,\mu_0}(x,y)^{\prime}=B_2\,e^{(B_1+B_2)(t-s)}\quad\mbox{\rm and}\quad
 \forall n\geq 2\quad q_{s,t}^{[n],\mu_1,\mu_0}(x,y)=0  
$$
}

In this notation, we have the following proposition.
\begin{prop}\label{prop-d-n-Q}
For any $n\geq 0$ the $n$-th differential of the operator $Q^{\mu_1,\mu_0}_{s,t}$ is given by the formula
 $$
  \nabla^n Q^{\mu_1,\mu_0}_{s,t}(f)= \Qa^{[n],\mu_1,\mu_0}_{s,t}(\nabla f)
$$
with the $(n,1)$-tensor integral operator given by 
\begin{equation}\label{def-Qa-n}
\Qa^{[n],\mu_1,\mu_0}_{s,t}(\nabla f)(x):= \int_{\Delta^1_{s,t}}~\Phi_{s,u}(\mu_1)(d(u,y))~q^{[n],\mu_1,\mu_0}_{s,u}(x,y)~\Pa^{\mu_0}_{u,t}(\nabla f)(y)
\end{equation}
In addition, when condition $(H)$ is satisfied for any $n\geq 1$ we have the exponential estimates
  \begin{equation}\label{def-lambda-hat-2}
 \vertiii{ \Qa^{[n],\mu_1,\mu_0}_{s,t}}\leq c_n~e^{- \lambda(t-s)}\quad \mbox{for some $\lambda>0$}
\end{equation}
\end{prop}
\proof
The proof of the first assertion follows from (\ref{p-def-al}). More precisely, using (\ref{p-def-al}) we have
$$
\nabla^n_x \,p_{s,t}^{\mu_1,\mu_0}(x,y)
=q^{[n],\mu_1,\mu_0}_{s,t}(x,y)
$$
On the other hand, by proposition~\ref{prop-ref-rev} we also have
\begin{eqnarray*}
  \nabla^n Q^{\mu_1,\mu_0}_{s,t}(f)(x)&=&  \nabla^n\Qa^{\mu_1,\mu_0}_{s,t}(\nabla f)(x)\\
  &=&\int_{\Delta_{s,t}^1}~\Phi_{s,u}(\mu_1)(d(u,y))~\nabla^n_x\,p^{\mu_1,\mu_0}_{s,u}(x,y)~\Pa^{\phi_{s,u}(\mu_0)}_{u,t}(\nabla f)(y)= \Qa^{[n],\mu_1,\mu_0}_{s,t}(\nabla f)(x)
\end{eqnarray*}
  This ends the proof of the first assertion.
When condition $(H)$ is satisfied, for any $x\in \RR^d$ and $(u,y)\in \Delta_{s,t}^n$ we have
\begin{equation}\label{estimate-BB}
\Vert \BB^{\mu}_{s,t}(y_0,y_1)\Vert_2\leq \Vert b^{[2]}\Vert_2~e^{-\lambda_1(t-s)}\quad\mbox{\rm and}\quad
\Vert \BB_{s,u}^{\mu}(x,y)\Vert_2\leq  \Vert b^{[2]}\Vert_2^n~e^{-\lambda_1(u_n-s)}
\end{equation}
Using (\ref{tau-n-estimates}) we also check the uniform estimate
\begin{equation}\label{estimate-q-2}
 \Vert 
  q_{s,t}^{[n],\mu_1,\mu_0}(x,y)  
  \Vert\leq c_n~e^{- \lambda_{1,2}(t-s)}
\end{equation}
The end of the proof is now a consequence of (\ref{tau-1-estimates}). 
\cqfd

\begin{prop}
For any $n\geq 0$ any bounded function $f$ on $\RR^d$ and for  any function $\omega$ of the form (\ref{bismut-omega-varphi}) we have the Bismut-Elworthy-Li formula
  \begin{equation}\label{def-Q-bismut-ref}
  \nabla^n Q^{\mu_1,\mu_0}_{s,t}(f)= \int_{\Delta^1_{s,t}}~\Phi_{s,u}(\mu)(d(u,y))~q^{[n],\mu_1,\mu_0}_{s,u}(x,y)~ \EE\left(f(X^{\mu_0}_{u,t}(y))~  \tau^{\mu_0,\omega}_{u,t}(y)\right)
\end{equation}
In the above display, $\tau^{\mu,\omega}_{u,t}(y)$ stands for the stochastic process defined in (\ref{bismut-omega}). In addition, when condition $(H)$ is satisfied we have the exponential estimates
  \begin{equation}\label{def-lambda-hat-bismut}
\Vert  \nabla^n Q^{\mu_1,\mu_0}_{s,t}(f)\Vert\leq c_n~e^{- \lambda(t-s)}~\Vert f\Vert\quad \mbox{for some $\lambda>0$}
\end{equation}
\end{prop}
\proof
The proof of the first assertion is a direct application of the Bismut-Elworthy-Li formula (\ref{bismut-omega}). More precisely, using  (\ref{bismut-omega}) we have
$$
\Pa^{\mu_0}_{u,t}(\nabla f)(y)= \EE\left(f(X^{\mu_0}_{u,t}(y))~  \tau^{\mu_0,\omega}_{u,t}(y)\right)
$$
The formula (\ref{def-Q-bismut-ref}) is now a direct consequence of (\ref{def-Qa-n}).

We check (\ref{def-lambda-hat-bismut}) combining (\ref{bismut-est}) with (\ref{estimate-q-2}).
This ends the proof of the proposition.
\cqfd

When $n=1$ we drop the upper index and we write $\left(\BB_{s,u}^{\mu},q_{s,t}^{\mu_1,\mu_0}\right)$ instead of  $\left(
 \BB_{s,u}^{[1],\mu},q_{s,t}^{[1],\mu_1,\mu_0}\right)$.

 The operators discussed above are indexed by a pair of measures $(\mu_0,\mu_1)$. To simplify notation, when $\mu_1=\mu_0=\mu$ we suppress one of the indices and
 we write
 $
( Q^{\mu}_{s,t},\Qa^{[n],\mu}_{s,t})$
 and
  $(p^{\mu}_{s,t},q^{[n],\mu}_{s,t})$
  instead of  $(Q^{\mu,\mu}_{s,t},\Qa^{[n],\mu,\mu}_{s,t})$ and $(p^{\mu,\mu}_{s,t}, q^{[n],\mu,\mu}_{s,t})$.

\section{Tangent processes}\label{tangent-sec}

The tangent process associated with the diffusion flow $\psi_{s,t}(Y)$ introduced in (\ref{eq-diff}) is given for any $U\in \HH_s(\RR^d)$ by the evolution equation
 \begin{equation}\label{eq-diff-lin-eq}
\partial_t(\partial \psi_{s,t}(Y)\cdot U)=\partial B_t(\psi_{s,t}(Y))\cdot  (\partial \psi_{s,t}(Y)\cdot U)
 \end{equation}
 In the above display,  $\partial B_t(X)\in \mbox{\rm Lin}(\HH_t(\RR^d),\HH_t(\RR^d))$ stands for the Fr\'echet differential
 of the drift function $B_t$ defined for any $Z\in \HH_t(\RR^d)$ by
$$
\partial B_t(X)\cdot Z=\EE\left(
\nabla_{x_1}\,b_t(X,\overline{X})^{\prime}~Z+\nabla_{x_2}\,b_t(X,\overline{X})^{\prime}~\overline{Z}~|~\FF_t\right)
$$
where $(\overline{X},\overline{Z})$ stands for an independent copy of $(X,Z)$.

 \subsection{Spectral estimate}\label{spectral-sec}
 This section is mainly concerned with the proof of theorem~\ref{theo-intro-1}.
 
 For any pair of random variables $Z_1,Z_2\in \HH_t(\RR^d)$ we have the duality formula
 $$
\langle Z_1, \partial B_t(X)\cdot  Z_2\rangle_{\,\HH_t(\RR^d)}=\langle \partial B_t(X)^{\star}\cdot  Z_1, Z_2\rangle_{\,\HH_t(\RR^d)}
$$
with the dual operator $\partial B_t(X)^{\star}$ defined by the formula
$$
\partial B_t(X)^{\star}\cdot  Z_1:=\EE\left(b_t^{[1]}(X,\overline{X})~Z_1+b_t^{[2]}(\overline{X},X)~\overline{Z}_1~\vert~\FF_t\right)
$$
 In the above display, $(\overline{X},\overline{Z}_1)$ stands for an independent copy of $(X,Z_1)$. The symmetric part of $ \partial B_t(X)$ is given by the formula
  $$
 \partial B_t(X)_{\tiny sym}:=
 \frac{1}{2}\left[\partial B_t(X)+\partial B_t(X)^{\star}\right]
 $$

We are now in position to prove theorem~\ref{theo-intro-1}. 

The first assertion is a direct consequence of the evolution equation
 \begin{eqnarray*}
2^{-1} \partial_t\,\Vert \partial \psi_{s,t}(Y)\cdot U\Vert^2_{\,\HH_t(\RR^d)}&=&\langle (\partial \psi_{s,t}(Y)\cdot U), \partial B_t(\psi_{s,t}(Y))_{\tiny sym}\cdot (\partial \psi_{s,t}(Y)\cdot U)\rangle_{\,\HH_t(\RR^d)}
 \end{eqnarray*}

Whenever $(H)$ is met we have $  \partial B_t(X)_{\tiny sym}\leq -\lambda_0~I$ for some $\lambda_0>0$. In this situation, the r.h.s. estimate in (\ref{theo-intro-1-eq})
is a direct consequence of (\ref{theo-intro-1-eq-log-norm}). 
Given an independent copy  $(\overline{X},\overline{Z}_2) $ of $(X,Z_2)$ we have
\begin{eqnarray*}
2~\langle Z_1,\partial B_t(X)^{\star}\cdot   Z_2\rangle_{\,\HH_t(\RR^d)}&=&\EE\left(\left\langle
\left[\begin{array}{c}
Z_1\\
\overline{Z}_1
\end{array}\right],A_t(X,\overline{X})\left[\begin{array}{c}
 Z_2\\
\overline{Z}_2
\end{array}\right]\right\rangle\right)\\&=&2~\langle \partial B_t(X)\cdot  Z_1,  Z_2\rangle_{\,\HH_t(\RR^d)}
\end{eqnarray*}
This yields the log-norm estimate
$$
A_t(X,\overline{X})_{\tiny sym}\leq -\lambda_0~I\Longrightarrow
~ \partial B_t(X)_{\tiny sym}\leq -\lambda_0~I
$$
The proof of   theorem~\ref{theo-intro-1} is now completed.
\cqfd

\subsection{Dyson-Phillips expansions}

 In the further development of this section we shall denote by
$$(\overline{\psi}_{s,t},\overline{U},\overline{X}^{\mu}_{s,t},\overline{Y})\quad\mbox{\rm and} \quad(\overline{\psi}^n_{s,t},\overline{U}^n,\overline{X}^{\mu,n}_{s,t},\overline{Y}^n)_{n\geq 0}$$  a collection of
independent copies of the stochastic flows $(\psi_{s,t},X^{\mu}_{s,t})$ and some given $U,Y\in\HH_s(\RR^d)$. To simplify notation, we also set 
$$X_{s,t}:=\psi_{s,t}(Y)\qquad\overline{X}_{s,t}:=\overline{\psi}_{s,t}(\overline{Y})\quad\mbox{\rm and}\quad\overline{X}^n_{s,t}:=\overline{\psi}^n_{s,t}(\overline{Y}^n)$$
We are now in position to state and prove the main result of this section.
\begin{theo}
The tangent process $\partial \psi_{s,t}$ is given  
for any $U\in\HH_s(\RR^d)$ and any $Y\in\HH_s(\RR^d)$ with distribution $\mu\in P_2(\RR^d)$ by the Dyson-Phillips
series
\begin{equation}\label{tangent-hilbert}
\begin{array}{l}
\displaystyle\partial \psi_{s,t}(Y)\cdot U=\nabla X^{\mu}_{s,t}(Y)^{\prime}~U\\
\\
\hskip.3cm\displaystyle+\sum_{n\geq 1}\int_{[s,t]_n}~
\left(\nabla X^{\phi_{s,u_n}(\mu)}_{u_n,t}\right)\left(X_{s,u_n}\right)^{\prime}
~
{\EE}\left(
\left[
\prod_{1\leq k\leq n}{\BB}_{u_{k-1},u_k}^{\phi_{s,u_{k-1}}(\mu)}\left(\overline{X}^{k-1}_{s,u_{k-1}},\overline{X}^k_{s,u_k}\right)
\right]^{\prime}
\overline{U}~\vert~\FF_{u_n}\right)~du
\end{array}
\end{equation}
with  the boundary conventions 
$$
 u_0=s\qquad \overline{X}^0_{s,u_1}=\overline{X}_{s,u_1}\quad \mbox{and}\quad
\overline{X}^n_{s,u_n}=X_{s,u_n}\quad \mbox{for any $n\geq 1$}
$$

\end{theo}

\proof
For any $s\leq u\leq t$ and $x\in\RR^d$ we have
$$
\partial_t \,\nabla X^{\mu}_{s,t}(x)^{-1}=-~b_t^{[1]}\left(X^{\mu}_{s,t}(x),\phi_{s,t}(\mu)\right)~\nabla X^{\mu}_{s,t}(x)^{-1}
$$
and
$$
\nabla X^{\mu}_{s,t}(x)=\nabla X^{\mu}_{s,u}(x)~\left(\nabla X_{u,t}^{\phi_{s,u}(\mu)}\right)\left(X^{\mu}_{s,u}(x)\right)
$$
In addition, for any $s\leq u\leq t$ and $x_0,x_1\in\RR^d$ we have
$$
\nabla_{x_0} b_t(x_1,X_{u,t}^{\phi_{s,u}(\mu)}(x_0))=\nabla X_{u,t}^{\phi_{s,u}(\mu)}(x)~b^{[2]}_t(x_1,X_{u,t}^{\phi_{s,u}(\mu)}(x_0))
$$
Combining the above with (\ref{eq-diff-lin-eq}) we check that
$$
\begin{array}{l}
\displaystyle\partial_t\left(\left(\nabla X^{\mu}_{s,t}(Y)^{-1}\right)^{\prime}~(\partial \psi_{s,t}(Y)\cdot U)\right)\\
\\
\displaystyle=\left(\nabla X^{\mu}_{s,t}(Y)^{-1}\right)^{\prime}~{\EE}\left(
\nabla b_t\left(\psi_{s,t}(Y),\overline{X}^{\mu}_{s,t}(\point)\right)(\overline{Y})^{\prime}~\left(\nabla \overline{X}^{\mu}_{s,t}(\overline{Y})^{-1}\right)^{\prime}
(\partial \overline{\psi}_{s,t}(\overline{Y})\cdot \overline{U})~\vert~\FF_t
\right)
\end{array}
$$
In the above display, $\nabla b_t\left(\psi_{s,t}(Y),\overline{X}^{\mu}_{s,t}(\point)\right)(\overline{Y})=\nabla h(\overline{Y})$ stands for the gradient of the random
function
$$h~:~x\mapsto h(x)=b_t\left(\psi_{s,t}(Y),\overline{X}^{\mu}_{s,t}(x)\right)
\quad\mbox{\rm
evaluated at $x=\overline{Y}$.}
$$
Equivalently, we have
$$
\begin{array}{l}
\displaystyle\left(\nabla X^{\mu}_{s,t}(Y)^{-1}\right)^{\prime}~(\partial \psi_{s,t}(Y)\cdot U)\\
\\
\displaystyle =U+\int_s^t\left(\nabla X^{\mu}_{s,u}(Y)^{-1}\right)^{\prime}~{\EE}\left(
\nabla b_u\left(\psi_{s,u}(Y),\overline{X}^{\mu}_{s,u}(\point)\right)(\overline{Y})^{\prime}~\left(\nabla \overline{X}^{\mu}_{s,u}(\overline{Y})^{-1}\right)^{\prime}
(\partial \overline{\psi}_{s,u}(\overline{Y})\cdot \overline{U})~\vert~\FF_u
\right)~du
\end{array}
$$
and therefore
$$
\begin{array}{l}
\displaystyle\partial \psi_{s,t}(Y)\cdot U=\left(\nabla X^{\mu}_{s,t}(Y)\right)^{\prime}U+\int_s^t~\left(\left(\nabla X^{\phi_{s,u}(\mu)}_{u,t}\right)\left(X_{s,u}^{\mu}(Y)\right)\right)^{\prime}\\
\\
\displaystyle~\hskip3cm\times {\EE}\left(
\nabla b_u\left(\psi_{s,u}(Y),\overline{X}^{\mu}_{s,u}(\point)\right)(\overline{Y})^{\prime}~\left(\nabla \overline{X}^{\mu}_{s,u}(\overline{Y})^{-1}\right)^{\prime}
(\partial \overline{\psi}_{s,u}(\overline{Y})\cdot \overline{U})~\vert~\FF_u
\right)~du
\end{array}
$$
Now, the end of the proof of (\ref{tangent-hilbert}) follows a simple induction, thus it is skipped.
\cqfd

\begin{cor}\label{cor-tangent}
For any $V\in\HH_t(\RR^d)$ and for any $Y\in\HH_s(\RR^d)$ with distribution $\mu\in P_2(\RR^d)$ we have
\begin{equation}\label{tangent-hilbert}
\begin{array}{l}
\displaystyle\partial \psi_{s,t}(Y)^{\star}\cdot V=\EE\left(\nabla X^{\mu}_{s,t}(Y) V~|~\FF_s\right)\\
\\
\hskip.3cm\displaystyle+\sum_{n\geq 1}\int_{[s,t]_n}~{\EE}\left(
\left[
\prod_{1\leq k\leq n}  \BB_{u_{k-1},u_k}^{\phi_{s,u_{k-1}}(\mu)}\left(\overline{X}^{k-1}_{s,u_{k-1}},\overline{X}^k_{s,u_k}\right)
\right]\left(\nabla \overline{X}^{\phi_{s,u_n}(\mu)}_{u_n,t}\right)\left(\overline{X}_{s,u_n}\right)
\overline{V}~|~\FF_s\right)~du
\end{array}
\end{equation}
with  the boundary conditions
$$
 u_0=s\quad \mbox{and}\quad \overline{X}^0_{s,u_1}=\psi_{s,u_1}(Y)
 \quad \mbox{and}\quad \overline{X}^n_{s,u_n}=\overline{X}_{s,u_n}
$$
\end{cor}
\subsection{Gradient semigroup analysis}\label{g-sg-sec}

This section is concerned with a gradient semigroup description of the dual of the tangent process.

 \begin{defi}
For any $\mu_0,\mu_1\in P_2(\RR^d)$ and $s\leq t$ we let $D_{\mu_1,\mu_0}\phi_{s,t}$  be  
the operator  defined  on differentiable functions $f$ on $\RR^d$ by 
\begin{equation}\label{def-D}
D_{\mu_1,\mu_0}\phi_{s,t}:=P^{\mu_0}_{s,t}+Q^{\mu_1,\mu_0}_{s,t} 
\end{equation}
In the above display, $Q^{\mu_1,\mu_0}_{s,t} $ stands for the operator defined in (\ref{def-Q-Qa}).
  \end{defi}
  {Rewritten in terms of expectation operators we have
  $$
  \begin{array}{l}
D_{\mu_1,\mu_0}\phi_{s,t}(f)(x)\\
\\
\displaystyle=\EE\left[(f\circ X^{\mu_0}_{s,t})(x)\right]+\sum_{n\geq 1}\int_{\Delta_{s,t}^n}~\Phi_{s,u}(\mu_1)(d(u,y))~\EE\left[b^{\mu_0}_{s,u}(x,y)^{\prime}~\nabla (f\circ X^{\phi_{s,u_n}(\mu_0)}_{u_n,t})(y_n)\right]
\end{array}
  $$}

Recall that $b_{t}(x,y)$ is differentiable at any order with uniformly bounded derivatives.
Thus, arguing as in the proof of (\ref{crude-est-rev}) and (\ref{est-Q-1})  for any  $m,n\geq 1$, $\mu_0,\mu_1\in P_{m\vee 2}(\RR^d)$ we have
\begin{equation}
\Vert D_{\mu_1,\mu_0}\phi_{s,t}\Vert_{\Ca^n_m(\RR^d)\rightarrow \Ca^n_m(\RR^d)}\leq c_{m,n}(t)~\rho_{m\vee 2}(\mu_0,\mu_1)\label{nabla-D}
\end{equation}
In the same vein, we check that
\begin{equation}\label{nabla-D-times-2}
\Vert \left(D_{\mu_1,\mu_0}\phi_{s,t}\right)^{\otimes 2}\Vert_{\Ca^{n+1}_m(\RR^{2d})\rightarrow \Ca^n_{m+1}(\RR^{2d})}\leq c_{m,n}(t)~\rho_{m\vee 2}(\mu_0,\mu_1)
\end{equation}
The proof of the above estimate is rather technical, thus it is housed in the appendix on page~\pageref{nabla-D-times-2-poof}.

\begin{rmk}\label{rmk-bismut}
Using the Bismut-Elworthy-Li formula
(\ref{def-Q-bismut-ref}), we extend the operators $D_{\mu_1,\mu_0}\phi_{s,t}$ with $s<t$ to non necessarily differentiable and bounded functions.
 
We also extend the operator $D_{\mu_1,\mu_0}\phi_{s,t}$ to tensor functions $f=(f_{i})_{i\in [n]}$ by considering the tensor function
with entries
\begin{equation}\label{def-D-extension}
D_{\mu_1,\mu_0}\phi_{s,t}(f)_i=D_{\mu_1,\mu_0}\phi_{s,t}(f_{i})
\end{equation}

In this situation, the function $p^{\mu_1,\mu_0}_{s,t}$  introduced in (\ref{p-def-al}) takes the form
$$
 p_{s,t}^{\mu_1,\mu_0}(x,z)=D_{\mu_1,\mu_0}\phi_{s,t}(b_t(z,\point))(x)
$$
\end{rmk}

Let $G_{t,\mu_1}$ be the  collection of integro-differential operators  indexed by $\mu_1\in P_2(\RR^d)$ defined by
$$
G_{t,\mu_1}(f)(x_2):=\int\mu_1(dx_1)~b_t(x_1,x_2)^{\prime}\,\nabla f(x_1)
$$
We also set
$$
H_{t,\mu_0,\mu_1}:=L_{t,\mu_0}+G_{t,\mu_1}\quad \mbox{\rm and}\quad H_{t,\mu_0}:=L_{t,\mu_0}+G_{t,\mu_0}
$$
In this notation, we have the first order expansion
\begin{equation}\label{ref-rev-key}
\mu_1L_{t,\mu_1}-\mu_0L_{t,\mu_0}=(\mu_1-\mu_0)L_{t,\mu_0}+(\mu_1-\mu_0)~G_{t,\mu_1}=(\mu_1-\mu_0)H_{t,\mu_0,\mu_1}
\end{equation}

  \begin{theo}\label{theo-1}
For any $m,n\geq 1$ and any $\mu_0,\mu_1\in P_{m\vee 2}(\RR^d)$
the operator $D_{\mu_1,\mu_0}\phi_{s,t}$ coincides with the evolution semigroup of the integro-differential operator
$
H_{t,\phi_{s,t}(\mu_0),\phi_{s,t}(\mu_1)}
$;
that is, we have the forward evolution equation
\begin{equation}\label{forward-eq}
\partial_tD_{\mu_1,\mu_0}\phi_{s,t}=D_{\mu_1,\mu_0}\phi_{s,t}\circ H_{t,\phi_{s,t}(\mu_0),\phi_{s,t}(\mu_1)}\quad \mbox{on}\quad \Ca^{n\vee 2}_m(\RR^d)
\end{equation}
In addition, for any $s\leq u< t$ we have the backward evolution equation 
\begin{equation}\label{backward-eq}
\partial_uD_{\phi_{s,u}(\mu_1),\phi_{s,u}(\mu_0)}\phi_{u,t}=-H_{u,\phi_{s,u}(\mu_0),\phi_{s,u}(\mu_1)}\circ D_{\phi_{s,u}(\mu_1),\phi_{s,u}(\mu_0)}\phi_{u,t}\quad \mbox{on}\quad \Ca^n_m(\RR^d)
\end{equation}
 \end{theo}
\proof
The proof of the forward equation (\ref{forward-eq}) is a direct consequence of the forward evolution equation 
$$
\partial_tP^{\mu_0}_{s,t}=P_{s,t}^{\mu_0}L_{t,\phi_{s,t}(\mu_0)}
$$
associated with the Markov semigroup $P^{\mu_0}_{s,t}$, thus it is skipped.
The semigroup property (\ref{sg-chain-rule}) yields
$$
\partial_u\left(D_{\mu_1,\mu_0}\phi_{s,u}\circ D_{\phi_{s,u}(\mu_1),\phi_{s,u}(\mu_0)}\phi_{u,t}\right)=0
$$
Combining the above with the forward equation (\ref{forward-eq}) we check that
$$
D_{\mu_1,\mu_0}\phi_{s,u}\circ \partial_uD_{\phi_{s,u}(\mu_1),\phi_{s,u}(\mu_0)}\phi_{u,t}=-D_{\mu_1,\mu_0}\phi_{s,u}\circ H_{u,\phi_{s,u}(\mu_0),\phi_{s,u}(\mu_1)}\circ D_{\phi_{s,u}(\mu_1),\phi_{s,u}(\mu_0)}\phi_{u,t}
$$
This implies that
$$
\left[\partial_uD_{\phi_{s,u}(\mu_1),\phi_{s,u}(\mu_0)}\phi_{u,t}\right]_{u=s}=-H_{s,\mu_0,\mu_1}D_{\mu_1,\mu_0}\phi_{u,t}
$$
from which we conclude that
$$
\begin{array}{l}
 \left[\partial_uD_{\phi_{s,u}(\mu_1),\phi_{s,u}(\mu_0)}\phi_{u,t}\right]_{u=v}\\
 \\
 =\left[\partial_uD_{\phi_{v,u}\left(\phi_{s,v}(\mu_1)\right),\phi_{v,u}\left(\phi_{s,v}(\mu_0)\right)}\phi_{u,t}\right]_{u=v}=-H_{s,\phi_{s,v}(\mu_0),\phi_{s,v}(\mu_1)}D_{\phi_{s,v}(\mu_1),\phi_{s,v}(\mu_0)}\phi_{v,t}
\end{array}
$$
This yields the backward evolution equation (\ref{backward-eq}).
This ends the proof of the theorem.
\cqfd

Next proposition is a direct consequence of (\ref{def-D}) combined with the formulae (\ref{commutation-tensor}) and (\ref{def-Qa-n}).

  \begin{prop}
We  have the commutation formula
\begin{equation}\label{commutation-D}
  \nabla \circ D_{\mu_1,\mu_0}\phi_{s,t}= \Da_{\mu_1,\mu_0}\phi_{s,t}\circ \nabla
\end{equation}
with the $(1,1)$-tensor integral operator given by the column vector function
  \begin{equation}\label{def-Da}
  \Da_{\mu_1,\mu_0}\phi_{s,t}(\nabla f)(x):=   \Pa^{\mu_0}_{s,t}(\nabla f)(x)+\int_{\Delta^1_{s,t}}\Phi_{s,v}(\mu_1)(d(v,y))~q^{\mu_1,\mu_0}_{s,v}(x,y)~\Pa^{\phi_{s,v}(\mu_0)}_{v,t}(\nabla f)(y)
\end{equation}
In addition, when condition $(H)$ is satisfied we have
\begin{equation}\label{commutation-D-etimate}
\vertiii{ \Da_{\mu_1,\mu_0}\phi_{s,t}}\leq c~e^{-\lambda(t-s)}\quad \mbox{for some $\lambda>0$}
  \end{equation}
  \end{prop}

\begin{rmk}\label{rmk-1}
Following remark~\ref{rmk-bismut}, using the Bismut-Elworthy-Li formula
(\ref{def-Q-bismut-ref}), we extend the gradient operators $\nabla D_{\mu_1,\mu_0}\phi_{s,t}$ with $s<t$ to measurable and bounded functions.   The exponential estimate stated in (\ref{def-lambda-hat-bismut}) are a direct consequence of the estimates presented in (\ref{def-lambda-hat-bismut}).

By (\ref{def-D-extension}) the commutation formula (\ref{commutation-D}) is also satisfied for  multivariate column functions $f$. In this situation
$  \Da_{\mu_1,\mu_0}\phi_{s,t}(\nabla f)$ is a $(d\times d)$-matrix valued function.

The proof of theorem ~\ref{theo-intro-2} is now a consequence of the estimate (\ref{commutation-D-etimate}) and the fact that
$$
  \partial_t \left[\phi_{s,t}(\mu_1)-   \phi_{s,t}(\mu_0)\right]=\left[\phi_{s,t}(\mu_1)-   \phi_{s,t}(\mu_0)\right]\circ H_{t,\phi_{s,t}(\mu_0),\phi_{s,t}(\mu_1)}
$$
More precisely, using  (\ref{ref-rev-key}) the above formula implies that
$$
\begin{array}{l}
  \partial_u\left( \left[\phi_{s,u}(\mu_1)-   \phi_{s,u}(\mu_0)\right]D_{\phi_{s,u}(\mu_1),\phi_{s,u}(\mu_0)}\phi_{u,t}\right)=0\\
  \\
  \Longrightarrow
\phi_{s,t}(\mu_1)-   \phi_{s,t}(\mu_0)  =(\mu_1-\mu_0)D_{\mu_1,\mu_0}\phi_{s,t}
  \end{array}
$$
\end{rmk}

 The operators discussed above are indexed by a pair of measures $(\mu_0,\mu_1)$. To simplify notation, when $\mu_1=\mu_0=\mu$ we suppress one of the parameter and
 we write
$(D_{\mu}\phi_{s,t},\Da_{\mu}\phi_{s,t})
  $
  instead of $(D_{\mu,\mu}\phi_{s,t},\Da_{\mu,\mu}\phi_{s,t})$.

\begin{theo}\label{theo-0}
For any $m,n\geq 1$, any  function $f\in \Ca^n_m(\RR^d)$ and any $Y\in\HH_s(\RR^d)$   with distribution $\mu\in P_2(\RR^d)$ we have
 the gradient formula
$$
\partial \psi_{s,t}(Y)^{\star}\cdot\nabla f(\psi_{s,t}(Y))=\nabla D_{\mu}\phi_{s,t}(f)(Y)=\Da_{\mu}\phi_{s,t}(\nabla f)(Y)
$$
\end{theo}
\proof
Given a smooth function $f$ on $\RR^d$ we have
$$
\langle \nabla f(\psi_{s,t}(Y)),\partial \psi_{s,t}(Y)\cdot U\rangle_{\,\HH_t(\RR^d)}=
\langle \partial \psi_{s,t}(Y)^{\star}\cdot\nabla f(\psi_{s,t}(Y)), U\rangle_{\,\HH_s(\RR^d)}
$$
Replacing $V$ by $\nabla f(\psi_{s,t}(Y))$ in (\ref{tangent-hilbert})
we check that
$$
\begin{array}{l}
\displaystyle\partial \psi_{s,t}(Y)^{\star}\cdot\nabla f(\psi_{s,t}(Y))\\
\\
\displaystyle=\nabla D_{\mu}\phi_{s,t}(f)(Y)=\EE\left(
\nabla \left(f\circ  X^{\mu}_{s,t}\right)(Y)~|~Y\right)\\
\\
\hskip.3cm\displaystyle+\sum_{n\geq 1}\int_{[s,t]_n}~\overline{\EE}\left(
\left[
\prod_{0\leq k< n}  \BB_{u_{k},u_{k+1}}^{\phi_{s,u_{k}}(\mu)}\left(\overline{X}^{k}_{s,u_{k}},\overline{X}^{k+1}_{s,u_{k+1}}\right)
\right]
\nabla \left(f\circ X^{\phi_{s,u_n}(\mu)}_{u_n,t}\right)\left(\overline{X}_{s,u_n}\right)~|~Y\right)~du
\end{array}
$$
This ends the proof of the theorem
\cqfd

\section{Taylor expansions}\label{sec-taylor}

This section is mainly concerned with the proof of the first and second order Taylor expansions stated in theorem~\ref{theo-intro-3} and 
theorem~\ref{theo-intro-4} . Section~\ref{sec-some-diff} presents some preliminary differential formulae used in the proof of the theorems.

\subsection{Some differential formulae}\label{sec-some-diff}

The commutation formula  (\ref{commutation-D}) takes the form
$$
\nabla D_{\mu_1,\mu_0}\phi_{s,t}(f)= \Da_{\mu_1,\mu_0}\phi_{s,t}(\nabla f)
$$
Combining (\ref{def-D}) with proposition~\ref{prop-d-n-Q}  and the second order formula (\ref{ref-nabla-Pa}) we also have
$$
  \nabla^2  D_{\mu_1,\mu_0}\phi_{s,t}(f)=  \nabla^2  P^{\mu_0}_{s,t}(f)+ \Qa^{[2],\mu_1,\mu_0}_{s,t}(\nabla f)
  =\Pa^{[2,1],\mu}_{s,t}(\nabla f)+\Pa^{[2,2],\mu}_{s,t}(\nabla^2 f)+ \Qa^{[2],\mu_1,\mu_0}_{s,t}(\nabla f)
$$
In summary, we have
 the first and second order differential formulae
  \begin{eqnarray}
    \nabla D_{\mu}\phi_{s,t}(f) &=& \Da_{\mu}\phi_{s,t}(\nabla f) \nonumber\\
  \nabla^2D_{\mu}\phi_{s,t}(f) &=&\Da_{\mu}\phi_{s,t}^{[2,1]}(\nabla f)+\Pa^{[2,2],\mu}_{s,t}(\nabla^2 f)\quad \mbox{with}\quad
  \Da_{\mu}\phi_{s,t}^{[2,1]}=\Pa^{[2,1],\mu}_{s,t}+\Qa^{[2],\mu}_{s,t}\label{dD-W}
\end{eqnarray}
Similar  formulae for $ \nabla D_{\mu_0,\mu_1}\phi_{s,t}$ and $ \nabla^2 D_{\mu_0,\mu_1}\phi_{s,t}$ can easily be found.
In the same vein, using (\ref{ref-nabla-3-Pa}) 
we check the third order differential formula
  \begin{equation}\label{dD-W-3}
  \begin{array}{l}
  \nabla^3D_{\mu}\phi_{s,t}(f) \\
  \\
  =\Da_{\mu}\phi_{s,t}^{[3,1]}(\nabla f)+\Pa^{[3,2],\mu}_{s,t}(\nabla^2 f)+\Pa^{[3,3],\mu}_{s,t}(\nabla^3 f)\quad \mbox{with}\quad
  \Da_{\mu}\phi_{s,t}^{[3,1]}:=\Pa^{[3,1],\mu}_{s,t}+\Qa^{[3],\mu}_{s,t}
  \end{array}
\end{equation}
In addition, when condition $(H)$ is satisfied we have the exponential estimates
  \begin{equation}\label{def-lambda-hat-3}
 \vertiii{ \Da_{\mu}\phi_{s,t}}\vee \vertiii{ \Da_{\mu}\phi_{s,t}^{[2,1]}}\vee\vertiii{ \Da_{\mu}\phi_{s,t}^{[3,1]}}\leq c~e^{- \lambda(t-s)}
 \quad \mbox{for some $\lambda>0$}
\end{equation}

\begin{defi}
 We let $S_{s,t}^{\mu}$ be  the operator  defined  for any differentiable function $f$ on $\RR^d$  by 
$$
S_{s,t}^{\mu}(f)=\Sa_{s,t}^{\mu}(\nabla f)
$$
with the $(0,1)$-tensor  integral operator $\Sa_{s,t}^{\mu}$ defined  by  the formula
\begin{eqnarray}
\Sa_{s,t}^{\mu}(\nabla f)(x_1,x_2)
&:=&
b_s(x_1,x_2)^{\prime}~\Da_{\mu}\phi_{s,t}(\nabla f)(x_1)+b_s(x_2,x_1)^{\prime}~\Da_{\mu}\phi_{s,t}(\nabla f)(x_2)\label{def-S}
\end{eqnarray}
\end{defi}

Using (\ref{nabla-D}) and (\ref{def-Da})   for any $m,n\geq 0$ and $\mu\in P_{m\vee 2}(\RR^d)$ we check that
\begin{equation}\label{estimate-S-m}
\Vert S^{\mu}_{s,t}\Vert_{\Ca^{n+1}_m(\RR^d)\rightarrow \Ca^{n}_{m+1}(\RR^{2d})}\leq c_{m,n}(t)~\rho_{m\vee 2}(\mu)
\end{equation}

We also have the differential formula
\begin{equation}\label{nabla-2-S}
(\nabla\otimes\nabla)\left(S^{\mu}_{s,t}(f)\right)=\SSS^{[2,1],\mu}_{s,t}(\nabla f)+\SSS^{[2,2],\mu}_{s,t}(\nabla^2 f)
\end{equation}
with the matrix valued functions
$$
\begin{array}{l}
\SSS^{[2,1],\mu}_{s,t}(\nabla f)(x_1,x_2)=b_s^{[1,2]}(x_1,x_2)~\Da_{\mu}\phi_{s,t}(\nabla f)(x_1)+b_s^{[2,1]}(x_2,x_1)~\Da_{\mu}\phi_{s,t}(\nabla f)(x_2)\\
\\
\displaystyle\hskip5cm+b_s^{[2]}(x_2,x_1)~\Da_{\mu}\phi_{s,t}^{[2,1]}(\nabla f)(x_2)^{\prime}+
\Da_{\mu}\phi_{s,t}^{[2,1]}(\nabla f)(x_1)~b_s^{[2]}(x_1,x_2)^{\prime}\\
\\
\SSS^{[2,2],\mu}_{s,t}(\nabla^2 f)(x_1,x_2):=b_s^{[2]}(x_2,x_1)~\Pa^{[2,2],\mu}_{s,t}(\nabla^2 f)(x_2)^{\prime}+
\Pa^{[2,2],\mu}_{s,t}(\nabla^2 f)(x_1)~b_s^{[2]}(x_1,x_2)^{\prime}
\end{array}$$
 When condition $(H)$ is satisfied we also have the exponential estimates
\begin{equation}\label{3s-2}
\vertiii{  \SSS^{[2,1],\mu}_{s,t}}\vee \vertiii{  \SSS^{[2,2],\mu}_{s,t}} \leq c~e^{-\lambda(t-s)}\quad \mbox{for some $\lambda>0$}
\end{equation}
In addition, using the Bismut-Elworthy-Li extension formulae and the estimates (\ref{def-lambda-hat-bismut-intro}) and
(\ref{def-lambda-hat-bismut-intro-hessian}), or any bounded measurable function $f$ on $\RR^d$ we check that
$$
\Vert  (\nabla\otimes\nabla)\left(S^{\mu}_{s,t}(f)\right)\Vert\leq c~\left(1\vee 1/(t-s)\right)~e^{- \lambda(t-s)}~\Vert f\Vert\quad \mbox{for some $\lambda>0$}
$$
\subsection{A first order expansion}\label{sec-first-order-remainder}
This section is mainly concerned with the proof of theorem~\ref{theo-intro-3}.
The next technical lemma is pivotal.
\begin{lem}
For any $m\geq 1$ for  any $\mu_0,\mu_1\in P_{m+1}(\RR^d)$ we have the second order expansion
\begin{equation}\label{a-2}
\begin{array}{l}
\phi_{s,t}(\mu_1)- \phi_{s,t}(\mu_0)\\
\\
\displaystyle=(\mu_1-\mu_0)D_{\mu_0}\phi_{s,t}+\frac{1}{2}~ \int_s^t~\left[\phi_{s,u}(\mu_1)-\phi_{s,u}(\mu_0)\right]^{\otimes 2}\circ S^{\phi_{s,u}(\mu_0)}_{u,t}~du
\quad \mbox{on}\quad \Ca^{n+1}_m(\RR^d)
\end{array}
\end{equation}
\end{lem}

\proof

Combining (\ref{ref-rev-key}) with the backward evolution equation (\ref{backward-eq}) 
we check that
$$
 \begin{array}{l}
\displaystyle
\partial_u\left\{\left[\phi_{s,u}(\mu_1)-\phi_{s,u}(\mu_0)\right]\circ D_{\phi_{s,u}(\mu_0)}\phi_{u,t}\right\}\\
\\
=\left[\phi_{s,u}(\mu_1)-\phi_{s,u}(\mu_0)\right]\circ \left[
H_{u,\phi_{s,u}(\mu_0),\phi_{s,u}(\mu_1)}-H_{u,\phi_{s,u}(\mu_0)}\right]\circ D_{\phi_{s,u}(\mu_0)}\phi_{u,t}\\
\\
=\displaystyle\left[\phi_{s,u}(\mu_1)-\phi_{s,u}(\mu_0)\right]\circ \left[
G_{u,\phi_{s,u}(\mu_1)}-G_{u,\phi_{s,u}(\mu_0)}\right]\circ D_{\phi_{s,u}(\mu)}\phi_{u,t}
 \end{array}
$$
On the other hand, we have
$$
 \left[
G_{u,\phi_{s,u}(\mu_1)}-G_{u,\phi_{s,u}(\mu_0)}\right](x_2):=\int \left(\phi_{s,u}(\mu_1)-\phi_{s,u}(\mu_0)\right)(dx_1)~b_u(x_1,x_2)^{\prime}\,\nabla f(x_1)
$$
Integrating $u$ from $u=s$ to $u=t$ we obtain the formula
$$
 \begin{array}{l}
\displaystyle
\left[\phi_{s,t}(\mu_1)-\phi_{s,t}(\mu_0)-(\mu_1-\mu_0)D_{\mu_0}\phi_{s,t}\right](f)\\
\\
\displaystyle=\frac{1}{2}
\int_s^t~\int~\left[\phi_{s,u}(\mu_1)-\phi_{s,u}(\mu_0)\right]^{\otimes 2}(d(x_1,x_2))~\\
\\
\hskip3cm\left[
b_u(x_1,x_2)^{\prime}~\nabla D_{\phi_{s,u}(\mu_0)}\phi_{u,t}(f)(x_1)+b_u(x_2,x_1)^{\prime}~\nabla D_{\phi_{s,u}(\mu_0)}\phi_{u,t}(f)(x_2)\right]~du
 \end{array}$$
The end of the lemma is now completed.
\cqfd

Combining the above lemma with (\ref{nabla-D-times-2}) and (\ref{estimate-S-m}) we check (\ref{s-o-r}) with the operator $D^2_{\mu_1,\mu_0}\phi_{s,t}$  defined for any $m,n\geq 0$ and $\mu_0,\mu_1\in P_{m+2}(\RR^d)$ by
\begin{equation}\label{def-D2}
D^2_{\mu_1,\mu_0}\phi_{s,t}:=\int_s^t~\left(D_{\mu_1,\mu_0}\phi_{s,u}\right)^{\otimes 2}\circ S^{\phi_{s,u}(\mu_0)}_{u,t}~du\in\mbox{\rm Lin}\left(\Ca^{n+2}_m(\RR^d),\Ca^{n}_{m+2}(\RR^{2d})\right)
\end{equation}

\begin{rmk}
The second order term in (\ref{s-o-r}) can alternatively be expressed in terms of the Hessian of the 
semigroup $D^2_{\mu_1,\mu_0}\phi_{s,t}$; that is, we have that
\begin{equation}\label{nabla-equiv-2}
\begin{array}{l}
(\mu_1-\mu_0)^{\otimes 2}D^2_{\mu_1,\mu_0}\phi_{s,t}(f)\\
\\
\displaystyle=\int_{[0,1]^2}~\EE\left(\langle \left[ (\nabla\otimes\nabla) D^2_{\mu_1,\mu_0}\phi_{s,t}(f)\right](Y_{\epsilon,\overline{\epsilon}}),(Y_1-Y_0)\otimes  (\overline{Y}_1-\overline{Y}_0)\rangle\right)~d\epsilon~ d\overline{\epsilon}
\end{array}
\end{equation}
with the interpolating path
$$
Y_{\epsilon,\overline{\epsilon}}:=(Y_0+\epsilon (Y_1-Y_0),
\overline{Y}_0+\overline{\epsilon}(\overline{Y}_1-\overline{Y}_0))
$$
In the above display, $ (\overline{Y}_1,\overline{Y}_0)$ stands for an independent copy of a pair of random variables $(Y_{0},Y_1)$ with distribution  $(\mu_{0},\mu_1)$. Also observe that
$$
(\mu_1-\mu_0)^{\otimes 2}D^2_{\mu_1,\mu_0}\phi_{s,t}=(\mu_1-\mu_0)^{\otimes 2}\overline{D}^2_{\mu_1,\mu_0}\phi_{s,t}
$$
with the centered second order operator
$$
\begin{array}{l}
\overline{D}^2_{\mu_1,\mu_0}\phi_{s,t}(f)(x_1,x_2)\\
\\
\displaystyle:=\left[(\delta_{x_1}-\mu_0)\otimes (\delta_{x_2}-\mu_0)\right]D^2_{\mu_0}\phi_{s,t}(f)\\
\\
\displaystyle=\int_{[0,1]^2}~\EE\left(\langle \left[ (\nabla\otimes\nabla) D^2_{\mu_1,\mu_0}\phi_{s,t}(f)\right](Y_{\epsilon,\overline{\epsilon}}(x_1,x_2)),(x_1-Y_0)\otimes  (x_2-\overline{Y}_0)\rangle\right)~d\epsilon ~d\overline{\epsilon}
\end{array}
$$
In the above display, $Y_{\epsilon,\overline{\epsilon}}(x_1,x_2)$ stands for the interpolating path
$$
Y_{\epsilon,\overline{\epsilon}}(x_1,x_2):=(Y_0+\epsilon (x_1-Y_0),
\overline{Y}_0+\overline{\epsilon}(x_2-\overline{Y}_0))
$$
\end{rmk}
  \begin{prop}
We have commutation formula
\begin{equation}\label{commutation-D-Da}
(\nabla\otimes\nabla)\circ\left(D_{\mu_1,\mu_0}\phi_{s,t}\right)^{\otimes 2}=\left(\Da_{\mu_1,\mu_0}\phi_{s,t}\right)^{\otimes 2}\circ (\nabla\otimes\nabla)
   \end{equation}
In addition, we have the estimate
\begin{equation}\label{estimate-D2}
      \vertiii{\left(\Da_{\mu_1,\mu_0}\phi_{s,t}\right)^{\otimes 2}} \leq c~
  e^{-\lambda(t-s)}\quad \mbox{for some $\lambda>0$}
   \end{equation}
\end{prop}
\proof
The proof of the first assertion is a consequence of the commutation formula (\ref{commutation-D}).
Letting $h=(\nabla\otimes\nabla) g$ we have 
  $$
  \begin{array}{l}
\left(\Da_{\mu_1,\mu_0}\phi_{s,t}\right)^{\otimes 2}( h)(x_1,x_2)= \left(  \Pa^{\mu_0}_{s,t}\right)^{\otimes 2}(h)(x_1,x_2)\\
  \\
   \displaystyle \hskip1cm+\int_{\Delta_{s,t}^1}\Phi_{s,v}(\mu_1)(d(u,y))~q^{\mu_1,\mu_0}_{s,u}(x_2,y)~\left(\Pa^{\mu_0}_{s,t}\otimes\Pa^{\phi_{s,u}(\mu_0)}_{u,t}\right)(h)(x_1,y)\\
   \\
      \displaystyle \hskip1cm+\int_{\Delta_{s,t}^1}\Phi_{s,v}(\mu_1)(d(u,y))~q^{\mu_1,\mu_0}_{s,u}(x_1,y)~\left(\Pa^{\phi_{s,u}(\mu_0)}_{u,t}\otimes \Pa^{\mu_0}_{s,t}\right)(h)(y,x_2)\\
      \\
     \displaystyle \hskip1cm+\int_{\Delta_{s,t}^1\times \Delta_{s,t}^1}\Phi_{s,u}(\mu_1)(d(u,y))~\Phi_{s,v}(\mu_1)(d(v,z))\\
  \\
   \displaystyle \hskip3cm~\times~\left[q^{\mu_1,\mu_0}_{s,u}(x_1,y)\otimes q^{\mu_1,\mu_0}_{s,v}(x_2,z)\right]\left(\Pa^{\phi_{s,u}(\mu_0)}_{u,t}\otimes\Pa^{\phi_{s,v}(\mu_0)}_{v,t}\right)(h)(y,z)
  \end{array} $$
The proof of (\ref{estimate-D2}) now follows the same arguments as the ones we used in the proof of (\ref{commutation-D-etimate}), thus it is skipped.
This ends the proof of the proposition.
\cqfd

Combining (\ref{nabla-2-S}) with the commutation formula (\ref{commutation-D-Da}),
for any twice differentiable function $f$  and any $s\leq t$ and $\mu_0,\mu_1\in P_2(\RR^d)$
we check that
\begin{equation}\label{nabla-2-D-2}
(\nabla\otimes\nabla)D^2_{\mu_0,\mu_1}\phi_{s,t}(f):=\int_s^t~\left(\Da_{\mu_0,\mu_1}\phi_{s,u}\right)^{\otimes 2}\left(
\SSS^{[2,1],\phi_{s,u}(\mu_0)}_{u,t}(\nabla f)+\SSS^{[2,2],\phi_{s,u}(\mu_0)}_{u,t}(\nabla^2 f)\right)
~du
\end{equation}
with the operators $\SSS^{[2,k],\mu}_{s,t}$ discussed in (\ref{nabla-2-S}).
The proof of (\ref{estimate-D2-nabla}) is a direct consequence of (\ref{3s-2}) and (\ref{estimate-D2}). The proof of theorem~\ref{theo-intro-3}
is now completed.

\subsection{Second order analysis}\label{sec-o-analysis}
This short section is mainly concerned with the proof of the first part of theorem~\ref{theo-intro-4}.
\begin{lem}\label{lem-tex}
For any $m\geq 1$ and  $\mu_0,\mu_1\in P_{m+3}(\RR^d)$ and $s\leq t$ we have the tensor product formula
$$
\begin{array}{l}
\left(\phi_{s,t}(\mu_1)-\phi_{s,t}(\mu_0)\right)^{\otimes 2}\\
\\
\displaystyle=(\mu_1-\mu_0)^{\otimes 2}\left(D_{\mu_0}\phi_{s,t}\right)^{\otimes 2}+(\mu_1-\mu_0)^{\otimes 3}\,\Ra_{\mu_1,\mu_0}\phi_{s,t}\quad\mbox{on}\quad \Ca^{n+2}_m(\RR^{2d})
\end{array}$$
for some third order linear operator $\Ra_{\mu_1,\mu_0}\phi_{s,t}$ such that
$$
\Vert \Ra_{\mu_1,\mu_0}\phi_{s,t}\Vert_{\Ca^{n+2}_m(\RR^{2d})\rightarrow \Ca^n_{m+3}(\RR^{3d})}\leq c_{m,n}(t)~\rho_{m+2}(\mu_0,\mu_1)
$$
\end{lem}
The proof of the above lemma is rather technical, thus it is housed in the appendix, on page~\pageref{lem-tex-proof}.

Combining the above lemma with (\ref{a-2}) we readily check the second order decomposition (\ref{taylor-2-intro})
with a the remainder linear operator $D^3_{\mu_0,\mu_1}\phi_{s,t}$
such that
$$
\Vert D^3_{\mu_0,\mu_1}\phi_{s,t}\Vert_{\Ca^{n+3}_m(\RR^{d})\rightarrow \Ca^n_{m+4}(\RR^{3d})}\leq c_{m,n}(t)~\rho_{m+3}(\mu_0,\mu_1)
$$
This ends  the proof of the first part of theorem~\ref{theo-intro-4}. The proof of the second part of the theorem is provided in the appendix, on page~\pageref{D3-estimation-lem-proof}.

\subsubsection*{Acknowledgements}

The authors are  supported by the ANR Quamprocs on quantitative analysis of metastable processes.
P. Del Moral is also supported in part from the Chair Stress Test, RISK Management and Financial Steering, led by the French Ecole polytechnique and its Foundation and sponsored by BNP Paribas.

\section*{Appendix}

\subsection*{Proof of (\ref{backward-D-phi-mf})}\label{backward-D-phi-mf-proof}
It is easy to check that this first assertion is true for any collection of generators $L_{t,\mu}$,  thus we  skip the details.
The proof of the second assertion is a also a direct consequence of a more general result
which is valid for any collection of generators and non necessarily symmetric functions.

For any $N\geq 2$ and $x=(x^i)_{1\leq i\leq N}\in  (\RR^d)^N$ we set
$$
m(x)^{\odot 2}:=\frac{1}{N(N-1)}\sum_{1\leq i\not j\leq N}~\delta_{(x^i,x^j)}\quad\mbox{\rm and}\quad \Fa(x)=m(x)^{\otimes 2}(F)
$$
We extend  $L_{t,\mu}$ to  functions $F(x^1,x^2)$ on $\RR^{2d}$ by setting
$$
\begin{array}{l}
L_{t,\mu}^{(2)}(F)(x^1,x^2)\\
\\
\displaystyle=\frac{1}{2}\left(L_{t,\mu}(F(x^1,\point))(x^2)+L_{t,\mu}(F(\point,x^2))(x^1)+L_{t,\mu}(F(x^2,\point))(x^1)+L_{t,\mu}(F(\point,x^1))(x^2)\right)
\end{array}$$
For any function $F(x^1,x^2)$ on $\RR^{2d}$ we have
$$
m(x)^{\otimes 2}(F)=\left(1-\frac{1}{N}\right)~m(x)^{\odot 2}(F)+\frac{1}{N}~m(x)(C(F))=
m(x)^{\odot 2}\left(\left(1-\frac{1}{N}\right)~F+\frac{1}{N}~C^{(2)}(F)\right)
$$
with
$$
C(F)(x)=F(x,x)\quad\mbox{\rm and}\quad C^{(2)}(F)(x^1,x^2)=\frac{1}{2}~\left(C(F)(x^1)+C(F)(x^2)\right)
$$
This implies that
$$
\Lambda_t(\Fa)(x)=\left(1-\frac{1}{N}\right)~m(x)^{\odot 2}\left(L_{t,m(x)}^{(2)}(F)\right)+\frac{1}{N}~m(x)^{\odot 2}\left(L_{t,m(x)}^{(2)}(C^{(2)}(F))\right)
$$
Recalling that
$$
\displaystyle m(x)^{\odot 2}(F)=\frac{N}{N-1}~m(x)^{\otimes 2}(F)-\frac{1}{N-1}~m(x)^{\odot 2}(C^{(2)}(F)) 
$$
we conclude that
$$
\Lambda_t(\Fa)(x)=m(x)^{\otimes 2}\left(L_{t,m(x)}^{(2)}(F)\right)+\frac{1}{N}~m(x)^{\odot 2}\left(\Gamma^{(2)}_{L_{t,m(x)}}(F)\right)
$$
with the operator
$$
\Gamma^{(2)}_{L_{t,m(x)}}=L_{t,m(x)}^{(2)}\circ C^{(2)}-C^{(2)}\circ L_{t,m(x)}^{(2)}
$$
Observe that
$$
\begin{array}{l}
\displaystyle\Gamma^{(2)}_{L_{t,m(x)}}(F)(x^1,x^2)\\
\\
\displaystyle=\frac{1}{2}\left(L_{t,m(x)}(C(F))(x^1)+L_{t,m(x)}(C(F))(x^2)\right)\\
\\
\displaystyle-\frac{1}{2}~\left(L_{t,m(x)}(F(x^1,\point))(x^1)+L_{t,m(x)}(F(\point,x^1))(x^1)+
L_{t,m(x)}(F(x^2,\point))(x^2)+L_{t,m(x)}(F(\point,x^2))(x^2)
\right)
\end{array}
$$
This yields the formula
$$
\Gamma^{(2)}_{L_{t,m(x)}}(F)(x^1,x^2)=\frac{1}{2}\left(C(\Gamma^{(2)}_{L_{t,m(x)}}(F))(x^1)+C(\Gamma^{(2)}_{L_{t,m(x)}}(F))(x^2)\right)
$$
from which we conclude that
$$
\Lambda_t(\Fa)(x)=m(x)^{\otimes 2}\left(L_{t,m(x)}^{(2)}(F)\right)+\frac{1}{N}~m(x)\left(\Gamma_{L_{t,m(x)}}(F)\right)
$$
with the function $\Gamma_{L_{t,m(x)}}(F)$ defined for any $y\in\RR^d$ by
\begin{eqnarray*}
\Gamma_{L_{t,m(x)}}(F)(y)&=&C\left(\Gamma^{(2)}_{L_{t,m(x)}}(F)\right)(y,y)\\
&=&
L_{t,m(x)}(C(F))(y)-L_{t,m(x)}(F(y,\point))(y)-L_{t,m(x)}(F(\point,y))(y)
\end{eqnarray*}
The above formula
readily implies (\ref{backward-D-phi-mf}) as soon as $L_{t,\mu}$ is the collection of generators associated with the stochastic flow defined in (\ref{diff-st-ref-general}). This ends the proof of (\ref{backward-D-phi-mf}).\cqfd

\subsection*{Proof of (\ref{tau-n-estimates})}\label{tau-n-estimates-proof}

For any given  $1\leq m\leq n$, we denote by 
$\Pi_{n,m}$ the set of partitions $\pi=\{\pi_1,\ldots,\pi_m\}$  of the set $\{1,\ldots,n\}$ with $m$ blocks $\pi_i$
of size $\vert \pi_i\vert$, with $i\in \{1,\ldots, m\}$. We also let $\Pi_{n}$ the set of partitions of the set $\{1,\ldots,n\}$
and $\flat(\pi)$ the number of blocks in a given partition $\pi$, and $\Pi^+_{n}$ the subset of partitions $\pi$ s.t. $\flat(\pi)>1$.

Let $[n]$ be the set of $m$ multiple indexes $i=(i_1,\ldots,i_n)\in \{1,\ldots,d\}^n$. 
For any given $i\in [n]$ and any subset $S=\{j_1,\ldots,j_s\}\subset \{1,\ldots,n\}$ we set
$$
i_{S}=(i_{j_1},\ldots,i_{j_s})
$$

For any $x=(x^1,\ldots,x^d)\in\RR^d$ and any multiple index $i\in [n]$ we write 
$
\partial_{i}
$ instead of $\partial_{x^{i_1},\ldots x^{i_n}}=\partial_{x^{i_1}}\ldots\partial_{ x^{i_n}}$ the $n$-th partial derivatives w.r.t. the coordinates $(x^{i_1},\ldots x^{i_n})$. 

Let $f$ and $X$ be a couple of smooth function from $\RR^d$ into itself.
In this notation for any $i\in [n]$ and $1\leq j\leq d$ we have the multivariate Fa\`a di Bruno derivation formula
$$
\partial_i(f^j\circ X)=\sum_{1\leq m\leq n}~\sum_{k\in [m]}~\partial_{k}f^j(X)~\sum_{\pi\in \Pi_{n,m}}(\nabla^{\pi}X)_{i,k}
$$
with the $\pi$-gradient tensor
$$
(\nabla^{\pi}X)_{i,k}:=(\nabla^{\vert\pi_1\vert}X)_{i_{\pi_1},k_1}\ldots (\nabla^{\vert\pi_m\vert}X)_{i_{\pi_m},k_m}
$$
We check the above formula by induction w.r.t. the parameter $n$. In a more compact we have checked the following lemma.
\begin{lem}
For any $n\geq 1$ we have the Fa\`a di Bruno derivation formula
\begin{equation}\label{faa-di-bruno-ap}
\nabla^n(f\circ X)=\sum_{\pi\in \Pi_{n}}(\nabla^{\pi}X)~(\nabla^{\flat(\pi)} f)(X)
\end{equation}
\end{lem}
Whenever $X(x)$ is a random function we have
\begin{equation}\label{faa-di-bruno-ap-2}
P(f)(x):=\EE((f\circ X)(x))\Longrightarrow
\nabla^nP(f)=\sum_{1\leq m\leq n}~\Pa^{[n,m]}(\nabla^m f)
\end{equation}
with the collection of integral operators
$$
\Pa^{[n,m]}(\nabla^m f)(x):=\sum_{\pi\in \Pi_{n,m}}~\EE\left((\nabla^{\pi}X(x))~\nabla^m f(X(x))\right)
$$

Using the above  lemma we also check the stochastic tensor evolution equation
$$
\begin{array}{l}
\partial_t\, (\nabla^n X^{\mu}_{s,t}(x))_{i,j}\\
\\
\displaystyle= (\nabla^n X^{\mu}_{s,t}(x)~b_t^{[1]}(X^{\mu}_{s,t}(x),\phi_{s,t}(\mu)))_{i,j}+\sum_{1<m\leq n}\sum_{k\in [m]}\sum_{\pi\in \Pi_{n,m}}(\nabla^{\pi}X^{\mu}_{s,t}(x))_{i,k}~
b_t^{[1]_{\flat(\pi)}}(x,\phi_{s,t}(\mu))_{k,j}:
\end{array}
$$
with
$$
b_t^{[1]_m}(x,\mu)_{(k_1,\ldots,k_m),j}:=\partial_{k_1,\ldots,k_m}b^j_t(x,\mu)
$$
In a more compact form we have
$$
\partial_t\, \nabla^n X^{\mu}_{s,t}(x)= \nabla^n X^{\mu}_{s,t}(x)~b_t^{[1]}(X^{\mu}_{s,t}(x),\phi_{s,t}(\mu))+\sum_{\pi\in \Pi^+_{n}}\nabla^{\pi}X^{\mu}_{s,t}(x)
~b_t^{[1]_{\flat(\pi)}}(x,\phi_{s,t}(\mu))
$$
This implies that
$$
\begin{array}{l}
   \displaystyle\partial_t\, \nabla^n X^{\mu}_{s,t}(x)\nabla^n X^{\mu}_{s,t}(x)^{\prime}\\
\\
   \displaystyle=\nabla^n X^{\mu}_{s,t}(x)~\left(b_t^{[1]}(X^{\mu}_{s,t}(x),\phi_{s,t}(\mu))+b_t^{[1]}(X^{\mu}_{s,t}(x),\phi_{s,t}(\mu))^{\prime}\right)\nabla^n X^{\mu}_{s,t}(x)^{\prime}\\
   \\
     \displaystyle +\sum_{\pi\in \Pi^+_{n}}\nabla^{\pi}X^{\mu}_{s,t}(x)
~\left(b_t^{[1]_{\flat(\pi)}}(x,\phi_{s,t}(\mu))+b_t^{[1]_{\flat(\pi)}}(x,\phi_{s,t}(\mu))^{\prime}\right)~\nabla^n X^{\mu}_{s,t}(x)^{\prime}
\end{array}
$$
Taking the trace in the above display, we check that
$$
\partial_t\Vert\nabla^n X^{\mu}_{s,t}(x)\Vert^2_{\tiny Frob}\leq~-2\lambda_1~\Vert\nabla^n X^{\mu}_{s,t}(x)\Vert_{\tiny Frob}^2+2\Vert\nabla^n X^{\mu}_{s,t}(x)\Vert_{\tiny Frob}
~
\sum_{\pi\in \Pi^+_{n}}\Vert b^{[1]_{\flat(\pi)}}\Vert_{\tiny Frob}~\Vert\nabla^{\pi} X^{\mu}_{s,t}(x)\Vert_{\tiny Frob}$$
This yields the rather crude estimate
$$
\begin{array}{l}
\partial_t\Vert\nabla^n X^{\mu}_{s,t}(x)\Vert^2_{\tiny Frob}\\
\\
    \displaystyle \leq~-2\lambda_1~\Vert\nabla^n X^{\mu}_{s,t}(x)\Vert_{\tiny Frob}^2+c_{n}~\Vert\nabla^n X^{\mu}_{s,t}(x)\Vert_{\tiny Frob}
~
\sum_{\pi\in \Pi^+_{n}}~\Vert\nabla^{\vert{\pi_1\vert}} X^{\mu}_{s,t}(x)\Vert_{\tiny Frob}\ldots \Vert\nabla^{\vert{\pi_{\flat(\pi)}\vert}} X^{\mu}_{s,t}(x)\Vert_{\tiny Frob}
\end{array}$$
from which we check that
$$
\begin{array}{l}
\partial_t\Vert\nabla^n X^{\mu}_{s,t}(x)\Vert_{\tiny Frob}\\
\\
    \displaystyle \leq~-\lambda_1~\Vert\nabla^n X^{\mu}_{s,t}(x)\Vert_{\tiny Frob}+c_{n}~
~
\sum_{}~\Vert\nabla X^{\mu}_{s,t}(x)\Vert^{l_1}_{\tiny Frob}~\Vert\nabla^2 X^{\mu}_{s,t}(x)\Vert^{l_2}_{\tiny Frob}\ldots \Vert\nabla^{n-1} X^{\mu}_{s,t}(x)\Vert_{\tiny Frob}^{l_{n-1}}
\end{array}$$
The summation in the above display is taken over all indices $l_1,\ldots,l_{n-1}$ such that
$l_1+\ldots+l_{n-1}=m$ and $l_1+2l_2+\ldots+(n-1)l_{n-1}=n$ and $1<m\leq n$. Assume that (\ref{tau-n-estimates}) has been checked up to rank $(n-1)$.
In this case, we have
$$
\Vert\nabla^n X^{\mu}_{s,t}(x)\Vert_{\tiny Frob}\leq c_{n,1}~e^{-\lambda_1(t-s)}~\int_s^t~e^{\lambda_1(u-s)}~e^{-2\lambda_1(u-s)}~du\leq c_{n,2}~e^{-\lambda_1(t-s)}
$$
This ends the proof of (\ref{tau-n-estimates}).\cqfd
\subsection*{Proof of (\ref{bismut-omega}) and (\ref{bismut-omega-2})}\label{bismut-omega-proof}
We recall the backward formula
$$
P^{\mu}_{s,t}(f)(x)=f(x)+\int_s^tL_{u,\phi_{s,u}(\mu)}\left(P^{\phi_{s,u}(\mu)}_{u,t}(f)\right)(x)~du
$$
A detailed proof of the above formula based on backward stochastic flows can be found in theorem 3.1 in the article~\cite{mp-var-18}.
This implies that
$$
d\left(P^{\phi_{s,u}(\mu)}_{u,t}(f)(X^{\mu}_{u,t}(x))\right)=\left(\nabla P^{\phi_{s,u}(\mu)}_{u,t}(f)\right)(X^{\mu}_{u,t}(x))^{\prime}~dW_u
$$
from which we check that
$$
f(X^{\mu}_{s,t}(x))=
P^{\mu}_{s,t}(f)(x)+\int_s^t\left(\nabla P^{\phi_{s,u}(\mu)}_{u,t}(f)\right)(X^{\mu}_{u,t}(x))^{\prime}~dW_u
$$
This yields the formula
$$
\begin{array}{l}
   \EE\left(f(X_{s,t}^{\mu}(x))~  \tau^{\mu,\omega}_{s,t}(x)\right)\\
   \\
   \displaystyle=\EE\left(\left(\int_s^t\left(\nabla P^{\phi_{s,u}(\mu)}_{u,t}(f)\right)(X^{\mu}_{u,t}(x))^{\prime}~dW_u\right)\left(\int_s^t~   \partial_u \omega_{s,t}(u)~\nabla X^{\mu}_{s,u}(x)~dW_u
\right)\right)\\
\\
   \displaystyle=\EE\left(\int_s^t~\nabla \left(P^{\phi_{s,u}(\mu)}_{u,t}(f)\circ X^{\mu}_{u,t}\right)(x)~ \partial_u \omega_{s,t}(u)~du\right)
\end{array}
$$
We conclude that
$$
\begin{array}{l}
   \EE\left(f(X_{s,t}^{\mu}(x))~  \tau^{\mu,\omega}_{s,t}(x)\right)\\
   \\
   \displaystyle=\nabla P^{\mu}_{s,t}(f)(x)~\EE\left(\int_s^t \partial_u \omega_{s,t}(u)~du\right)=\nabla P^{\mu}_{s,t}(f)(x)~(\omega_{s,t}(t)-\omega_{s,t}(s))=\nabla P^{\mu}_{s,t}(f)(x)
\end{array}
$$
This ends the proof of  (\ref{bismut-omega}).  For any $s\leq u\leq t$ applying (\ref{bismut-omega}) to the function $P^{\phi_{s,u}(\mu)}_{u,t}(f)$ we have
\begin{eqnarray*}
\nabla P^{\mu}_{s,t}(f)(x)&=&\nabla P^{\mu}_{s,u}\left(P^{\phi_{s,u}(\mu)}_{u,t}(f)\right)(x)\\
&=&\EE\left(P^{\phi_{s,u}(\mu)}_{u,t}(f)(X^{\mu}_{s,u}(x))~\int_s^u~ \partial_v \omega_{s,u}(v)~\nabla X^{\mu}_{s,v}(x)~dW_v
\right)
\end{eqnarray*}
This implies that
$$
\begin{array}{l}
\partial_{x_j,x_i} P^{\mu}_{s,t}(f)(x)\\
\\
=\displaystyle\sum_{1\leq l\leq d}\EE\left(\partial_{x_j}X^{\mu,l}_{s,u}(x)~\partial_{x_l}\left(P^{\phi_{s,u}(\mu)}_{u,t}(f)\right)(X^{\mu}_{s,u}(x))~

\int_s^u~ \partial_v \omega_{s,u}(v)~\partial_{x_i} X^{\mu,k}_{s,v}(x)~dW^k_v\right)\\
\\
\hskip3cm\displaystyle+\EE\left(P^{\phi_{s,u}(\mu)}_{u,t}(f)(X^{\mu}_{s,u}(x))~\int_s^u~ \partial_v \omega_{s,u}(v)~\partial_{x_j,x_i} X^{\mu,k}_{s,v}(x)~dW^k_v\right)\
\end{array}$$
Applying (\ref{bismut-omega}) to the first term we check that
$$
\begin{array}{l}
\displaystyle\EE\left(\partial_{x_j}X^{\mu,l}_{s,u}(x)~\partial_{x_l}\left(P^{\phi_{s,u}(\mu)}_{u,t}(f)\right)(X^{\mu}_{s,u}(x))~
\tau^{\mu,\omega}_{s,u}(x)_i\right)\\
\\
=\displaystyle
\displaystyle\EE\left(f(X^{\mu}_{s,t}(x))~~\partial_{x_j}X^{\mu,l}_{s,u}(x)~
\tau^{\mu,\omega}_{s,u}(x)_i
\left(\sum_{1\leq m\leq d}
\int_u^t~\partial_v \omega_{u,t}(v)~\left(\partial_{x_l}X_{u,v}^{\phi_{s,u}(\mu),m}\right)(X_{s,u}^{\mu}(x))~dW^m_v
\right)
\right)\\
\\
=\displaystyle
\EE\left(f(X^{\mu}_{s,t}(x))~~\partial_{x_j}X^{\mu,l}_{s,u}(x)~
~\tau^{\phi_{s,u}(\mu),\omega}_{u,t}(X_{s,u}^{\mu}(x))_l~\tau^{\mu,\omega}_{s,u}(x)_i\right)
\end{array}$$
We conclude that
$$
\begin{array}{l}
\nabla^2P^{\mu}_{s,t}(f)(x)_{i,j}=\nabla^2P^{\mu}_{s,t}(f)(x)_{j,i}\\
\\
=\displaystyle\EE\left(f(X^{\mu}_{s,t}(x))~~\nabla X^{\mu}_{s,u}(x)_{j,l}~
~\tau^{\phi_{s,u}(\mu),\omega}_{u,t}(X_{s,u}^{\mu}(x))_l~\tau^{\mu,\omega}_{s,u}(x)_i\right)\\
\\
\hskip3cm\displaystyle+\EE\left(P^{\phi_{s,u}(\mu)}_{u,t}(f)(X^{\mu}_{s,u}(x))~\int_s^u~ \partial_v \omega_{s,u}(v)~\nabla^2 X^{\mu}_{s,v}(x)_{(i,j),k}~dW^k_v\right)\
\end{array}$$
This ends the proof of (\ref{bismut-omega-2}).\cqfd

\subsection*{Proof of (\ref{linear-example-d})}\label{linear-example-d-proof}
We have
 $$
 \begin{array}{l}
\displaystyle 
p_{s,t}^{\mu_1,\mu_0}(x,z)=b^{\mu_0}_{s,t}(x,z) +\sum_{n\geq 1}
\int_{[s,t]_n}~B_2~e^{(t-u_n)B_1} B_2~e^{(u_n-u_{n-1})B_1}\ldots~B_2~e^{(u_2-u_1)B_1}~\\
\\
\displaystyle\times\left(B_1~\phi_{s,u_1}(\mu_1)(e)+
B_2\left(e^{(u_1-s)B_1}(x-\mu_0(e))+e^{(u_1-s)(B_1+B_2)}\mu_0(e)\right)
\right)~du_1\ldots du_n~\end{array}
$$
Recalling that
\begin{eqnarray*}
\phi_{s,u_1}(\mu_1)(e)&=&
e^{(u_1-s)[B_1+B_2]}~\mu_1(e)\\
b_{s,t}^{\mu_0}(x,z)&=&B_1z+B_2\left[e^{(t-s)B_1}(x-\mu_0(e))+e^{(t-s)[B_1+B_2]}~\mu_0(e)\right]
\end{eqnarray*}
and using the rather well known exponential formulae
\begin{eqnarray*}
e^{(t-s)(B_1+B_2)}&=&e^{(t-s)B_1}+\int_s^t~e^{(t-u)B_1}\,B_2\,e^{(u-s)(B_1+B_2)}~du\\
&=&e^{(t-s)B_1}+\int_s^t~e^{(t-u)(B_1+B_2)}\,B_2\,e^{(u-s)B_1}~du
\end{eqnarray*}
we check that
 $$
 \begin{array}{l}
\displaystyle 
p_{s,t}^{\mu_1,\mu_0}(x,z)=b^{\mu_0}_{s,t}(x,z) \\
\\
\displaystyle+B_2~\int_s^t~e^{(t-u_1)(B_1+B_2)}~\left(B_1~\phi_{s,u_1}(\mu_1)(e)+
B_2\left(e^{(u_1-s)B_1}(x-\mu_0(e))+e^{(u_1-s)(B_1+B_2)}\mu_0(e)\right)
\right)~du_1~\end{array}
$$
from which we find that
 $$
 \begin{array}{l}
\displaystyle 
p_{s,t}^{\mu_1,\mu_0}(x,z)=B_1z+B_2\left[e^{(t-s)B_1}(x-\mu_0(e))+e^{(t-s)[B_1+B_2]}~\mu_0(e)\right]\\
\\
\displaystyle +B_2~\left[\int_s^t e^{(t-u_1)(B_1+B_2)}~B_1~e^{(u_1-s)[B_1+B_2]}~du_1\right]
~\mu_1(e)~~\\
\\
\displaystyle +B_2~\left[\int_s^t
e^{(t-u_1)(B_1+B_2)}~B_2~e^{(u_1-s)(B_1+B_2)}~du_1\right]~
\mu_0(e) +B_2~\left[e^{(t-s)(B_1+B_2)}-e^{(t-s)B_1}\right]
~(x-\mu_0(e))
\end{array}
$$
This ends the proof of  (\ref{linear-example-d}).\cqfd
\subsection*{Proof of (\ref{nabla-D-times-2})}\label{nabla-D-times-2-poof}
We have the tensor product formula
$$
\left(D_{\mu_1,\mu_0}\phi_{s,t}\right)^{\otimes 2}:=\left(P^{\mu_0}_{s,t}\right)^{\otimes 2}+\left(Q^{\mu_1,\mu_0}_{s,t}\right)^{\otimes 2}+Q^{\mu_1,\mu_0}_{s,t}\otimes P^{\mu_0}_{s,t}+P^{\mu_0}_{s,t}\otimes Q^{\mu_1,\mu_0}_{s,t}
$$

We also have
$$
\begin{array}{l}
\displaystyle
\left(Q^{\mu_1,\mu_0}_{s,t}\otimes P^{\mu_0}_{s,t}\right)(g)(x,\overline{x})\\
\\
\displaystyle=\int_{\Delta_{s,t}}~\Phi_{s,u}(\mu_1)(d(u,y))~(b^{\mu_0}_{s,u}(x,y)^{\prime}\otimes I)~\left(\Pa^{\,\phi_{s,u}(\mu_0)}_{u,t}\otimes P^{\mu_0}_{s,t}\right)(\nabla_{x_1}g)(y,\overline{x})\end{array}$$
Recall that $b_{t}(x,y)$ is differentiable at any order with uniformly bounded derivatives. Thus all differentials of the above function
w.r.t. the coordinate $x$ have uniformly bounded derivatives. On the other hand, the mapping $x\mapsto b_{t}(x,y)$ has at most linear growth.
Thus, using the estimates (\ref{ref-moments}) and (\ref{tau-n-estimates}), for any $m\geq 0$ we check that
$$
\Vert Q^{\mu_1,\mu_0}_{s,t}\otimes P^{\mu_0}_{s,t}\Vert_{\Ca^{n+1}_m(\RR^{2d})\rightarrow \Ca^n_{m+1}(\RR^{2d})}\leq c_{m,n}(t)~\rho_{m\vee 2}(\mu_0,\mu_1)$$
In the same vein, we have the tensor product formula
$$
\begin{array}{l}
\displaystyle
\left(Q^{\mu_1,\mu_0}_{s,t}\right)^{\otimes 2}(g)(x,\overline{x})=\left(\Qa^{\mu_1,\mu_0}_{s,t}\right)^{\otimes 2}((\nabla\otimes\nabla) g)(x,\overline{x})\\
\\
\displaystyle:=\int_{\Delta_{s,t}\times \Delta_{s,t}}~\left[\Phi_{s,u}(\mu_1)\otimes \Phi_{s,u}(\mu_1)\right](d((u,y),(\overline{u},\overline{y})))~\\
\displaystyle\hskip7cm~\widehat{b}^{\,\mu_0}_{s,u,\overline{u}}((x,\overline{x}),(y,\overline{y}))^{\prime}~\widehat{\Pa}^{\,\phi_{s,u,\overline{u}}(\mu_0)}_{u,\overline{u},t}((\nabla\otimes\nabla) g)(y,\overline{y})
\end{array}
$$
with 
$$
\widehat{b}^{\,\mu_0}_{s,u,\overline{u}}((x,\overline{x}),(y,\overline{y}))^{\prime}:=b^{\,\mu_0}_{s,u}(x,y)^{\prime}\otimes  b^{\,\mu_0}_{s,\overline{u}}(\overline{x},\overline{y})^{\prime}
\quad\mbox{\rm and}\quad
\widehat{\Pa}^{\,\phi_{s,u,\overline{u}}(\mu_0)}_{u,\overline{u},t}:=\Pa^{\,\phi_{s,u}(\mu_0)}_{u,t}\otimes
\Pa^{\,\phi_{s,\overline{u}}(\mu_0)}_{\overline{u},t}
$$

Arguing as above and using the estimates (\ref{ref-moments}) and (\ref{tau-n-estimates}) for any $m\geq 0$ we check that
$$
\Vert \left(Q^{\mu_1,\mu_0}_{s,t}\right)^{\otimes 2}\Vert_{\Ca^{2}_m(\RR^{2d})\rightarrow \Ca^n_{2}(\RR^{2d})}\leq c_{m,n}(t)~\rho_{m\vee 2}(\mu_0,\mu_1)
$$

\subsection*{Proof of lemma~\ref{lem-tex}}\label{lem-tex-proof}

Using  the decomposition
$$
\begin{array}{l}
\phi_{s,u}(\mu_1)-\phi_{s,u}(\mu_0)\\
\\
\displaystyle=\sum_{1\leq l\leq n}\left[\phi_{s,u_1}(\mu_1)\otimes\ldots\otimes \phi_{s,u_{l-1}}(\mu_1) \right]\otimes \left[\phi_{s,u_l}(\mu_1)-\phi_{s,u_l}(\mu_0)\right]\otimes\left[\phi_{s,u_{l+1}}(\mu_{0})\otimes\ldots\otimes \phi_{s,u_{n}}(\mu_0) \right]
\end{array}
$$
which is valid for any $\mu_0,\mu_1\in P_2(\RR^d)$ and any    $u=(u_1,\ldots,u_n)\in [s,t]_n$ with $n\geq 1$, for any function
$$
(u,y)\in \Delta_{s,t}\mapsto h_u(y)\in \RR
$$
we check that
\begin{equation}\label{decomp-Phi-k}
\int_{\Delta_{s,t}}~\left[\Phi_{s,u}(\mu_1)-\Phi_{s,u}(\mu_0)\right](d(u,y))~h_u(y)=\int_{\Delta^1_{s,t}}~\left[\Phi_{s,v}(\mu_1)-\Phi_{s,v}(\mu_0)\right](d(v,z))~\overline{h}_{v}(z)
\end{equation}
with the function
$$  \begin{array}{l}
\displaystyle\overline{h}_{v}(z):=h_v(z)+\int_{\Delta_{s,v}}~\Phi_{s,u}(\mu_1)(d(u,y))~h_{u,v}(y,z)+\int_{\Delta_{v,t}}~\Phi_{v,u}\left(\phi_{s,v}(\mu_0)\right)(d(u,y))
~h_{v,u}(z,y)
\\
\\
    \hskip3cm \displaystyle+\int_{\Delta_{s,v}\times \Delta_{v,t}}~\Upsilon^{\mu_1,\mu_0}_{s,t}((v,z),d((u,y),(\overline{u},\overline{y})))~h_{(u,v,\overline{u})}(y,z,\overline{y})
     \end{array} $$
In the above display, $\Upsilon^{\mu_1,\mu_0}_{s,t}$ stands for the tensor product measures
$$
\Upsilon^{\mu_1,\mu_0}_{s,t}((v,z),d((u,y),(\overline{u},\overline{y})))=\Phi_{s,u}(\mu_1)(d(u,y))~\Phi_{v,\overline{u}}\left(\phi_{s,v}(\mu_0)\right)(d(\overline{u},\overline{y}))
$$
We also have the tensor product formula
$$
\begin{array}{l}
\displaystyle\left(D_{\mu_1,\mu_0}\phi_{s,t}\right)^{\otimes 2}-\left(D_{\mu_0}\phi_{s,t}\right)^{\otimes 2}\\
\\
\displaystyle=\left(Q^{\mu_1,\mu_0}_{s,t}\right)^{\otimes 2}-\left(Q^{\mu_0}_{s,t}\right)^{\otimes 2}+\left(Q^{\mu_1,\mu_0}_{s,t}-Q^{\mu_0}_{s,t}\right)\otimes P^{\mu_0}_{s,t}+P^{\mu_0}_{s,t}\otimes\left(Q^{\mu_1,\mu_0}_{s,t}-Q^{\mu_0}_{s,t}\right)
\end{array}$$
This yields the decomposition
$$
\begin{array}{l}
\displaystyle
\left(\left[Q^{\mu_1,\mu_0}_{s,t}-Q^{\mu_0}_{s,t}\right]\otimes P^{\mu_0}_{s,t}\right)(g)(x,\overline{x})\\
\\
\displaystyle:=\int_s^t~\int~\left[\phi_{s,v}(\mu_1)-\phi_{s,v}(\mu_0)\right](d\widehat{x})~\Ia_{s,v,t}^{\mu_0,\mu_1}(g)(x,\overline{x},\widehat{x})~dv\end{array}$$
with the integral operator
$$
\begin{array}{l}
\displaystyle\Ia_{s,v,t}^{\mu_0,\mu_1}(g)(x,\overline{x},\widehat{x})\\
\\
:=b^{\mu_0}_{s,v}(x,\widehat{x})^{\prime}~\left(\Pa^{\,\phi_{s,v}(\mu_0)}_{v,t}\otimes P^{\mu_0}_{s,t}\right)(\nabla_{x_1}g)(\widehat{x},\overline{x})\\
\\
\displaystyle+\int_{\Delta_{s,v}}~\Phi_{s,u}(\mu_1)(d(u,y))~b^{\mu_0}_{s,u,v}(x,y,\widehat{x})^{\prime}~\left(\Pa^{\,\phi_{s,v}(\mu_0)}_{v,t}\otimes P^{\mu_0}_{s,t}\right)(\nabla_{x_1}g)(\widehat{x},\overline{x})\\
\\
\displaystyle+\int_{\Delta_{v,t}}~\Phi_{v,u}\left(\phi_{s,v}(\mu_0)\right)(d(u,y))~b^{\mu_0}_{s,v,u}(x,\widehat{x},y)^{\prime}~\left(\Pa^{\,\phi_{s,u}(\mu_0)}_{u,t}\otimes P^{\mu_0}_{s,t}\right)(\nabla_{x_1}g)(y,\overline{x})\\
\\
+\displaystyle\int_{\Delta_{s,v}\times \Delta_{v,t}}~\Upsilon^{\mu_1,\mu_0}_{s,t}((v,z),d((u,y),(\overline{u},\overline{y})))~b^{\mu_0}_{s,u,v,\overline{u}}(x,y,\widehat{x},\overline{y})^{\prime}~\left(\Pa^{\,\phi_{s,\overline{u}}(\mu_0)}_{\overline{u},t}\otimes P^{\mu_0}_{s,t}\right)(\nabla_{x_1}g)(\overline{y},\overline{x})\end{array}$$
Arguing as in the proof of (\ref{crude-est-rev}) and (\ref{nabla-D}) we check that
$$
\Vert \Ia_{s,v,t}^{\mu_0,\mu_1}\Vert_{\Ca^{n+1}_{m}(\RR^{2d})\rightarrow \Ca^n_{m+2}(\RR^{3d})}\leq c_{m,n}(t)~\rho_{m\vee 2}(\mu_0,\mu_1)
$$

In the same vein, we have
$$
\begin{array}{l}
\displaystyle
\left[\left(Q^{\mu_1,\mu_0}_{s,t}\right)^{\otimes 2}-\left(Q^{\mu_0}_{s,t}\right)^{\otimes 2}\right](g)(x,\overline{x})\\\
\\
=\displaystyle\int_{\Delta_{s,t}}~\left[\Phi_{s,u}(\mu_1)-\Phi_{s,u}(\mu_0)\right](d(u,y))~\left[\Theta_{s,u,t}^{\mu_1,\mu_0}+\overline{\Theta}_{s,u,t}^{\mu_1,\mu_0}\right](g)(x,\overline{x},y)~dv
\end{array}
$$
with
$$
\overline{\Theta}_{s,u,t}^{\mu_1,\mu_0}(g)(x,\overline{x},y):=
\int_{\Delta_{s,t}}~\Phi_{s,\overline{u}}(\mu_1)(d(\overline{u},\overline{y}))~
\widehat{b}^{\,\mu_0}_{s,u,\overline{u}}((x,\overline{x}),(y,\overline{y}))^{\prime}~\widehat{\Pa}^{\,\phi_{s,u,\overline{u}}(\mu_0)}_{u,\overline{u},t}((\nabla\otimes\nabla) g)(y,\overline{y})
$$
and
$$
\Theta_{s,u,t}^{\mu_1,\mu_0}(g)(x,\overline{x},y):=
\int_{\Delta_{s,t}}~\Phi_{s,\overline{u}}(\mu_0)(d(\overline{u},\overline{y}))~
\widehat{b}^{\,\mu_0}_{s,\overline{u},u}((x,\overline{x}),(\overline{y},y))^{\prime}~\widehat{\Pa}^{\,\phi_{s,\overline{u},u}(\mu_0)}_{\overline{u},u,t}((\nabla\otimes\nabla) g)(\overline{y},y)
$$
This yields the formula
$$
\begin{array}{l}
\displaystyle
\left[\left(Q^{\mu_1,\mu_0}_{s,t}\right)^{\otimes 2}-\left(Q^{\mu_0}_{s,t}\right)^{\otimes 2}\right](g)(x,\overline{x})\\
\\
\displaystyle=\int_s^t~\left[\phi_{s,v}(\mu_1)-\phi_{s,v}(\mu_0)\right](d\widehat{x})~\Ja_{s,v,t}^{\mu_0,\mu_1}(g)(x,\overline{x},\widehat{x})~dv
\end{array}
$$
with the integral operator
$$
\begin{array}{l}
\displaystyle\Ja_{s,v,t}^{\mu_0,\mu_1}(g)(x,\overline{x},\widehat{x})\\
\\
\displaystyle:=\left[\Theta_{s,v,t}^{\mu_1,\mu_0}+\overline{\Theta}_{s,v,t}^{\mu_1,\mu_0}\right](g)(x,\overline{x},\widehat{x})\\
\\
\displaystyle+\int_{\Delta_{s,v}}~\Phi_{s,u}(\mu_1)(d(u,y))~\left[\Theta_{s,(u,v),t}^{\mu_1,\mu_0}+\overline{\Theta}_{s,(u,v),t}^{\mu_1,\mu_0}\right](g)(x,\overline{x},(y,\widehat{x}))\\
\\
\displaystyle+\int_{\Delta_{v,t}}~\Phi_{v,u}\left(\phi_{s,v}(\mu_0)\right)(d(u,y))~\left[\Theta_{s,(v,u),t}^{\mu_1,\mu_0}+\overline{\Theta}_{s,(v,u),t}^{\mu_1,\mu_0}\right](g)(x,\overline{x},(\widehat{x},y))\\
\\
\displaystyle+\int_{\Delta_{s,v}\times \Delta_{v,t}}~\Upsilon^{\mu_1,\mu_0}_{s,t}((v,z),d((u,y),(\overline{u},\overline{y})))~\left[\Theta_{s,(u,v,\overline{u}),t}^{\mu_1,\mu_0}+\overline{\Theta}_{s,(u,v,\overline{u}),t}^{\mu_1,\mu_0}\right](g)(x,\overline{x},(y,\widehat{x},\overline{y}))
\end{array}$$
Arguing as above, we check that
$$
\Vert \Ja_{s,v,t}^{\mu_0,\mu_1}\Vert_{\Ca^{2}_m(\RR^{2d})\rightarrow \Ca^n_{2}(\RR^{3d})}\leq c_{m,n}(t)~\rho_{m\vee 2}(\mu_0,\mu_1)$$
Combining the above decompositions we find that 
$$
\begin{array}{l}
\left[\left(D_{\mu_1,\mu_0}\phi_{s,t}\right)^{\otimes 2}-\left(D_{\mu_0}\phi_{s,t}\right)^{\otimes 2}\right](g)(x,\overline{x})\\
\\
\displaystyle=\int_s^t~\left[\phi_{s,v}(\mu_1)-\phi_{s,v}(\mu_0)\right](d\widehat{x})~\Ka_{s,v,t}^{\mu_0,\mu_1}(g)(x,\overline{x},\widehat{x})~dv\quad\mbox{\rm with}\quad\Ka_{s,v,t}^{\mu_0,\mu_1}:=
2~\Ia_{s,v,t}^{\mu_0,\mu_1}+\Ja_{s,v,t}^{\mu_0,\mu_1}
\end{array}$$
For any $n\geq 2$ and $m\geq 0$ we have
$$
\Vert \Ka_{s,v,t}^{\mu_0,\mu_1}\Vert_{\Ca^{n+1}_{m}(\RR^{2d})\rightarrow \Ca^n_{m+2}(\RR^{3d})}\leq c_{m,n}(t)~\rho_{m\vee 2}(\mu_0,\mu_1)$$
We conclude that
$$
\left(\phi_{s,t}(\mu_1)-\phi_{s,t}(\mu_0)\right)^{\otimes 2}=(\mu_1-\mu_0)^{\otimes 2}\left(D_{\mu_0}\phi_{s,t}\right)^{\otimes 2}+(\mu_1-\mu_0)^{\otimes 3}\Ra_{\mu_1,\mu_0}\phi_{s,t}
$$
with the operator
$$
\begin{array}{l}
\Ra_{\mu_1,\mu_0}\phi_{s,t}(g)(x,\overline{x},\widehat{x})\\
\\
\displaystyle :=\int_s^t
\left[\int~P^{\mu_0}_{s,v}(\widehat{x},dz)~\Ka_{s,v,t}^{\mu_0,\mu_1}(g)(x,\overline{x},z)+\int_{\Delta_{s,v}}~\Phi_{s,u}(\mu_1)(d(u,y))~b^{\mu_0}_{s,u}(\widehat{x},y)^{\prime}~\La_{s,u,v,t}^{\mu_0,\mu_1}(g)(x,\overline{x},y)\right]~dv
\end{array}$$
In the above display, $\La_{s,u,v,t}^{\mu_0,\mu_1}$ stands for the integral operator operator
\begin{eqnarray*}
\La_{s,u,v,t}^{\mu_0,\mu_1}(g)(x,\overline{x},y)
&=&\Pa^{\,\phi_{s,u}(\mu_0)}_{u,t}\left(\nabla_{x_3}\Ka_{s,v,t}^{\mu_0,\mu_1}(g)(x,\overline{x},\point)\right)(y)
\end{eqnarray*}
We also check that
$$
\Vert \Ra_{\mu_1,\mu_0}\phi_{s,t}\Vert_{\Ca^{n+2}_m(\RR^{2d})\rightarrow \Ca^n_{m+3}(\RR^{3d})}\leq c_{m,n}(t)~\rho_{m+2}(\mu_1,\mu_2)
$$
This ends the proof of the lemma.
\cqfd

\subsection*{Proof of the estimate (\ref{D3-estimation})}\label{D3-estimation-lem-proof}

For any $x=(x_1,x_2)\in \RR^{2d}$ we set $\sigma(x_1,x_2):=\sigma(x_2,x_1)$. In this notation, for any matrix valued function $h(x)=(h_{i,j}(x))_{1\leq i,j\leq d}$ we have the tensor product formula
  $$
  \begin{array}{l}
     \displaystyle\left(\Da_{\mu_1,\mu_0}\phi_{s,t}\right)^{\otimes 2}( h)(x)\\
     \\
       \displaystyle = \left(  \Pa^{\mu_0}_{s,t}\right)^{\otimes 2}(h)(x)  +\int_{\Delta_{s,t}}\Phi_{s,v}(\mu_1)(d(u,y))~ \left[ \II^{\,\mu_0}_{s,u,t}(h)(x,y)+\II^{\,\mu_0}_{s,u,t}(h)(\sigma(x),y)\right]\\
      \\
     \displaystyle \hskip3cm+\int_{\Delta_{s,t}\times \Delta_{s,t}}\Phi_{s,u}(\mu_1)(d(u,y))~\Phi_{s,v}(\mu_1)(d(v,z))~  \JJ^{\,\mu_0}_{s,u,v,t}(h)(x,y,z)
  \end{array} $$
  with the matrix valued functions $ \II^{\,\mu_0}_{s,u,t}(h)$ and $  \JJ^{\,\mu_0}_{s,u,v,t}(h)$ given for any $(u,y)\in \Delta^n_{s,t}$ and $(v,z)\in \Delta^m_{s,t}$ by the formula
  \begin{eqnarray*}
  \II^{\,\mu_0}_{s,u,t}(h)(x,y)&:=&\BB^{\mu_0}_{s,u}(x_1,y)~\left(\Pa^{\phi_{s,u_n}(\mu_0)}_{u_n,t}\otimes\Pa^{\mu_0}_{s,t}\right)(h)(y_n,x_2)\\
  \JJ^{\,\mu_0}_{s,u,v,t}(h)(x,y,z)&:=&\left[\BB^{\mu_0}_{s,u}(x_1,y)\otimes \BB^{\mu_0}_{s,v}(x_2,z)\right]\left(\Pa^{\phi_{s,u_n}(\mu_0)}_{u_n,t}\otimes\Pa^{\phi_{s,v_m}(\mu_0)}_{v_m,t}\right)(h)(y_n,z_m)
  \end{eqnarray*}

Using (\ref{ref-nabla-Pa}) we have
$$
 \nabla\Pa^{\mu}_{s,t}(g)=\Pa^{[2,1],\mu}_{s,t}(g)+\Pa^{[2,2],\mu}_{s,t}(\nabla g)
$$
from which we check the formula
  $$
      \begin{array}{l}
  \nabla_{y_n} \left(\Pa^{\phi_{s,u_n}(\mu_0)}_{u_n,t}\otimes\Pa^{\phi_{s,v_m}(\mu_0)}_{v_m,t}\right)(h)(y_n,z_m)\\
  \\
    =\left[\Pa^{[2,1],\phi_{s,u_n}(\mu_0)}_{u_n,t}\otimes\Pa^{\phi_{s,v_n}(\mu_0)}_{v_n,t}\right](h)(y_n,z_m)+\left[\Pa^{[2,2]\phi_{s,u}(\mu_0)}_{u_n,t}\otimes\Pa^{\phi_{s,v_m}(\mu_0)}_{v_m,t}\right](\nabla_{x_1}h)(y_n,z_m)
 \end{array}   $$
 By symmetry arguments, we also have
  $$
      \begin{array}{l}
  \nabla_{z_m} \left(\Pa^{\phi_{s,u_n}(\mu_0)}_{u_n,t}\otimes\Pa^{\phi_{s,v_m}(\mu_0)}_{v_m,t}\right)(h)(y_n,z_m)\\
  \\
    =\left[\Pa^{[2,1],\phi_{s,v_m}(\mu_0)}_{v_m,t}\otimes\Pa^{\phi_{s,u_n}(\mu_0)}_{u_n,t}\right](h)(z_m,y_n)+\left[\Pa^{[2,2],\phi_{s,v_m}(\mu_0)}_{v_m,t}\otimes\Pa^{\phi_{s,u_n}(\mu_0)}_{u_n,t}\right](\nabla_{x_1}h)(z_m,y_n)
 \end{array}   $$
 Using (\ref{estimates-U-V}) for any differentiable matrix valued function $h(x_1,x_2)$ such that $\Vert h\Vert\vee \Vert \nabla_{x_1}h\Vert\leq 1$ 
  we have the uniform estimate
 $$
 \Vert   \nabla_{y_n} \left(\Pa^{\phi_{s,u_n}(\mu_0)}_{u_n,t}\otimes\Pa^{\phi_{s,v_m}(\mu_0)}_{v_m,t}\right)(h)(y_n,z_m)\Vert\leq  c_1~e^{-\lambda_1[(t-u_n)+(t-v_m)]}
 $$
 In the same vein, we have
 $$
 \begin{array}{l}
\BB^{\mu_0}_{s,u}(x,y)= \BB_{s,u}^{[1],\mu_0}(x,y)\\
\\
\displaystyle=\EE\left[\nabla X_{s,u_1}^{\mu_0}(x)~b^{[2]}_{u_1}(y_1,X_{s,u_1}^{\mu_0}(x))\right]\prod_{1\leq l<n}\EE\left[\nabla X_{u_l,u_{l+1}}^{\phi_{s,u_l}(\mu_0)}(y_l)~b^{[2]}_{u_{l+1}}(y_{l+1},X_{u_l,u_{l+1}}^{\phi_{s,u_l}(\mu_0)}(y_l))\right]
\end{array}
 $$
Using the gradient and the Hessian estimates (\ref{tau-1-estimates}) and (\ref{tau-2-estimates}) for any $1\leq k\leq n$ we check that
 $$
  \Vert   \nabla_{y_k}\BB^{\mu_0}_{s,u}(x_1,y)\Vert\leq  c_2~\Vert b^{[2]}\Vert_2^n~e^{-\lambda_1(u_n-s)}~
 $$
Combining the above estimates with (\ref{estimate-BB}) we check that
 $$
       \begin{array}{l}
 \Vert   \nabla_{y_n}   \II^{\,\mu_0}_{s,u,t}(h)(x,y)\Vert\\
 \\
 \displaystyle\leq  c_{3}~\Vert b^{[2]}\Vert_2^n~\left[
 e^{-\lambda_1(u_n-s)}~e^{-\lambda_1[(t-u_n)+(t-s)]}+e^{-\lambda_1[(u_n-s)]}~e^{-\lambda_1[(t-u_n)+(t-s)]}\right]\\
 \\
  \displaystyle\leq  c_{4}~\Vert b^{[2]}\Vert_2^n~ e^{-2\lambda_1(t-s)}~ \end{array}  $$
In addition, for any $1\leq k<n$ we have
  $$
         \begin{array}{l}
 \Vert   \nabla_{y_k}   \II^{\,\mu_0}_{s,u,t}(h)(x,y)\Vert\\
 \\
  \displaystyle \leq   c_{5}~\Vert b^{[2]}\Vert_2^n~e^{-\lambda_1(u_n-s)}~e^{-\lambda_1[(t-u_n)+(t-s)]}  \leq c_{5}~\Vert  b^{[2]}\Vert_2^n~e^{-2\lambda_1(t-s)}
 \end{array}   $$
 We conclude that 
   \begin{equation}\label{estimate-II}
 \sup_{1\leq k\leq n}\Vert   \nabla_{y_k}   \II^{\,\mu_0}_{s,u,t}(h)(x,y)\Vert \leq c~\Vert b^{[2]}\Vert_2^n~e^{-2\lambda_1(t-s)}
 \end{equation}
 Arguing as above, for any $1\leq k<n$ we have
  $$
         \begin{array}{l}
 \Vert   \nabla_{y_k}    \JJ^{\,\mu_0}_{s,u,v,t}(h)(x,y,z)\Vert\\
 \\
  \displaystyle \leq   c_{1}~ \Vert b^{[2]}\Vert_2^{m+n}~~e^{-\lambda_1(u_n-s)}~e^{-\lambda_1(v_m-s)}~e^{-\lambda_1[(t-u_n)+(t-v_m)]}  \leq   c_{2}~ \Vert b^{[2]}\Vert_2^{m+n}~~e^{-2\lambda_1(t-s)} \end{array}   $$
 In addition, for $k=n$ we have
   $$
         \begin{array}{l}
 \Vert   \nabla_{y_n}    \JJ^{\,\mu_0}_{s,u,v,t}(h)(x,y,z)\Vert\\
 \\
  \displaystyle \leq   c_{3}~ \Vert\nabla_{x_2}b\Vert_2^{m+n}~\left[
  e^{-\lambda_1(u_n-s)}~e^{-\lambda_1(v_m-s)}~e^{-\lambda_1[(t-u_n)+(t-v_m)]}\right.\\
  \\
  \hskip3cm\left.+e^{-\lambda_1(u_n-s)}~e^{-\lambda_1(v_m-s)}~e^{-\lambda_1[(t-u_n)+(t-v_m)]}
\right]
 \end{array}   $$
 This implies that
   \begin{equation}\label{estimate-JJ}
\sup_{1\leq k\leq n} \Vert   \nabla_{y_k}    \JJ^{\,\mu_0}_{s,u,v,t}(h)(x,y,z)\Vert \leq    c~ \Vert b^{[2]}\Vert_2^{m+n}~~e^{-2\lambda_1(t-s)} 
 \end{equation}
On the other hand, we have the decomposition
    $$
     \displaystyle\left[\left(\Da_{\mu_1,\mu_0}\phi_{s,t}\right)^{\otimes 2}-\left(\Da_{\mu_0}\phi_{s,t}\right)^{\otimes 2}\right]( h)(x) =\int_{\Delta_{s,t}}\left[\Phi_{s,v}(\mu_1)-\Phi_{s,v}(\mu_0)\right](d(u,y))~ \KK^{\,\mu_0,\mu_1}_{s,u,t}(h)(x,y) $$
  with the matrix valued function
  $$
  \begin{array}{l}
     \displaystyle\KK^{\,\mu_0,\mu_1}_{s,u,t}(h)(x,y):= \II^{\,\mu_0}_{s,u,t}(h)(x,y)+ \II^{\,\mu_0}_{s,u,t}(h)(\sigma(x),y) +\int_{\Delta_{s,t}}~\Phi_{s,v}(\mu_1)(d(v,z))
  ~ \JJ^{\,\mu_0}_{s,u,v,t}(h)(x,y,z)\\
   \\
      \hskip7cm  \displaystyle+\int_{\Delta_{s,t}}\Phi_{s,v}(\mu_0)(d(v,z))
  ~ \JJ^{\,\mu_0}_{s,v,u,t}(h)(x,z,y)
  \end{array} $$
Using the estimates (\ref{estimate-II}) and (\ref{estimate-JJ}),  for any $(u,y)\in \Delta^n_{s,t}$ we check that
 \begin{equation}\label{estimate-KK}
  \begin{array}{l}
     \displaystyle
\sup_{1\leq k\leq n} \Vert  \nabla_{y_k} \KK^{\,\mu_0,\mu_1}_{s,u,t}(h)(x,y)\Vert\\
\\
 \leq  
c_{1}~\Vert b^{[2]}\Vert_2^n~e^{-\lambda_1(t-s)}~\left[e^{-\lambda_1(t-s)}+\left(e^{\Vert  b^{[2]}\Vert_2(t-s)}-1\right)
e^{-\lambda_1(t-s)}\right]\\
\\
 \leq  
c_{2}~\Vert  b^{[2]}\Vert_2^n~e^{-\lambda_1(t-s)}~e^{- \lambda_{1,2}(t-s)}  \end{array}
 \end{equation}
Using  the decomposition (\ref{decomp-Phi-k})
we also check that
    $$
  \begin{array}{l}
     \displaystyle\left[\left(\Da_{\mu_1,\mu_0}\phi_{s,t}\right)^{\otimes 2}-\left(\Da_{\mu_0}\phi_{s,t}\right)^{\otimes 2}\right]( h)(x)  =\int_{s}^t~\left[\phi_{s,v}(\mu_1)-\phi_{s,v}(\mu_0)\right](dz)~  \overline{\KK}^{\,\mu_0,\mu_1}_{s,v,t}(h)(x,z)  ~dv
       \end{array} $$
       with the matrix valued function
       $$   
         \begin{array}{l}
    \displaystyle \overline{\KK}^{\,\mu_0,\mu_1}_{s,v,t}(h)(x_1,x_2,x_3)\\
    \\
  \displaystyle={\KK}^{\,\mu_0,\mu_1}_{s,v,t}(h)(x_1,x_2,x_3)  + \int_{\Delta_{s,v}}~\Phi_{s,u}(\mu_1)(d(u,y))~
     \KK^{\,\mu_0,\mu_1}_{s,u,v,t}(h)(x_1,x_2,(y,x_3)) \\
\\    
   \hskip4cm     \displaystyle   +\int_{\Delta_{v,t}}~\Phi_{v,u}\left(\phi_{s,v}(\mu_0)\right)(d(u,y))
~ \KK^{\,\mu_0,\mu_1}_{s,v,u,t}(h)(x_1,x_2,x_3,y) 
 \\
    \\
  \hskip4cm     \displaystyle+\int_{\Delta_{s,v}\times \Delta_{v,t}}~\Upsilon^{\mu_1,\mu_0}_{s,t}((v,z),d((u,y),(\overline{u},\overline{y})))~   \KK^{\,\mu_0,\mu_1}_{s,u,v,\overline{u},t}(h)(x_1,x_2,(y,x_3,\overline{y}))       \end{array} $$
     Using  (\ref{estimate-KK}) we find the uniform estimates
    \begin{equation}\label{estimate-over-KK}
      \begin{array}{l}
\Vert  \nabla_{x_3}\overline{\KK}^{\,\mu_0,\mu_1}_{s,v,t}(h)(x_1,x_2,x_3)\Vert\\
\\
\displaystyle \leq   
c_{1}~\left[
e^{-2 \lambda_{1,2}(t-s)} +\left(e^{\Vert b^{[2]}\Vert_2 (t-s)}-1\right)
~e^{-\lambda_1(t-s)}~e^{- \lambda_{1,2}(t-s)} \right]\leq  c_{2}~e^{-2 \lambda_{1,2}(t-s)} 
     \end{array}
     \end{equation}
     On the other hand, using (\ref{def-D}) and (\ref{a-1}) we have
    $$
      \left[\phi_{s,t}(\mu_1)- \phi_{s,t}(\mu_0)\right](f) =   (\mu_1-\mu_0)P^{\mu_0}_{s,t}(f)+(\mu_1-\mu_0)\Qa^{\mu_1,\mu_0}_{s,t}(\nabla f)
     $$   
         Thus, recalling that
         $$
         \Qa^{\mu_1,\mu_0}_{s,t}(\nabla f)(z):=\int_{\Delta_{s,t}}~\Phi_{s,u}(\mu_1)(d(u,y))~b^{\mu_0}_{s,u}(z,y)^{\prime}~\Pa^{\,\phi_{s,u}(\mu_0)}_{u,t}(\nabla f)(y)
         $$
          we check that
       $$
  \begin{array}{l}
     \displaystyle\left[\left(\Da_{\mu_1,\mu_0}\phi_{s,t}\right)^{\otimes 2}-\left(\Da_{\mu_0}\phi_{s,t}\right)^{\otimes 2}\right]( h)(x)  =\int~(\mu_1-\mu_0)(dz)~\int_{s}^t~P^{\mu_0}_{s,v}\left(\overline{\KK}^{\,\mu_0,\mu_1}_{s,v,t}(h)(x,\point)\right)(z)  ~dv\\
        \\
        \hskip.3cm \displaystyle+
        \int~(\mu_1-\mu_0)(dz)~\int_{s}^t~        
        \int_{\Delta_{s,v}}~\Phi_{s,u}(\mu_1)(d(u,y))~b^{\mu_0}_{s,u}(z,y)^{\prime}~\Pa^{\,\phi_{s,u}(\mu_0)}_{u,v}(\nabla_{x_3}\overline{\KK}^{\,\mu_0,\mu_1}_{s,v,t}(h)(x,\point))(y)~dv
       \end{array} $$
       This implies that
       $$
         \begin{array}{l}
(\nabla\otimes\nabla)D^2_{\mu_0,\mu_1}\phi_{s,t}(f)(x_1,x_2)-
(\nabla\otimes\nabla)D^2_{\mu_0}\phi_{s,t}(f)(x_1,x_2)\\
\\
  \displaystyle= \int~(\mu_1-\mu_0)(dx_3)~\int_s^t~ \LL^{\mu_1,\mu_0}_{s,u}\left(
\SSS^{[2,1],\phi_{s,u}(\mu_0)}_{u,t}(\nabla f)+\SSS^{[2,2],\phi_{s,u}(\mu_0)}_{u,t}(\nabla^2 f)\right) (x_1,x_2,x_3)
~du
     \end{array} 
$$
       with the tensor integral operator
          $$
  \begin{array}{l}
 \LL^{\mu_1,\mu_0}_{s,t}(h) (x_1,x_2,x_3)
     \displaystyle:=~\int_{s}^t~P^{\mu_0}_{s,v}\left(\overline{\KK}^{\,\mu_0,\mu_1}_{s,v,t}(h)(x_1,x_2,\point)\right)(x_3)  ~dv\\
        \\
        \hskip.3cm \displaystyle+
 ~\int_{s}^t~        
        \int_{\Delta_{s,v}}~\Phi_{s,u}(\mu_1)(d(u,y))~b^{\mu_0}_{s,u}(x_3,y)^{\prime}~\Pa^{\,\phi_{s,u}(\mu_0)}_{u,v}(\nabla_{x_3}\overline{\KK}^{\,\mu_0,\mu_1}_{s,v,t}(h)(x_1,x_2,\point))(y)~dv
       \end{array} $$  
    On the other hand, using (\ref{nabla-equiv-2})
    $$
\begin{array}{l}
(\mu_1-\mu_0)^{\otimes 2}D^2_{\mu_1,\mu_0}\phi_{s,t}(f)-(\mu_1-\mu_0)^{\otimes 2}D^2_{\mu_0}\phi_{s,t}(f)\\
\\
\displaystyle=\int_{[0,1]^3}~\int_s^t~\EE\left(\langle \nabla_{x_3}\LL^{\mu_1,\mu_0}_{s,u}\left(
\SSS^{[2,1],\phi_{s,u}(\mu_0)}_{u,t}(\nabla f)+\SSS^{[2,2],\phi_{s,u}(\mu_0)}_{u,t}(\nabla^2 f)\right)(\Ya_{\epsilon}),(\Ya_1-\Ya_0)^{\otimes 3}\rangle\right)~du~d\epsilon\end{array}  
$$ 
   with the interpolating path
$$
\epsilon=(\epsilon_1,\epsilon_2,\epsilon_3)\mapsto
\Ya_{\epsilon}:=\left(\overline{Y}^1_0+\epsilon_1(\overline{Y}^1_1-\overline{Y}^1_0),
\overline{Y}^2_0+\epsilon_2(\overline{Y}^2_1-\overline{Y}^2_0),\overline{Y}^3_0+\epsilon_3(\overline{Y}^3_1-\overline{Y}^3_0)\right)
$$
and
$$
(\Ya_1-\Ya_0)^{\otimes 3}:=(\overline{Y}^1_1-\overline{Y}^1_0)\otimes (\overline{Y}^2_1-\overline{Y}^2_0)\otimes (\overline{Y}^3_1-\overline{Y}^2_0)
$$
In the above display, $ (\overline{Y}^i_1,\overline{Y}^i_0)_{i=1,2,3}$ stands for independent copies of a pair of random variables $(Y_{0},Y_1)$ with distribution  $(\mu_{0},\mu_1)$.    
       
Using the commutation formula (\ref{commutation-tensor})  we check that
      $$
  \begin{array}{l}
\nabla_{x_3} \LL^{\mu_1,\mu_0}_{s,t}(h) (x_1,x_2,x_3)
     \displaystyle:=~\int_{s}^t~\Pa^{\mu_0}_{s,v}\left(\nabla_{x_3}\overline{\KK}^{\,\mu_0,\mu_1}_{s,v,t}(h)(x_1,x_2,\point)\right)(x_3)  ~dv\\
        \\
        \hskip.3cm \displaystyle+
 ~\int_{s}^t~        
        \int_{\Delta_{s,v}}~\Phi_{s,u}(\mu_1)(d(u,y))~\BB^{\mu_0}_{s,u}(x_3,y)~\Pa^{\,\phi_{s,u}(\mu_0)}_{u,v}(\nabla_{x_3}\overline{\KK}^{\,\mu_0,\mu_1}_{s,v,t}(h)(x_1,x_2,\point))(y)~dv
       \end{array} $$  
       Using (\ref{estimate-over-KK}) for any differentiable matrix valued function $h(x_1,x_2)$ such that $\Vert h\Vert\vee \Vert \nabla_{x_1}h\Vert\leq 1$ and for any $\epsilon\in ]0,1[$ we check that
       $$
         \begin{array}{l}
      \Vert \nabla_{x_3} \LL^{\mu_1,\mu_0}_{s,t}(h) (x_1,x_2,x_3)\Vert\\
       \\
   \displaystyle    \leq c_{1}~ e^{-2 \lambda_{1,2}(t-s)}~\left[\int_s^te^{-\lambda_1(v-s)}~dv+    ~\int_{s}^t~        
       \left(e^{\Vert b^{[2]}\Vert_2(v-s)}-1\right)e^{-\lambda_1(v-s)}~dv\right]    \leq  c_{2}~ e^{-2 \lambda_{1,2}(t-s)}~
          \end{array} $$
          On the other hand, we have
$$
         \begin{array}{l}
       \nabla_{x_1}\left[\SSS^{[2,1],\mu}_{s,t}(\nabla f)+    \nabla_{x_1}\SSS^{[2,2],\mu}_{s,t}(\nabla^2 f)\right](x_1,x_2)\\
       \\
   \displaystyle=   b^{[1,1,2]}_s(x_1,x_2)~\nabla D_{\mu}\phi_{s,t}(f)(x_1)+b^{[2,2,1]}_s(x_2,x_1)~\nabla D_{\mu}\phi_{s,t}(f)(x_2)\\
   \\
   +\nabla^3D_{\mu}\phi_{s,t}(f)(x_1)~ b^{[2]}_s(x_1,x_2)^{\prime}+b^{[2,2]}_s(x_2,x_1)~\nabla^2D_{\mu}\phi_{s,t}(f)(x_2)+\nabla^2D_{\mu}\phi_{s,t}(f)(x_1)\star b^{[1,2]}_s(x_1,x_2)
       \end{array} $$    
       with the $\star$-tensor product
       $$
           \begin{array}{l}
       \left[\nabla^2D_{\mu}\phi_{s,t}(f)(x_1)\star b^{[1,2]}_s(x_1,x_2)\right]_{k,i,j}\\
       \\
       \displaystyle   =\sum_{1\leq l\leq d}\left[
       \nabla^2D_{\mu}\phi_{s,t}(f)(x_1)_{k,l}~ b^{[1,2]}_s(x_1,x_2)^{\prime}_{l,i,j}+ b^{[1,2]}_s(x_1,x_2)^{\prime}_{k,j,l}  \nabla^2D_{\mu}\phi_{s,t}(f)(x_1)^{\prime}_{l,j}\right]
        \end{array}
       $$    
Using (\ref{def-lambda-hat-3}) we check that
$$
\Vert   \nabla_{x_1}\left[\SSS^{[2,1],\mu}_{s,t}(\nabla f)+    \nabla_{x_1}\SSS^{[2,2],\mu}_{s,t}(\nabla^2 f)\right]\Vert\leq c\,e^{- \lambda(t-s)}~\sup_{k=1,2,3}\Vert\nabla^k f\Vert\quad \mbox{for some $\lambda>0$}
$$

We conclude that for any  function  $f\in \Ca^3(\RR^d)$ s.t. $\sup_{k=1,2,3}\Vert\nabla^k f\Vert\leq 1$
    $$
\vert (\mu_1-\mu_0)^{\otimes 2}D^2_{\mu_1,\mu_0}\phi_{s,t}(f)-(\mu_1-\mu_0)^{\otimes 2}D^2_{\mu_0}\phi_{s,t}(f)\vert\leq   c~ e^{-\lambda(t-s)}~\WW_2(\mu_0,\mu_1)^3\quad \mbox{for some $\lambda>0$}
$$ 
The last assertion comes from the formula
$$
\frac{1}{2}~(\mu_1-\mu_0)^{\otimes 2}D^2_{\mu_1,\mu_0}\phi_{s,t}=\frac{1}{2}
(\mu_1-\mu_0)^{\otimes 2}D^2_{\mu_0}\phi_{s,t}+(\mu_1-\mu_0)^{\otimes 3}D^3_{\mu_0,\mu_1}\phi_{s,t}
$$
\cqfd
\subsection*{Proof of theorem~\ref{theo-ae-taylor}}\label{theo-ae-taylor-proof}

We extend the operators $D^k_{\mu_1,\mu_0}\phi_{s,t}$ introduced in theorem~\ref{theo-intro-4} to tensor functions $f=(f_{i})_{i\in [n]}$ by considering the tensor function
with entries
\begin{equation}\label{def-D-k-extension}
D_{\mu_1,\mu_0}^k\phi_{s,t}(f)_i=D_{\mu_1,\mu_0}^k\phi_{s,t}(f_{i})
\end{equation}
By theorem~\ref{theo-intro-4} we have
\begin{equation}\label{k-deco}
\begin{array}{l}
\displaystyle\left[\phi_{s,u}(\mu_1)-\phi_{s,u}(\mu_0)\right](b_u(X^{\mu_0}_{s,u}(x),\point))\\
\\
\displaystyle=\int~(\mu_1-\mu_0)(dy)~d_{s,u}^{[1],\mu_1,\mu_0}(X^{\mu_0}_{s,u}(x),y)\\
\\
\displaystyle=\int~(\mu_1-\mu_0)(dy)~d_{s,u}^{[1],\mu_0}(X^{\mu_0}_{s,u}(x),y)+\frac{1}{2}~\int~(\mu_1-\mu_0)^{\otimes 2}(dz)~d_{s,u}^{[2],\mu_1,\mu_0}(X^{\mu_0}_{s,u}(x),z)\\
\\
\displaystyle=\int~(\mu_1-\mu_0)(dy)~d_{s,u}^{[1],\mu_0}(X^{\mu_0}_{s,u}(x),y)+\frac{1}{2}~\int~(\mu_1-\mu_0)^{\otimes 2}(dz)~d_{s,u}^{[2],\mu_0}(X^{\mu_0}_{s,u}(x),z)\\
\\
\displaystyle\hskip3cm+\int~(\mu_1-\mu_0)^{\otimes 3}(dz)~d_{s,u}^{[3],\mu_1,\mu_0}(X^{\mu_0}_{s,u}(x),z)
\end{array}
\end{equation}
with the functions
\begin{eqnarray*}
d_{s,t}^{[1],\mu_1,\mu_0}(X^{\mu_0}_{s,t}(x),y)&:=&D_{\mu_1,\mu_0}\phi_{s,t}(b_t(X^{\mu_0}_{s,t}(x),\point))(y)\\
d_{s,t}^{[2],\mu_1,\mu_0}(X^{\mu_0}_{s,t}(x),(z_1,z_2))&:=&D_{\mu_1,\mu_0}^2\phi_{s,t}(b_t(X^{\mu_0}_{s,t}(x),\point))(z_1,z_2)\\
d_{s,u}^{[3],\mu_1,\mu_0}(X^{\mu_0}_{s,u}(x),(z_1,z_2,z_3))&:=&D_{\mu_1,\mu_0}^3\phi_{s,t}(b_t(X^{\mu_0}_{s,t}(x),\point))(z_1,z_2,z_3)
\end{eqnarray*}
We also write $d_{s,t}^{[k],\mu}$ instead of $d_{s,t}^{[k],\mu,\mu}$.
Using  (\ref{estimate-D2-nabla}) and (\ref{commutation-D-etimate}) we check that
$$
\Vert \nabla_y\, d_{s,t}^{[1],\mu_1,\mu_0}(X^{\mu_0}_{s,t}(x),y)\Vert \leq c_1~~e^{-  \lambda(t-s)}
$$
as well as
\begin{equation}\label{R-ok-4-3}
\Vert  (\nabla_{z_1}\otimes\nabla_{z_2})\,d_{s,t}^{[2],\mu_1,\mu_0}(X^{\mu_0}_{s,t}(x),z_1,z_2)\Vert\leq c_2~e^{-\lambda(t-s)}\quad \mbox{for some $\lambda>0$}
\end{equation}
Using (\ref{D3-estimation}) we also have
\begin{equation}\label{R-ok-3}
\vert\int~(\mu_1-\mu_0)^{\otimes 3}(dz)~d_{s,t}^{[3],\mu_1,\mu_0}(X^{\mu_0}_{s,t}(x),z)\vert
\leq c_3~ e^{-\lambda(t-s)}~\WW_2(\mu_0,\mu_1)^3\quad \mbox{for some $\lambda>0$}
\end{equation}

On the other hand, we have the second order expansions
$$
\begin{array}{l}
\left[\nabla X^{\phi_{s,u}(\mu_0)}_{u,t}\right]({X}^{\mu_1}_{s,u}(x))^{\prime}-\left[\nabla X^{\phi_{s,u}(\mu_0)}_{u,t}\right]({X}^{\mu_0}_{s,u}(x))^{\prime}\\
\\
\displaystyle=\int_0^1
\left[\nabla^2 X^{\phi_{s,u}(\mu_0)}_{u,t}\right]\left(X^{\mu_0}_{s,u}(x)+\epsilon(X^{\mu_1}_{s,u}(y)-X^{\mu_0}_{s,u}(x))\right)^{\prime}~
[X^{\mu_1}_{s,u}(x)-X^{\mu_0}_{s,u}(x)]~d\epsilon\\
\\
=\left[\nabla^2 X^{\phi_{s,u}(\mu_0)}_{u,t}\right]({X}^{\mu_0}_{s,u}(x))^{\prime}~[X^{\mu_1}_{s,u}(x)-X^{\mu_0}_{s,u}(x)]\\
\\
\displaystyle+\int_0^1
(1-\epsilon)~\left[\nabla^3 X^{\phi_{s,u}(\mu_0)}_{u,t}\right]\left(X^{\mu_0}_{s,u}(x)+\epsilon(X^{\mu_1}_{s,u}(y)-X^{\mu_0}_{s,u}(x))\right)^{\prime}~
[X^{\mu_1}_{s,u}(x)-X^{\mu_0}_{s,u}(x)]^{\otimes 2}~d\epsilon
\end{array}$$
In the same vein, we have
$$
\begin{array}{l}
b_u(X^{\mu_1}_{s,u}(x),y)-b_u(X^{\mu_0}_{s,u}(x),y)\\
\\
\displaystyle=\int_0^1
~b_u^{[1]}\left(X^{\mu_0}_{s,u}(x)+\epsilon(X^{\mu_1}_{s,u}(x)-X^{\mu_0}_{s,u}(x)),y\right)^{\prime}~
[X^{\mu_1}_{s,u}(x)-X^{\mu_0}_{s,u}(x)]~d\epsilon\\
\\
\displaystyle=b_u^{[1]}(X^{\mu_0}_{s,u}(x),y)^{\prime}~[X^{\mu_1}_{s,u}(x)-X^{\mu_0}_{s,u}(x)]\\
\\
\displaystyle+\int_0^1
(1-\epsilon)~b_u^{[1,1]}\left(X^{\mu_0}_{s,u}(x)+\epsilon(X^{\mu_1}_{s,u}(x)-X^{\mu_0}_{s,u}(x)),y\right)^{\prime}~
[X^{\mu_1}_{s,u}(x)-X^{\mu_0}_{s,u}(x)]^{\otimes 2}~d\epsilon
\end{array}$$
This implies that
$$
\begin{array}{l}
 \displaystyle X^{\mu_1}_{s,t}(x)-X^{\mu_0}_{s,t}(x)\\
 \\
  \displaystyle=\int_s^t\left[\nabla X^{\phi_{s,u}(\mu_0)}_{u,t}\right]({X}^{\mu_0}_{s,u}(x))^{\prime}~\left[\phi_{s,u}(\mu_1)-\phi_{s,u}(\mu_0)\right](b_u(X^{\mu_0}_{s,u}(x),\point))~du+\sum_{k=2,3}R^{[k],\mu_0,\mu_1}_{s,t}(x)
\end{array}$$
with
 the second order remainder term
 $$
\begin{array}{l}
\displaystyle R^{[2],\mu_0,\mu_1}_{s,t}(x)\\
\\
\displaystyle:= \int_s^t\left[\nabla^2 X^{\phi_{s,u}(\mu_0)}_{u,t}\right]({X}^{\mu_0}_{s,u}(x))^{\prime}~[X^{\mu_1}_{s,u}(x)-X^{\mu_0}_{s,u}(x)]~\left[\phi_{s,u}(\mu_1)-\phi_{s,u}(\mu_0)\right](b_u(X^{\mu_0}_{s,u}(x),\point))~du\\
\\
\displaystyle+\int_s^t\left[\nabla X^{\phi_{s,u}(\mu_0)}_{u,t}\right]({X}^{\mu_0}_{s,u}(x))^{\prime}~
\left[\phi_{s,u}(\mu_1)-\phi_{s,u}(\mu_0)\right]\left(b_u^{[1]}(X^{\mu_0}_{s,u}(x),\point)^{\prime}\right)~[X^{\mu_1}_{s,u}(x)-X^{\mu_0}_{s,u}(x)]~du
\end{array}$$
and  the third order remainder term
$$
\begin{array}{l}
\displaystyle R^{[3],\mu_0,\mu_1}_{s,t}(x)\\
\\
\displaystyle:=\int_0^1\int_s^t\left[\nabla^2 X^{\phi_{s,u}(\mu_0)}_{u,t}\right]({X}^{\mu_0}_{s,u}(x))^{\prime}~[X^{\mu_1}_{s,u}(x)-X^{\mu_0}_{s,u}(x)]~\\
\\
\displaystyle\hskip.3cm\left[\phi_{s,u}(\mu_1)-\phi_{s,u}(\mu_0)\right]\left(b_u^{[1]}\left(X^{\mu_0}_{s,u}(x)+\epsilon(X^{\mu_1}_{s,u}(x)-X^{\mu_1}_{s,u}(x)),\point\right)^{\prime}~\right)
[X^{\mu_1}_{s,u}(x)-X^{\mu_0}_{s,u}(x)]~d\epsilon~du\\\
\\
\displaystyle+\int_0^1
(1-\epsilon)~\int_s^t\left[\nabla X^{\phi_{s,u}(\mu_0)}_{u,t}\right]({X}^{\mu_0}_{s,u}(x))^{\prime}~\\
\\
\displaystyle\hskip.3cm\left[\phi_{s,u}(\mu_1)-\phi_{s,u}(\mu_0)\right]
\left(b_u^{[1,1]}\left(X^{\mu_0}_{s,u}(x)+\epsilon(X^{\mu_1}_{s,u}(x)-X^{\mu_1}_{s,u}(x)),\point\right)^{\prime}~
\right)~[X^{\mu_1}_{s,u}(x)-X^{\mu_0}_{s,u}(x)]^{\otimes 2}~du~d\epsilon\\
\\
\displaystyle+\int_0^1~(1-\epsilon)~ \int_s^t
\left[\nabla^3 X^{\phi_{s,u}(\mu_0)}_{u,t}\right]\left(X^{\mu_0}_{s,u}(x)+\epsilon(X^{\mu_1}_{s,u}(y)-X^{\mu_1}_{s,u}(x))\right)^{\prime}~\\
\\
\displaystyle\hskip3cm [X^{\mu_1}_{s,u}(x)-X^{\mu_0}_{s,u}(x)]^{\otimes 2}
~\left[\phi_{s,u}(\mu_1)-\phi_{s,u}(\mu_0)\right](b_u(X^{\mu_1}_{s,u}(x),\point))~du~~d\epsilon\\
\\
\end{array}$$
Combining  (\ref{tau-n-estimates}) with (\ref{ref-stab-W-1}) and (\ref{estimate-X-mu}) for any $k=1,2$ we check the uniform estimate
\begin{equation}\label{R-k-e}
\Vert R^{[k],\mu_0,\mu_1}_{s,t}(x)\Vert\leq c~e^{-\lambda(t-s)}~\WW_2(\mu_0,\mu_1)^k\quad \mbox{for some $\lambda>0$}
\end{equation}
We check (\ref{a-e-T1})  using (\ref{R-ok-4-3}) and (\ref{k-deco}).

Using (\ref{def-lambda-hat-3}) we also have the estimate
$$
\Vert \nabla_y\,D_{\mu_0}X^{\mu_0}_{s,t}(x,y)\Vert\leq c_3~e^{-\lambda(t-s)}\quad \mbox{for some $\lambda>0$}
$$
Observe that
$$
\begin{array}{l}
\displaystyle\left[\phi_{s,u}(\mu_1)-\phi_{s,u}(\mu_0)\right](b_u^{[1]}(X^{\mu_0}_{s,u}(x),\point)^{\prime})=\int~(\mu_1-\mu_0)(dy)~d_{s,u}^{[1,1],\mu_1,\mu_0}(X^{\mu_0}_{s,u}(x),y)\\
\\
\displaystyle=\int~(\mu_1-\mu_0)(dy)~d_{s,u}^{[1,1],\mu_0}(X^{\mu_0}_{s,u}(x),y)+\frac{1}{2}~\int~(\mu_1-\mu_0)^{\otimes 2}(dz)~d_{s,u}^{[2,1],\mu_1,\mu_0}(X^{\mu_0}_{s,u}(x),z)
\end{array}$$
with the matrix valued functions
\begin{eqnarray*}
d_{s,t}^{[1,1],\mu_1,\mu_0}(X^{\mu_0}_{s,t}(x),y)&:=&D_{\mu_1,\mu_0}\phi_{s,t}(b_t^{[1]}(X^{\mu_0}_{s,t}(x),\point)^{\prime})(y)\\
d_{s,t}^{[2,1],\mu_1,\mu_0}(X^{\mu_0}_{s,t}(x),z_1,z_2)&:=&D_{\mu_1,\mu_0}^2\phi_{s,t}(b_t^{[1]}(X^{\mu_0}_{s,t}(x),\point)^{\prime})(z_1,z_2)
\end{eqnarray*}
We also write $d_{s,t}^{[1,1],\mu}$ instead of $d_{s,t}^{[1,1],\mu,\mu}$.
Observe that
 $$
\begin{array}{l}
\displaystyle R^{[2],\mu_0,\mu_1}_{s,t}(x)\\
\\
\displaystyle=\frac{1}{2}~\int~(\mu_1-\mu_0)^{\otimes 2}(dz)~\int_s^t\left[\nabla^2 X^{\phi_{s,u}(\mu_0)}_{u,t}\right]({X}^{\mu_0}_{s,u}(x))^{\prime}~D^{[2,1]}_{\mu_0}X^{\mu_0}_{s,u}(x,z)~du\\
\\
\displaystyle+\frac{1}{2}~\int~(\mu_1-\mu_0)^{\otimes 2}(dy)~\int_s^t\left[\nabla X^{\phi_{s,u}(\mu_0)}_{u,t}\right]({X}^{\mu_0}_{s,u}(x))^{\prime}~D^{[1,1]}_{\mu_0}X^{\mu_0}_{s,u}(x,z)~du+R^{[3,2],\mu_0,\mu_1}_{s,t}(x)
\end{array}$$
with
 $$
\begin{array}{l}
\displaystyle R^{[3,2],\mu_0,\mu_1}_{s,t}(x)\\
\\
\displaystyle=\frac{1}{2}~\int~(\mu_1-\mu_0)^{\otimes 3}(dy)~\int_s^t\left[\nabla^2 X^{\phi_{s,u}(\mu_0)}_{u,t}\right]({X}^{\mu_0}_{s,u}(x))^{\prime}~
D_{\mu_0}X^{\mu_0}_{s,u}(x,y_1)~d_{s,u}^{[2],\mu_1,\mu_0}(X^{\mu_0}_{s,u}(x),(y_2,y_3))~du\\
\\
\displaystyle+\frac{1}{2}~\int~(\mu_1-\mu_0)^{\otimes 3}(dy)~\int_s^t\left[\nabla X^{\phi_{s,u}(\mu_0)}_{u,t}\right]({X}^{\mu_0}_{s,u}(x))^{\prime}~
d_{s,u}^{[2,1],\mu_1,\mu_0}(X^{\mu_0}_{s,u}(x),(y_2,y_3))~
D_{\mu_0}X^{\mu_0}_{s,u}(x,y_1)~du\\
\\
\displaystyle+ \int_s^t\left[\nabla^2 X^{\phi_{s,u}(\mu_0)}_{u,t}\right]({X}^{\mu_0}_{s,u}(x))^{\prime}~
\Ra^{[2],\mu_0,\mu_1}_{s,u}(x)~\left[\phi_{s,u}(\mu_1)-\phi_{s,u}(\mu_0)\right](b_u(X^{\mu_0}_{s,u}(x),\point))~du\\
\\
\displaystyle+\int_s^t\left[\nabla X^{\phi_{s,u}(\mu_0)}_{u,t}\right]({X}^{\mu_0}_{s,u}(x))^{\prime}~
\left[\phi_{s,u}(\mu_1)-\phi_{s,u}(\mu_0)\right]\left(b_u^{[1]}(X^{\mu_0}_{s,u}(x),\point)^{\prime}\right)~\Ra^{[2],\mu_0,\mu_1}_{s,u}(x)~du
\end{array}$$
Observe that
\begin{equation}\label{R-k-32-e}
\Vert R^{[3,2],\mu_0,\mu_1}_{s,t}(x)\Vert\leq c~e^{-\lambda(t-s)}~\WW_2(\mu_0,\mu_1)^3
\quad \mbox{for some $\lambda>0$}
\end{equation}
This yields the second order decompositionn (\ref{ae-taylor})
with the remainder term
$$
\begin{array}{l}
 \displaystyle \Ra^{\mu_1,\mu_0}_{s,t}(x):=R^{[3],\mu_0,\mu_1}_{s,t}(x)+R^{[3,2],\mu_0,\mu_1}_{s,t}(x)\\
  \\
   \displaystyle \hskip3cm+\int~(\mu_1-\mu_0)^{\otimes 3}(dz)~\int_s^t\left[\nabla X^{\phi_{s,u}(\mu_0)}_{u,t}\right]({X}^{\mu_0}_{s,u}(x))^{\prime}~d_{s,u}^{[3],\mu_1,\mu_0}(X^{\mu_0}_{s,u}(x),z)~du
\end{array}$$
The end of the proof of is now a consequence of the estimates  (\ref{R-ok-3}), (\ref{R-k-e}) and (\ref{R-k-32-e}). The proof of the theorem is completed.
\cqfd

\end{document}